\documentclass[12pt,nopreprintline]{elsarticle}

\usepackage{framed,multirow}
\usepackage[mathlines]{lineno}

\usepackage{amssymb,amsthm}
\usepackage{latexsym}
\usepackage{bbm}
\usepackage{hyperref}
\usepackage{amsfonts}
\usepackage{graphicx}
\usepackage{subfig}
\usepackage{multicol}
\usepackage{epstopdf}
\usepackage{overpic}
\usepackage{amsmath,esint}
\usepackage{epsfig}
\usepackage{dsfont}
\usepackage{tikz}
\usepackage{cleveref}
\usetikzlibrary{arrows}
\usetikzlibrary{fadings}
\usepackage{pgfplots}
\pgfplotsset{compat=1.10}
\usepgfplotslibrary{fillbetween}
\usetikzlibrary{patterns}
\usepackage{appendix}
\usepackage[utf8]{inputenc}
\usepackage{enumitem}
\usepackage{mathrsfs}
\usepackage{verbatim}
\usepackage{algorithm}
\usepackage{tcolorbox}
\usepackage{algorithmicx}
\usepackage{algpseudocode}
\usepackage{mathtools}
\usepackage{pifont}
\usepackage{calc}
\usepackage{accents}
\usepackage{cancel}
\usepackage[margin=1in]{geometry}
\usepackage{ifpdf}
\usepackage{array}

\newcommand{\RN}[1]{%
\textup{\uppercase\expandafter{\romannumeral#1}}%
}

\setenumerate{label=(\roman*),itemsep=3pt,topsep=3pt}

\newtheorem{remark}{Remark}

 \newcommand{\bbE}{\mathbb E}

\newcommand{\bbR}{\mathbb R}

 \newcommand{\cJ}{\mathcal J}

 \newcommand{\cR}{\mathcal R}
 
 \newcommand{\cV}{\mathcal V}
 \newcommand{\cX}{\mathcal X}

\newcommand{\Rthree}{{\bbR^3}}

\begin{document}

\begin{frontmatter}
\title{Adjoint DSMC Method for Spatially Inhomogeneous Boltzmann Equation with General Boundary Conditions}

\author[address1]{Russel  Caflisch}
\ead{caflisch@courant.nyu.edu}
\address[address1]{Courant Institute of Mathematical Sciences, New York University, New York, NY 10012, USA}
\author[address2]{Linglai Chen}
\ead{linglai_chen@u.nus.edu}
\address[address2]{Department of Computer Science, School of Computing, National University of Singapore, Singapore 117417, Singapore}
\author[address3]{Yunan Yang\corref{mycorrespondingauthor}}
\ead{yunan.yang@cornell.edu}
\address[address3]{Department of Mathematics, Cornell University, Ithaca, NY 14853, USA}
\cortext[mycorrespondingauthor]{Corresponding author}

\begin{abstract}
We develop adjoint Direct Simulation Monte Carlo (DSMC) formulations for the spatially inhomogeneous Boltzmann equation with periodic, specular reflecting, diffuse thermal, and prescribed inflow boundary conditions. Periodic and specular boundaries are treated using a pathwise particle adjoint conditional on the realized event history. For diffuse thermal boundaries, we introduce a randomized-time regularization of wall-crossing events and use score-function terms to differentiate the resulting boundary probabilities. Reparameterization of the outgoing half-Maxwellian samples provides sensitivities with respect to wall temperatures and tangential wall velocities. Prescribed inflow requires a different construction because perturbations of incoming particles affect subsequent cell populations, local collision frequencies, and collision schedules. We therefore derive an ensemble adjoint based on the locally linearized Boltzmann collision operator and evaluate boundary sensitivities using local inflow-source scores, with injection counts held fixed. For a scalar objective, the dominant adjoint cost is largely independent of the number of parameters. Numerical experiments for Maxwell molecules validate the formulations against centered finite differences for thermal, mixed thermal-specular, two-sided inflow, and high-Mach Couette-flow configurations. The results demonstrate consistent gradient estimates, Monte Carlo convergence, stability with respect to the thermal regularization parameter, and accurate sensitivity calculation in a regime with limited relative statistical noise.
\end{abstract}


\end{frontmatter}

\section{Introduction}
\label{sec:Intro}

The Boltzmann equation is a fundamental kinetic model describing the statistical behavior of dilute gases through the combined effects of particle transport and binary collisions. Numerical methods for solving the Boltzmann equation play a central role in rarefied gas dynamics, plasma physics, radiative transfer, and related fields~\cite{cercignani1988boltzmann,cercignani2000rarefied}. Among these methods, Direct Simulation Monte Carlo (DSMC), originally introduced by Bird~\cite{bird1970direct,bird1994molecular}, has emerged as a particularly flexible and robust particle-based approach. DSMC approximates the Boltzmann equation by simulating stochastic particle trajectories that alternate between particle transport, probabilistic collision events, and boundary interactions; see~\cite{nanbu1980direct,babovsky1986simulation,pareschi2001introduction} for foundational analyses and algorithmic variants.

Adjoint methods are indispensable tools for sensitivity analysis, inverse problems, and partial differential equation (PDE)-constrained optimization~\cite{biegler2003large,hinze2008optimization}. In deterministic settings, adjoint formulations for kinetic equations and their discretizations are well understood and widely used. Extending these techniques to stochastic particle methods such as DSMC is more difficult because the dependence of a realized particle trajectory on a model or boundary parameter is generally only piecewise smooth. In a spatially inhomogeneous DSMC simulation, a parameter perturbation may change the spatial cell containing a particle, the collision partners selected in that cell, whether a particle reaches a boundary during a time step, or, under an inflow boundary condition, the set of particles present in the domain. Conditional on a fixed sequence of such events, the particle updates may be differentiated pathwise away from event-switching parameter values. The derivative of this fixed realization, however, does not necessarily represent the derivative of the expected DSMC objective when the perturbation changes the probabilities or distributions associated with the discrete events. For example, under a diffuse thermal boundary condition, a small perturbation may change whether a particle reaches the wall, after which its outgoing velocity is sampled from a half-Maxwellian flux. Similar difficulties occur in adjoint Monte Carlo methods for radiative transfer, neutron transport, and stochastic chemical kinetics~\cite{hoogenboom1977adjoint,bal2011importance,wang2016efficiency}.

Adjoint and gradient-based methods for kinetic equations have attracted growing interest in optimization, control, and inverse problems~\cite{mohamed2020monte,lovbak2022reversible}. For radiative transport and kinetic models, Monte Carlo gradient and adjoint formulations have clarified when pathwise differentiation is valid and when likelihood-ratio or score-function corrections are required to account for parameter-dependent sampling and discrete event selection~\cite{li2022monte,yang2023adjoint,ni2025ergodic}. In the context of the Boltzmann equation, adjoint DSMC formulations have been derived for spatially homogeneous problems with general collision models~\cite{caflisch2021adjoint,yang2023adjoint}, and broader perspectives on adjoint Monte Carlo methods for kinetic equations are provided in~\cite{caflisch2024adjoint}. Parallel developments in adjoint-based optimization for rarefied gas flows have employed deterministic kinetic solvers, including discrete-velocity and BGK-type methods, to enable topology and shape optimization across rarefied and continuum regimes~\cite{guan2024topology,yuan2024design,yuan2025adjoint}.

{Automatic differentiation (AD) computes derivatives of numerical programs by applying the chain rule, either in forward mode using dual numbers or in reverse mode through a computational graph~\cite{GriewankWalther2008,Baydin2018,Sapienza2024}. However, black-box AD applied to a realized DSMC trajectory treats discrete choices, such as cell assignments, collision-pair selections, boundary events, and parameter-dependent sampling laws, as fixed. It therefore does not generally recover the derivative of the expected DSMC objective when a parameter changes the probabilities of these events. This limitation motivates the specialized pathwise, score-function, and ensemble-adjoint formulations developed here.}

The present work develops complementary adjoint DSMC formulations for the spatially inhomogeneous Boltzmann equation with periodic, specular reflecting, diffuse thermal, and prescribed inflow boundary conditions. For periodic and specular boundaries, we employ a pathwise discrete adjoint that differentiates the realized particle updates while holding the sampled random variables and the associated discrete event history fixed. Away from spatial-cell and boundary interfaces, the hard cell-membership indicators are locally constant, and their derivatives therefore vanish almost surely along the realized trajectory. The resulting adjoint propagates sensitivities with respect to both the initial velocity distribution and the initial spatial distribution.

Diffuse thermal boundaries require additional treatment because a parameter perturbation may switch the boundary hit-or-miss event. We introduce a selective randomized-time regularization for particles that may interact with a thermal wall. For a fixed regularization width $\varepsilon$, the resulting adjoint computes the derivative of an $\varepsilon$-regularized objective $J_\varepsilon$. Score-function terms capture the sensitivity of the boundary-crossing probabilities, while reparameterizations of the half-Maxwellian flux yield derivatives with respect to wall temperatures and prescribed tangential wall velocities. The regularized scheme converges to the original DSMC update as $\varepsilon\to0$, and the numerical experiments examine the stability of both the objective and the adjoint gradient over the range of $\varepsilon$ used in the calculations.

Prescribed inflow boundaries require a different formulation. A perturbation of an incoming-particle distribution changes not only the injected particle states but also their subsequent cell memberships, the local particle densities, and the collision schedules. A pathwise adjoint that freezes one realized sequence of cell memberships and collision pairs therefore does not, in general, recover the derivative of the expected DSMC objective. We instead differentiate the ensemble evolution and obtain a scalar backward equation involving the adjoint of the locally linearized Boltzmann collision operator. This construction extends the continuous linearized-collision framework of our earlier work~\cite{caflisch2021adjoint} to a spatially inhomogeneous setting and couples it to free transport, absorbing outflow, and parameter-dependent boundary sources. The derivative of each inflow source is evaluated using a local likelihood-ratio score of the incoming flux distribution. A Monte Carlo discretization that retains a prescribed subset of particle states in each cell limits the storage required by the backward calculation. For a scalar objective, the dominant forward and backward costs are largely independent of the number of inflow parameters.

The proposed formulations are examined in four numerical examples. The first two examples consider two-sided thermal boundaries and mixed thermal-specular boundaries, respectively, and include sensitivities with respect to boundary temperatures and parameters in the initial velocity and spatial distributions. The third example validates the ensemble-adjoint formulation for a two-sided inflow problem and studies its dependence on the number of retained particle states. The fourth example considers a high-Mach-number planar Couette flow driven by moving thermal walls. It tests the wall-speed derivative in a regime where the sensitivity signal is clearly resolved relative to the statistical fluctuations. In all cases, the adjoint estimates are compared with centered finite-difference estimates using common random numbers.

The remainder of the paper is organized as follows: Section~\ref{sec:background} reviews background on the Boltzmann equation and the forward and adjoint DSMC algorithms. Section~\ref{sec:adjoint_DSMC_space} derives adjoint equations for DSMC with periodic and specular reflecting boundary conditions. Section~\ref{sec:adjoint_DSMC_thermal} develops the adjoint formulation for thermal boundary conditions using randomized time stepping. Section~\ref{sec:adjoint_DSMC_inflow} introduces the ensemble-adjoint formulation and the local boundary-source score for prescribed inflow conditions. Several numerical examples are presented in Section~\ref{sec:numerics}. We conclude in Section~\ref{sec:conclusion} with a discussion of implications for adjoint particle methods for transport-dominated PDEs.

\section{Background}\label{sec:background}
In this section, we briefly review the essential background for the spatially inhomogeneous Boltzmann equation with the bilinear collision operator, the DSMC method, and the adjoint DSMC method proposed in~\cite{caflisch2021adjoint,yang2023adjoint,caflisch2024adjoint}.

We adopt the argument order $f(x,v,t)$ for the distribution function, and $r(x,v)$ for observables. The projection $\mathcal P$ maps velocities to spatial dimensions (e.g., $\mathcal P=e_1^\top$ in 1D). Primes denote post-collision quantities, and $(v_{k,i},x_{k,i})$ denotes the $i$th particle state at time $t_k$. Adjoint variables associated with $(v_{k,i},x_{k,i})$ are $(\beta_{k,i},\alpha_{k,i})$.

We remark that the adjoint variables introduced throughout this work are mathematically equivalent to Lagrange multipliers in constrained optimization. In this setting, the adjoint variable enforces the kinetic evolution equation and enables efficient computation of gradients of objective functionals with respect to parameters. Thus, adjoint variables serve the same fundamental purpose as Lagrange multipliers by encoding sensitivity of the objective to constraints.

\subsection{Boltzmann Equation}
We consider the spatially inhomogeneous Boltzmann equation on a finite spatial domain $\Omega$,
\begin{equation}
\label{eq:Boltz1}
\frac{\partial f}{\partial t} + v \cdot \nabla_x f = Q(f,f),
\end{equation}
subject to the initial condition
\begin{equation}
\label{eq:Boltzmann_Initial}
f(x,v,0) = f_0(x,v),
\end{equation}
and appropriate boundary conditions (BCs) prescribed on $\partial \Omega$, which will be specified in \Cref{sec:adjoint_DSMC_space,sec:adjoint_DSMC_thermal,sec:adjoint_DSMC_inflow}. Here, $f(x,v,t)$ is a nonnegative density function describing the time evolution of particles located at position $x \in \Omega$ with velocity $v \in \mathbb{R}^3$ at time $t>0$.

The bilinear collision operator $Q(f,f)$ models binary particle collisions and is given by
\begin{equation}
\label{eq:nonlinearCollision}
Q(f,f)
=
\int_{\mathbb{R}^3}
\int_{\mathbb{S}^2}
q(v-v_1,\sigma)
\bigl(
f(v_1') f(v') - f(v_1) f(v)
\bigr)
\, d\sigma \, dv_1,
\end{equation}
where $(v', v_1')$ denote the post-collisional velocities corresponding to the pre-collisional velocities $(v, v_1)$, and the integration in $\sigma$ is taken over the unit sphere $\mathbb{S}^2$. Note that $Q(f,f)$ and  $f$ in~\eqref{eq:nonlinearCollision} are functions of the spatial variable $x$, velocity variable $v$ and time variable $t$.

In this work, we focus on the effects arising from spatial transport and BCs. Accordingly, for simplicity we restrict attention to Maxwell molecules, for which the collision kernel $q(v-v_1,\sigma)$ is constant. The proposed adjoint formulation, however, extends directly to more general collision kernels; see~\cite{yang2023adjoint}.

Because of conservation of the momentum $v+v_1$ and the energy $v^2 + v_1^2$, the elastic collision formulae satisfy
\begin{align}
v' &= 1/2 (v+v_1) + 1/2 |v-v_1| \sigma ,  \label{BoltzmannSolution0} \\
v_1'&=1/2 (v+v_1) - 1/2  |v-v_1| \sigma ,
\label{BoltzmannSolution1}
\end{align}
where $\sigma$ is a collision parameter representing a unit direction of the relative velocity of particles after collision. We will hereafter use the shorthand notation $f, f_1, f', f_1'$ to denote $f(v), f(v_1), f(v')$ and $f(v_1')$. Let
\begin{equation}\label{eq:alpha-sigma}
\zeta = (v -v_1){\hat~}, \quad \sigma = (v' -v_1'){\hat~},
\end{equation}
where the notation $x{\hat~} = x/|x|.$ The collision rules~\eqref{BoltzmannSolution0} and~\eqref{BoltzmannSolution1} can be presented in an operator formulation as
\begin{equation}\label{eq:AB_vel_General}
        \begin{bmatrix}
    v'\\
   v_1'
    \end{bmatrix}    = A(\sigma,\zeta) \begin{bmatrix}
    v\\
    v_1
    \end{bmatrix}, \quad
   \begin{bmatrix}
    v\\
     v_1
    \end{bmatrix}     = B(\sigma,\zeta)    \begin{bmatrix}
    v'\\
    v_1'
    \end{bmatrix},
\end{equation}
where
\begin{equation} \label{eq:AB_def}
A(\sigma ,\zeta )= \frac{1}{2}\begin{bmatrix}
    I+\sigma  \zeta ^\top & I-\sigma  \zeta ^\top\\
    I-\sigma  \zeta ^\top & I+\sigma  \zeta ^\top
    \end{bmatrix},\quad B(\sigma ,\zeta )= \frac{1}{2}\begin{bmatrix}
    I+\zeta  \sigma ^\top & I-\zeta  \sigma ^\top\\
    I-\zeta  \sigma ^\top & I+\zeta  \sigma ^\top
    \end{bmatrix},
\end{equation}
where $I$ is the identity matrix in $\Rthree$ and $B= A^\top = A^{-1}$. It was shown in~\cite{caflisch2021adjoint} that for Maxwell molecules, the first-order variations (with collision parameters fixed) in pre- and post-collision velocities satisfy
\begin{equation}\label{eq:fwd_perturb_General}
        \begin{bmatrix}
    \delta v'\\
    \delta v_1'
    \end{bmatrix}   =
    A(\sigma,\zeta)
 \begin{bmatrix}
    \delta  v\\
    \delta v_1
    \end{bmatrix}.
\end{equation}
Note that when the collision parameters (i.e., $\sigma$ and $\zeta$) are fixed, \eqref{eq:AB_vel_General} defines a linear operation, namely a matrix-vector product.

\subsection{The DSMC Method}\label{sec:intro_dsmc}

In this section, we describe the classical Direct Simulation Monte Carlo (DSMC) method~\cite{bird1970direct,nanbu1980direct,babovsky1986simulation,babovsky1989convergence} for the spatially inhomogeneous Boltzmann equation~\eqref{eq:Boltz1}, following the tutorial presentation in~\cite{pareschi2001introduction}. A standard approach for treating spatial inhomogeneity is operator splitting, which separates the collision dynamics from the advection (transport) dynamics. Specifically, one first solves the spatially homogeneous Boltzmann equation
\begin{equation}
\label{eq:homoBoltz}
\frac{\partial f}{\partial t} = Q(f,f)
\end{equation}
over a single time step, using the initial condition~\eqref{eq:Boltzmann_Initial}. The resulting solution is then used as the initial condition for the transport equation
\begin{equation}
\label{eq:transport}
\frac{\partial f}{\partial t} + v \cdot \nabla_x f = 0,
\end{equation}
which is likewise solved over one time step. The composition of these two substeps yields an approximation to the solution of the full Boltzmann equation~\eqref{eq:Boltz1} after one time step, and the procedure is iterated in time.

\begin{algorithm}
  \caption{Nanbu--Babovsky Algorithm for Maxwell Molecules\label{alg:DSMC}}
  {
\begin{algorithmic}[1]
\State Sample the initial particle velocities and positions, given the initial particle density $f_0(x,v)$, to get $\mathcal V_0 = \{v_{0,1},\dots,v_{0,N}\}$ and $\mathcal X_0 = \{x_{0,1},\dots,x_{0,N}\}$. Set $M=T/\Delta t$ for final time $T$. Divide the spatial domain $\Omega$ into $N_c$ spatial cells: $\Omega = \Omega_1 \cup \ldots  \cup  \Omega_{N_c}$.
\For{$k=1$ to $M$}
\State  Given the set of particle velocities $\mathcal V_{k-1}$ and positions $\mathcal X_{k-1}$ from time step $t_{k-1}$.
\For{$j=1$ to $N_c$}
\State Identify $\mathcal V_{k-1}^j \subseteq \mathcal V_{k-1}$ where the particle positions are in the $j$-th spatial cell $\Omega_j$.
\State Let $N^j_{k-1}=|\mathcal V_{k-1}^j|$ and define the local density and collision frequency by $\rho_{k-1}^j= N_{k-1}^j/(N|\Omega_j|)$ and $
\mu_{k-1}^j=\rho_{k-1}^j \int_{\mathbb S^2}q(\sigma)\,d\sigma$.
Select $\lceil
N_{k-1}^j\Delta t\,\mu_{k-1}^j/2 \rceil$ collision pairs uniformly without replacement from
$\mathcal V_{k-1}^j$.
\State For each selected pair $(v_{k-1,i},  v_{k-1,i_1})$, perform collision based on~\eqref{BoltzmannSolution0}-\eqref{BoltzmannSolution1}; obtain post-collision velocities $( v_{k-1,i}',  v_{k-1,i_1}')$. Set $v_{k-1,i}' = {v}_{k-1,i}$ for particles that do not collide.
\EndFor
\State Compute the new spatial position $x_{k-1,i}'$ based on~\eqref{eq:x_update}.
\If{$x_{k-1,i}'\in\Omega$}
\State Set $v_{k,i} = v_{k-1,i}'$ and $x_{k,i} = x_{k-1,i}'$.
\Else
\State Enforce the specific boundary conditions and obtain $v_{k,i}$ and $x_{k,i}$.
\EndIf
\EndFor
\end{algorithmic}
}
\end{algorithm}

The DSMC method for~\eqref{eq:Boltz1} is also based on such an operator splitting. We consider a collection of $N$ Monte Carlo particles evolving over the time interval $[0,T]$. The interval is divided into $M$ subintervals of equal length $\Delta t = T/M$. At the $k$th time level $t_k = k \Delta t$, the particle velocities are given by
\begin{equation}
\label{eq:vk}
\mathcal V_k = \{ v_1, \ldots, v_N \}(t_k),
\end{equation}
where we denote the velocity of the $i$th particle by $v_i(t_k)$ or equivalently $v_{k,i}$. Each particle is also associated with a spatial position
\begin{equation}
\label{eq:xk}
\mathcal X_k = \{ x_1, \ldots, x_N \}(t_k),
\end{equation}
with $x_{k,i}$ denoting the position of the $i$th particle at time $t_k$. We assume that all particles lie within the spatial domain $\Omega$ at the initial time $t=0$.

As in~\cite{caflisch2021adjoint}, we focus on a Maxwellian gas with a constant collision kernel, i.e., $ q(v-v_1,\sigma)=C $. The algorithm described in this subsection can, however, be extended to more general collision kernels; see, for example,~\cite{bird1970direct,pareschi2001introduction}. A summary of the DSMC algorithm is provided in~\Cref{alg:DSMC}. In contrast to the spatially homogeneous setting, each collision step is followed by an advection step.

According to~\eqref{eq:transport}, the post-collision spatial position of each particle is updated as
\begin{equation}
\label{eq:x_update}
x_{k,i}' = x_{k,i} + \Delta t\,\mathcal P v_{k,i}',
\end{equation}
where $\mathcal P$ denotes a projection operator mapping velocities in $\mathbb R^3$ to the spatial domain. For example, if $\Omega\subset\mathbb R$, then $\mathcal P=\mathbf e_1^\top=[1,0,0]$; for a two-dimensional spatial domain, $\mathcal P=[\mathbf e_1,\mathbf e_2]^\top$; and for a three-dimensional spatial domain, $\mathcal P$ is the identity matrix. Following this update, particle positions are required to remain within the spatial domain $\Omega$ through the enforcement of appropriate BCs.

While the spatial domain $\Omega$ may be one-, two-, or three-dimensional, the velocity space is always taken to be the full $\mathbb R^3$. Viewed in this way, \Cref{alg:DSMC} implements DSMC for the homogeneous Boltzmann equation~\eqref{eq:homoBoltz}, followed by the advection step~\eqref{eq:x_update} and the enforcement of BCs on the spatial domain.

Finally, we note that under inflow BCs, the total number of Monte Carlo particles $N$ may vary in time, whereas for the other BCs considered in this work, including periodic boundaries, specular reflection, and thermalization, the particle number remains constant.

\subsection{The Adjoint DSMC Method}\label{sec:adjoint_DSMC}
The adjoint DSMC method was first proposed in~\cite{caflisch2021adjoint} as an efficient particle-based method to compute the gradient for Boltzmann-equation constrained optimization problems for Maxwell molecules. It was later generalized to adapt to more general collision kernels, such as the variable hard sphere models~\cite{yang2023adjoint}. The idea of using adjoint Monte Carlo-type methods for kinetic equation-constrained optimization is more general and applies to the radiative transport equation~\cite{li2022monte}. We refer interested readers to~\cite{caflisch2024adjoint} for a survey.

For the Boltzmann equation with a constant collision kernel, the pathwise adjoint used for periodic, specular, and thermal boundary conditions can be interpreted as the formal adjoint of the forward DSMC method in~\Cref{alg:DSMC}, conditional on the random variables and discrete event history sampled during the forward simulation. Adjoint DSMC methods were previously developed for the spatially homogeneous Boltzmann equation~\eqref{eq:homoBoltz}. The present work extends this framework to spatial transport and several boundary conditions. The prescribed-inflow problem requires a different ensemble-adjoint formulation, developed in Section~\ref{sec:adjoint_DSMC_inflow}, because freezing the realized cell memberships and collision schedule does not generally capture the derivative of the expected inflow objective.

Consider a simple optimization problem for the spatially homogeneous Boltzmann equation~\eqref{eq:homoBoltz}. The initial condition is
\begin{equation}
f(v,0)=f_0(v;m), \label{BoltzmannIC}
\end{equation}
in which $f_0$ is the prescribed initial data depending on a parameter $m$. The goal is to find $m$ which optimizes the objective function at time $t=T$,
\begin{equation}\label{eq:OTD_obj}
J_1 (m) = \int_{\bbR^3} r(v) f(v,T) dv,
\end{equation}
where $f(v,T)$ is the solution to~\eqref{eq:homoBoltz} given the initial condition~\eqref{BoltzmannIC}, and thus $f(v,T)$ depends on $m$ through the initial condition. Here, we take parameters in the initial distribution as an illustrative example of the optimization problem, but the proposed approach readily extends to other types of parameters appearing in the Boltzmann equation; the same applies to the choice of objective function~\eqref{eq:OTD_obj}.

To compute the gradient of $J_1$ with respect to a parameter $m$,~\cite{caflisch2021adjoint}~considered the following total objective function $\cJ$ (i.e., Lagrangian) with Lagrangian multipliers $\{\beta_{k,i}\}$:
{
\begin{equation} \label{eq:DTO_obj}
\cJ = \underbracket{\frac{1}{N}\sum_{i=1}^N r(v_{{M},i})}_{\cJ_1} +\underbracket{ \frac{1}{N}\sum_{i=1}^N \beta_{0,i} \cdot \left(v_{0, i}- v_{0,i}(m) \right)}_{\cJ_2}+\underbracket{\frac{1}{N} \sum_{k=0}^{M-1} \sum_{i=1}^N \beta_{k,i} \cdot (v_{k+1,i}-v_{k,i}')}_{\cJ_3}.
\end{equation}
}
Here, $\cJ_1$ is the Monte Carlo quadrature of the objective function~\eqref{eq:OTD_obj} by particle velocities $\{v_{M,i}\}_{i=1}^N$ at the final time, $\cJ_2$ is the constraint on the DSMC initial condition $v=v_0(m)$ using the adjoint variable $ \{\beta_{0,i}\}_{i=1}^N$, and $\cJ_3$ is the constraint that enforces the binary collision law~\eqref{BoltzmannSolution0} and~\eqref{BoltzmannSolution1} using the adjoint variable $\beta_{k,i}$ for each particle $i$ at the $k$th time interval. In particular, $v_{k,i}'$ represents the post-collision velocity of particle $i$ if it participates in a collision at the $k$th time interval. Otherwise, $v_{k,i}' = v_{k,i}$, which means the particle velocity remains the same at the $(k+1)$th time interval. The advantage of the Lagrangian approach is that $\{v_{k,i}\}$ is a general set of velocities, and its dependence on the collision rules~\eqref{BoltzmannSolution0}-\eqref{BoltzmannSolution1} is imposed through the adjoint variables.

By setting the derivatives of $\cJ$ with respect to the adjoint variables $\{\beta_{k,i}\}$ as zero, we get the collision rule~\eqref{BoltzmannSolution0}-\eqref{BoltzmannSolution1} and the initial condition~\eqref{BoltzmannIC}. By setting the  derivatives of $\cJ$ with respect to the state variables $\{v_{k,i}\}$, we get the adjoint equations, which can be rewritten using the operator notation~\eqref{eq:AB_def}:
\begin{equation} \label{eq:adjoint_law}
   \begin{bmatrix}
    \beta_{k,i}\\
    \tilde \beta_{k,i_1}
    \end{bmatrix}     = C \begin{bmatrix}
    \beta_{k+1,i}\\
    \tilde \beta_{k+1,i_1}
    \end{bmatrix},\quad
C  = \begin{cases}
    B(\sigma_{k,i},\zeta_{k,i}) ,  & \text{if $(v_{k,i},v_{k,i_1})$ collided at $t_k$,}\\
    I , & \text{otherwise},
    \end{cases}
\end{equation}
in which $\zeta_{k,i}=(v_{k,i} - v_{k,i_1} ){\hat ~}$, $\sigma_{k,i}$ is the collision parameter for the pair, and $I\in\bbR^{6\times 6}$ is the identity matrix. The adjoint equation evolves backward in time with a given final condition
\begin{equation}\label{gammaFinal}
\beta_{M,i} = - \nabla_v r(v_{M,i}), \quad 1\leq i \leq N.
\end{equation}
The gradient of the objective function with respect to the parameter $m$ is
\begin{equation*}
    \partial_m \cJ =  - \frac{1}{N}\sum_{i=1}^N \beta_{0,i} \cdot  \partial_m v_{0,i}(m),
\end{equation*}
which can be computed using the adjoint variables $\{\beta_{0,i}\}$. They are solved backward in time following~\eqref{eq:adjoint_law} with the final condition~\eqref{gammaFinal}.

For more general optimization problems involving  the Boltzmann equation, that are different from~\eqref{BoltzmannIC} and~\eqref{eq:OTD_obj}, this procedure can be directly followed to obtain a more general version of the total objective function~\eqref{eq:DTO_obj}, which can then be optimized as above.

The adjoint DSMC method enables efficient computation of gradients for Boltzmann-constrained optimization problems at a cost that is largely independent of the number of optimization parameters, making it particularly well suited for high-dimensional and stochastic settings. Moreover, adjoint variable backpropagation closely mirrors the structure of the forward DSMC simulation, allowing gradient information to be propagated backward along particle trajectories using the same random samples generated in the forward run, thereby avoiding additional sampling and incurring only minimal computational overhead.

\section{Adjoint DSMC Method With Periodic and Specular Reflecting Boundary Conditions}\label{sec:adjoint_DSMC_space}

The adjoint DSMC method reviewed in~\Cref{sec:adjoint_DSMC} is based on the spatially homogeneous Boltzmann equation~\eqref{eq:homoBoltz}, as developed in~\cite{caflisch2021adjoint}. We now extend this framework to the spatially inhomogeneous Boltzmann equation~\eqref{eq:Boltz1}.

{
To illustrate the methodology, we first consider a simple Boltzmann-constrained optimization problem associated with~\eqref{eq:Boltz1}. The initial condition is given by
\begin{equation}
\label{eq:BoltzmannIC_space}
f(x,v,0) = f_0(x,v; m),
\end{equation}
where $f_0$ is prescribed initial data depending on the parameter $m$. For a general parameterized initial distribution, the initial particles may be generated from parameter-independent random seeds through a deterministic differentiable sampling map $\Psi$:
\[
(x_{0,i},v_{0,i})=\Psi(\xi_i;m),
\]
where $\xi_i$ is a random variable whose distribution does not depend on $m$. In what follows, we write $x_{0,i}=x_{0,i}(m)$ and $v_{0,i}=v_{0,i}(m)$ to emphasize the dependence of the initial samples on the parameter.
}

The objective is to determine parameters $m$ that optimize a terminal-time functional of the form
\begin{equation}
\label{eq:OTD_obj_space}
J_1(m) = \int_{\Omega} \int_{\mathbb{R}^3} r(x,v)\, f(x,v,T)\, dv\, dx,
\end{equation}
where $f(x,v,T)$ denotes the solution of~\eqref{eq:Boltz1} at time $t=T$ corresponding to the initial condition~\eqref{eq:BoltzmannIC_space}. The terminal distribution $f(x,v,T)$ depends on the parameter $m$ through the initial condition.

As before, the specific choices of parameterization and objective function are purely illustrative, and the proposed adjoint DSMC framework applies to more general parameterizations and objective functions.

Using adjoint variables, we define
\begin{eqnarray}
\cJ_1 &=& \frac{1}{N} \sum_{i=1}^N r\!\left(v_{M,i}, x_{M,i}\right), \label{J1}\\
\cJ_2 &=& \frac{1}{N} \sum_{i=1}^{N}  \beta_{0,i} \cdot \bigl(v_{0,i} - v_{0,i}(m)\bigr)
        + \frac{1}{N} \sum_{i=1}^{N} \alpha_{0,i} \cdot \bigl(x_{0,i} - x_{0,i}(m)\bigr), \label{J2}\\
\cJ_3 &=& \frac{1}{N} \sum_{k=0}^{M-1}\sum_{i\in N^k_{\mathrm{in}}}
          \beta_{k+1,i} \cdot \bigl(v_{k+1,i}-v_{k,i}'\bigr)
        + \frac{1}{N} \sum_{k=0}^{M-1} \sum_{i\in N^k_{\mathrm{in}}}
          \alpha_{k+1,i} \cdot \bigl(x_{k+1,i}-x'_{k,i}\bigr), \label{J3}
\end{eqnarray}
where
\[
N^k_{\mathrm{in}} = \{ i \in \{1,\ldots,N\} : x'_{k,i} \in \Omega \}
\]
denotes the set of particle indices whose post-collision positions remain inside the domain $\Omega$, and $v'_{k,i}$ denotes the post-collision velocity associated with $v_{k,i}$. We denote $\beta_{k,i}$ as the adjoint variable corresponding to the particle velocity $v_{k,i}$ and $\alpha_{k,i}$ as the adjoint variable corresponding to the particle position $x_{k,i}$, for $i=1,\ldots,N$ and $k=0,\ldots,M$.

Here, $\cJ_1$ represents the Monte Carlo approximation of the objective function~\eqref{eq:OTD_obj_space} at the final time, where $r(x,v)$ denotes the observable. The terms in $\cJ_2$ enforce the constraints on the initial particle velocities and positions through the adjoint variables $\{\beta_{0,i}\}$ and $\{\alpha_{0,i}\}$, respectively, with the joint initial samples $(x_{0,i}(m),v_{0,i}(m))=\Psi(\xi_i;m)$ generated by the general sampling map introduced above; no independence between the initial position and velocity is required. The terms in $\cJ_3$ impose the collision and advection constraints using the adjoint variables $\{\beta_{k,i}\}$ and $\{\alpha_{k,i}\}$ for particles that remain inside the domain after the update ($x'_{k,i} \in \Omega$). For particles that exit the domain, BCs are applied, and the corresponding particle variables $\{v_{k,i}\}$ and $\{x_{k,i}\}$ are no longer subject to the constraints in $\cJ_3$.

We conclude this subsection with a few remarks on notation that will appear repeatedly. In the forward DSMC algorithm~(Algorithm~\ref{alg:DSMC}), at each time step $t_k$ particles are grouped into pairs $(v_{k,i}, v_{k,i_1})$, where $i_1$ denotes the pairing index of particle $i$ at time $t_k$. Depending on the collision kernel, a subset of these pairs is selected to undergo collisions.

\begin{itemize}
    \item If the pair $(v_{k,i}, v_{k,i_1})$ is selected for collision in the forward DSMC process at time $t_k$, one would sample the parameter $\sigma_{k,i}$, the post-collision relative velocity direction, and compute $\zeta_{k,i}$, the pre-collision relative velocity. In this case, we use the $6\times6$ matrix $B(\sigma_{k,i}, \zeta_{k,i})$ defined in~\eqref{eq:AB_def} in the adjoint equations.
    \item If the pair $(v_{k,i}, v_{k,i_1})$ is \emph{not} selected for collision, we slightly abuse notation and set $B_{k,i} = I$, where $I$ denotes the $6\times6$ identity matrix.
\end{itemize}

Combining both cases, we define the operator $B_{k,i}$ by
\begin{equation}
\label{eq:B_ki}
B_{k,i} =
\begin{cases}
B(\sigma_{k,i}, \zeta_{k,i}), & \text{if } (v_{k,i}, v_{k,i_1}) \text{ collides at time } t_k,\\
I, & \text{otherwise.}
\end{cases}
\end{equation}
An important observation is that $B_{k,i} = B_{k,i_1}$.

Next, we discuss periodic and specular reflection BCs. Additional Lagrangian terms associated with these BCs must be introduced.

\subsection{Periodic Boundary Condition}
\label{subsec:P}

We begin by considering a 1D spatial domain $\Omega = [L,R] \subset \mathbb{R}$. The periodic BC enforces the following update whenever $x'_{k,i} \notin \Omega$:
\begin{equation}
x_{k+1,i} = L + \bigl( x'_{k,i} \bmod (R-L) \bigr),
\end{equation}
where $a \bmod b$ denotes the smallest nonnegative real number $c$ such that $(a-c)/b$ is an integer. Under periodic BCs, the particle velocity is unchanged by boundary crossing, i.e., $v_{k+1,i} = v'_{k,i}$ regardless of whether $x'_{k,i}$ lies inside $\Omega$. The above formulation readily extends to periodic BCs on more general spatial domains in $\mathbb{R}^d$.

For simplicity, we present the derivation in the 1D setting $\Omega = [L,R]$. In this case, the functional $\cJ_3$ in~\eqref{J3} is replaced by the following expression, which incorporates the periodic BC:
\begin{eqnarray}
\cJ_3^{p1} &=& \frac{1}{N} \sum_{k=0}^{M-1} \sum_{i=1}^N \beta_{k+1,i} \cdot \bigl(v_{k+1,i}-v'_{k,i}\bigr) \nonumber\\
&& + \frac{1}{N} \sum_{k=0}^{M-1} \sum_{i=1}^N \alpha_{k+1,i}
\biggl(x_{k+1,i}-\Bigl(L + (x'_{k,i} \bmod (R-L))\Bigr)\biggr).
\label{J5_p1}
\end{eqnarray}
We consider the unconstrained optimization problem defined by the augmented objective
\[
\cJ = \cJ_1 + \cJ_2 + \cJ_3^{p1},
\]
where $\cJ_1$ and $\cJ_2$ are given by~\eqref{J1}-\eqref{J2}. The superscript $p1$ indicates the 1D periodic setting.

The initial condition, collision rule, advection rule, and BC are recovered by differentiating $\cJ$ with respect to the adjoint variables
\[
\{\beta_{k,i}\}, \quad \{\alpha_{k,i}\}, \qquad k = 0,\ldots,M,\quad i = 1,\ldots,N.
\]
The adjoint equations are then obtained by setting the derivatives of $\cJ$ with respect to all state variables $\{v_{k,i}\}$ and $\{x_{k,i}\}$ to zero.

For any spatial dimension $d=1,2,3$, the terminal conditions for the adjoint variables $\beta_{M,i}$ and $\alpha_{M,i}$, $i=1,\ldots,N$, are given by
\begin{eqnarray}
\beta_{M,i} &=& -\partial_v r(x_{M,i},v_{M,i}), \label{eq:J_v_final}\\
\alpha_{M,i} &=& -\partial_x r(x_{M,i},v_{M,i}). \label{eq:J_x_final}
\end{eqnarray}
The adjoint variables $\{\beta_{k,i}\}$ and $\{\alpha_{k,i}\}$ are then propagated backward in time for $k=0,\ldots,M-1$ according to
\begin{eqnarray}
\begin{bmatrix}
\beta_{k,i} \\
\beta_{k,i_1}
\end{bmatrix}
&=&
B_{k,i}
\begin{bmatrix}
\beta_{k+1,i} + \Delta t\, \mathcal P^* \alpha_{k+1,i} \\
\beta_{k+1,i_1} + \Delta t\, \mathcal P^* \alpha_{k+1,i_1}
\end{bmatrix},
\quad \text{for all collision pairs $(i,i_1)$},
\label{eq:adjoint_gamma_P}\\
\alpha_{k,i} &=& \alpha_{k+1,i}, \quad \text{for all } i=1,\ldots,N,
\label{eq:adjoint_alpha_P}
\end{eqnarray}
where $B_{k,i}$ is defined in~\eqref{eq:B_ki}. The adjoint projection operator $\mathcal P^*$ is given by $\mathcal P^*=\mathbf e_1$ in one dimension, $\mathcal P^*=[\mathbf e_1,\mathbf e_2]$ in two dimensions, and $\mathcal P^*=I$ in three dimensions.

Once the adjoint variables $\{\beta_{0,i}\}$ and $\{\alpha_{0,i}\}$ at the initial time $t=0$ have been computed, the gradient of the objective function with respect to the parameters in the initial velocity and spatial distributions can be evaluated as{
\begin{align}\label{eq:grad_m}
\partial_{m} \mathcal J
=
-\frac{1}{N}\sum_{i=1}^N
\left[
\beta_{0,i}\cdot \nabla_{m} v_{0,i}(m)
+
\alpha_{0,i}\cdot \nabla_{m} x_{0,i}(m)
\right].
\end{align}
}

\subsection{Reflecting Boundary Condition}
\label{subsec:R}

Next, we consider the specular (mirror) reflection boundary condition, in which a particle undergoes mirror reflection upon reaching the boundary of the spatial domain $\Omega$.

Let $(v^*, x^*)$ denote the velocity and spatial position of a particle with $x^* \notin \Omega$. For the moment, we assume that $v^*$ and $x^*$ have the same dimension. Let $\mathbf n$ denote the outward unit normal vector at the boundary. The hyperplane defining the spatial reflection is given by $\{x \cdot \mathbf n = c\}$ for a constant $c$ determined by the boundary, while the corresponding hyperplane for velocity reflection is $\{v \cdot \mathbf n = 0\}$. The reflected particle has velocity $v$ and position $x$ given by
\begin{eqnarray*}
v &=& v^* - 2\mathbf n(\mathbf n \cdot v^*),\\
x &=& x^* - 2\mathbf n(\mathbf n \cdot x^* - c).
\end{eqnarray*}
Throughout, we neglect the possibility that a particle with sufficiently large velocity may cross the boundary more than once during a single time step, but this is easily included in the algorithm.

We now specialize to the setting in which the velocity domain is $\mathbb R^3$, while the spatial domain $\Omega$ is a subset of $\mathbb R$, $\mathbb R^2$, or $\mathbb R^3$. For clarity, we present the derivation in the 1D case and then state the final adjoint equations in full generality.

Consider the 1D spatial domain $\Omega = [L,R] \subset \mathbb R$, with $v^* \in \mathbb R^3$ and $x^* \in \mathbb R$. If $x^* \notin \Omega$, the reflecting BC takes the form
\begin{eqnarray*}
x &=& b - |x^* - b|\,\mathrm{sgn}(v^* \cdot \mathbf e_1) = 2b - x^*,\\
v &=& [v_x, v_y, v_z]^\top = [-v_x^*,\, v_y^*,\, v_z^*]^\top,
\end{eqnarray*}
where $b=L$ if $x^*<L$ and $b=R$ if $x^*>R$. Here, $\mathbf e_1 = [1,0,0]^\top$, and $\mathrm{sgn}(\cdot)$ is the sign function.

To enforce the reflecting BC within the adjoint formulation, we modify the Lagrangian to incorporate additional constraints beyond those in~\eqref{J1}-\eqref{J2}. For each time step $t_k$, the post-collision positions $x'_{k,i}$ are partitioned into three disjoint sets: $N^k_{\mathrm{in}}$, consisting of particles with $x'_{k,i} \in [L,R]$; $N^k_{\mathrm{out}_L}$, consisting of particles with $x'_{k,i}<L$; and $N^k_{\mathrm{out}_R}$, consisting of particles with $x'_{k,i}>R$. We replace~\eqref{J3} by
\begin{eqnarray*}
\cJ_3^{r1} &=& \frac{1}{N} \sum_{k=0}^{M-1} \sum_{i=1}^N
\beta_{k+1,i} \cdot \bigl(v_{k+1,i} - C_{k,i} v'_{k,i}\bigr) \\
&& + \frac{1}{N} \sum_{k=0}^{M-1} \sum_{i=1}^N
\alpha_{k+1,i} \Bigl(x_{k+1,i}
- \mathds{1}_{i\in N^k_{\mathrm{out}_L}}(2L-x'_{k,i})
- \mathds{1}_{i\in N^k_{\mathrm{out}_R}}(2R-x'_{k,i})
- \mathds{1}_{i\in N^k_{\mathrm{in}}}x'_{k,i}\Bigr),
\label{eq:J3_R1}
\end{eqnarray*}
where the diagonal matrix
\begin{equation}\label{eq:C_ki}
C_{k,i} = \mathrm{diag}\!\left([\,2\,\mathds{1}_{i\in N^k_{\mathrm{in}}}-1,\;1,\;1\,]\right)
\end{equation}
enforces the velocity reflection in the normal direction.

The full objective function is then given by
\[
\cJ = \cJ_1 + \cJ_2 + \cJ_3^{r1}.
\]
Proceeding as in~\Cref{subsec:P}, the adjoint equations are obtained by setting the derivatives of $\cJ$ with respect to all state variables equal to zero.

When the spatial domain is a connected interval (1D), a rectangle (2D), or a cuboid (3D), the resulting adjoint equations simplify as follows. The terminal conditions for the adjoint variables remain identical to~\eqref{eq:J_v_final} and~\eqref{eq:J_x_final}. For each particle index $i=1,\ldots,N$ and time step $k=0,\ldots,M-1$, the adjoint variables satisfy
\begin{eqnarray}
\begin{bmatrix}
\beta_{k,i} \\
\beta_{k,i_1}
\end{bmatrix}
&=&
B_{k,i}
\begin{bmatrix}
C_{k,i}\bigl(\beta_{k+1,i} + \Delta t\, \mathcal P^* \alpha_{k+1,i}\bigr) \\
C_{k,i_1}\bigl(\beta_{k+1,i_1} + \Delta t\, \mathcal P^* \alpha_{k+1,i_1}\bigr)
\end{bmatrix},
\quad \text{for all collision pairs $(i,i_1)$},
\label{eq:adjoint_gamma_R}\\
\mathcal P^* \alpha_{k,i}
&=&
C_{k,i}\,\mathcal P^* \alpha_{k+1,i},
\quad \text{for all } i=1,\ldots,N,
\label{eq:adjoint_alpha_R}
\end{eqnarray}
where $B_{k,i}$ is defined in~\eqref{eq:B_ki}, the adjoint projection operator $\mathcal P^*$ is defined dimension-wise as in~\Cref{subsec:P}, and the matrix $C_{k,i}$ depends on the spatial dimension according to~\eqref{eq:C_ki}.

Finally, once the adjoint variables have been fully back-propagated, the gradients with respect to the parameter $m$ are computed using~\eqref{eq:grad_m}.

\section{Adjoint DSMC Method With Thermal Boundary Conditions}
\label{sec:adjoint_DSMC_thermal}

In Section~\ref{sec:adjoint_DSMC_space}, we discussed strategies for handling periodic and specular (mirror) reflecting BCs. Here, we focus on the thermal BC, which is more complex than the previous cases and requires a slight modification of the forward DSMC algorithm. Note that different BCs for the forward DSMC algorithm only affect Line~13 of Algorithm~\ref{alg:DSMC}.

{For the thermal BC}, a particle that hits the boundary is absorbed and then re-emitted into the domain with a velocity randomly drawn from the thermal equilibrium flux at the wall, which is given by a half-Maxwellian distribution~\cite{pareschi2001introduction}. Given an  outward unit normal vector $n$, the reflected velocity $v$ is sampled from a half-Maxwellian flux satisfying $v \cdot n < 0$.

We next describe the Boltzmann equation with thermal BCs for the 1D spatial domain $\Omega = [L,R]$. The thermal BC is enforced whenever the post-collision spatial position $x'_{k,i}$ satisfies either $x'_{k,i} < L$ or $x'_{k,i} > R$. In this case, the original velocity is replaced by a sample from the half-Maxwellian flux with temperature parameter $T_L$ (for the left boundary) or $T_R$ (for the right boundary), respectively, and the spatial position is determined by the thermal BC rule~\cite[Alg.~6.3]{pareschi2001introduction}:

\begin{itemize}
    \item If $x'_{k,i} < L$, we sample $g_{k,i} \in \mathbb{R}^3$ from the half-Maxwellian flux
    $-v \cdot n_L \, M_{T_L}(v) \, \mathds{1}_{v \cdot n_L < 0}$ (denoted by $|M_{T_L}|$ hereafter) where $n_L = -\textbf{e}_1$.
    Thus, $e_1^\top g_{k,i} > 0$. The new velocity and spatial position are given by
    \begin{equation}
    \label{eq:thermal_v_L}
    v_{k+1,i} = g_{k,i}, \quad
    x_{k+1,i} = L + \frac{x'_{k,i} - L}{e_1^\top v'_{k,i}} \, e_1^\top g_{k,i}.
    \end{equation}

    \item If $x'_{k,i} > R$, we sample $h_{k,i} \in \mathbb{R}^3$ from the half-Maxwellian flux
    $-v \cdot n_R \, M_{T_R}(v) \, \mathds{1}_{v \cdot n_R < 0}$ (denoted by $-|M_{T_R}|$ hereafter) where $n_R = \textbf{e}_1$.
    Thus, $e_1^\top h_{k,i} < 0$. The new velocity and spatial position are given by
    \begin{equation}
    \label{eq:thermal_v_R}
    v_{k+1,i} = h_{k,i}, \quad
    x_{k+1,i} = R + \frac{x'_{k,i} - R}{e_1^\top v'_{k,i}} \, e_1^\top h_{k,i}.
    \end{equation}
\end{itemize}
The position and velocity updates~\eqref{eq:thermal_v_L}  and~\eqref{eq:thermal_v_R} do not include the possibility of multiple boundary crossings, but these can be easily added to the algorithm.

The algorithm in~\eqref{eq:thermal_v_L}-\eqref{eq:thermal_v_R} is straightforward, but it results in a \emph{random and  discontinuous} change in particle velocity at the instant a particle interacts with the boundary. While this discontinuity is not problematic for forward DSMC simulations, it poses significant challenges for the derivation and implementation of adjoint DSMC methods.

We emphasize that this difficulty is not specific to the Boltzmann equation or to DSMC. More generally, for any evolution equation involving an advection operator that is solved using particle-based methods (or, equivalently, the method of characteristics), thermal BCs induce an instantaneous resampling of particle velocities at the boundary. This resampling produces an inherent discontinuity in the particle velocity at the time of boundary interaction. While such discontinuities are entirely benign for forward simulations, they present a fundamental challenge for adjoint-based sensitivity analysis and gradient computation, as classical differentiation across the boundary is no longer well defined.

To address this issue, we introduce a \emph{stochastic modification} of the algorithm by replacing the deterministic time step $\Delta t$ in~\eqref{eq:x_update} with a random variable $\tau$. Specifically, we let $\tau$ follow a normal distribution $\mathcal{N}(\Delta t, \varepsilon^2)$, where $\Delta t$ is the target time step and $\varepsilon \ll \Delta t/3$ is chosen so that the probability of $\tau < 0$ is negligible.

For fixed $\varepsilon$, this construction defines an $\varepsilon$-regularized DSMC scheme and an associated objective $J_\varepsilon$. Because the randomized time increment changes the probability of reaching the wall and the post-boundary particle position, the regularized and deterministic forward objectives need not agree exactly, even in expectation. As $\varepsilon\to0$, however, the randomized update converges to the deterministic DSMC update. The numerical experiments in Section~\ref{sec:numerics} examine the dependence of the objective and gradient on $\varepsilon$ and show that the regularization effect is small over the range used in our calculations.

In the present work, this randomization of the boundary interaction serves as a numerical smoothing mechanism that enables the derivation of adjoint equations, which require differentiation of the particle velocity with respect to the particle position, despite the discontinuous dependence of velocity on position induced by the thermal BC. In other contexts, similar randomizations are introduced to model additional physical effects; here, however, the purpose is purely algorithmic.

The introduction of the random time step $\tau$ also has an important analytical consequence for the adjoint formulation. Because the boundary interaction is now described probabilistically, derivatives of expectations with respect to particle states can be computed using the score-function. In particular, sensitivities no longer require differentiating the discontinuous boundary map itself; instead, they involve differentiation of the log-density of $\tau$, which is smooth. This replaces an ill-defined pathwise derivative at the boundary by a well-defined weak derivative in expectation. We now make this statement precise.

Under this modification, Line~9 of \Cref{alg:DSMC}, i.e., Equation~\eqref{eq:x_update}, is replaced by
\begin{equation}
\label{eq:x_update_new}
x_{k,i}' = x_{k,i} + \tau_{k,i}\, \mathcal P v'_{k,i}, \qquad \mathcal P v = e_1^\top v,
\end{equation}
corresponding to the case of a 1D spatial domain and 3D velocity space. Since $\tau_{k,i}$ is a random variable, the three events $\{x'_{k,i} < L\}$, $\{x'_{k,i} \in \Omega\}$, and $\{x'_{k,i} > R\}$ occur with probabilities $p_l(x_{k,i}, v'_{k,i}; \Delta t, \varepsilon)$, $l=1,2,3$, respectively. They sum to one: $\sum_{l=1}^3 p_l(x_{k,i}, v'_{k,i}; \Delta t, \varepsilon) = 1$. Explicit formulas for $p_l(x,v;\Delta t,\varepsilon)$ are provided in~\ref{app:log}.

{Within each global time interval $[t_k,t_{k+1}]$, collisions are first performed using the particle states synchronized at $t_k$. The particles are then advanced through the advection and boundary update, with the randomized time used only in this latter step; after the update, all particles are assigned to the common discrete time $t_{k+1}$.}

Since $\Delta t$ and $\varepsilon$ are fixed hyperparameters, we henceforth write $p_l(x_{k,i}, v'_{k,i})$ for notational convenience. If $x'_{k,i} < L$ or $x'_{k,i} > R$, the updated velocity and position $(v_{k+1,i}, x_{k+1,i})$ are determined according to~\eqref{eq:thermal_v_L} and~\eqref{eq:thermal_v_R}, respectively. Otherwise, the particle remains inside the domain and we set $x_{k+1,i} = x'_{k,i}$ and $v_{k+1,i} = v'_{k,i}$.

For any test function $\varphi(x,v)$, the conditional expectation of the post-update state satisfies
\begin{eqnarray}
\label{eq:thermal_macro}
&& \mathbb E_{\tau_{k,i}}\!\left[\varphi(x_{k+1,i},v_{k+1,i}) \,\big|\, (x_{k,i}, v'_{k,i})\right]\\
&=& p_1\!\left(x_{k,i}, v'_{k,i}\right)
\mathbb E_{\tau_{k,i}} \mathbb E_{g_{k,i}}
\left[
\varphi\!\left(
L + \frac{x'_{k,i}-L}{e_1^\top v'_{k,i}}\, e_1^\top g_{k,i},
\, g_{k,i}
\right)
\,\big|\, (x_{k,i}, v'_{k,i})\right] +\nonumber \\
&& p_2\!\left(x_{k,i}, v'_{k,i}\right)
\mathbb E_{\tau_{k,i}}
\left[
\varphi(x'_{k,i}, v'_{k,i})
\,\big|\, (x_{k,i}, v'_{k,i})\right] + \nonumber \\
&&  p_3\!\left(x_{k,i}, v'_{k,i}\right)
\mathbb E_{\tau_{k,i}}\mathbb E_{h_{k,i}}
\left[
\varphi\!\left(
R + \frac{x'_{k,i}-R}{e_1^\top v'_{k,i}}\, e_1^\top h_{k,i},
\, h_{k,i}
\right)
\,\big|\, (x_{k,i}, v'_{k,i})\right],\nonumber
\end{eqnarray}
where $x'_{k,i}$ depends on $\tau_{k,i}$ through~\eqref{eq:x_update_new}. Equation~\eqref{eq:thermal_macro} highlights that, in a 1D spatial domain, each of the three possible post-update scenarios contributes two distinct sources of dependence on the pre-update state $(x_{k,i}, v'_{k,i})$:
\begin{enumerate}
    \item the probabilities $p_l(x_{k,i}, v'_{k,i})$, $l=1,2,3$, and
    \item the observable $\varphi$ evaluated at $(x_{k+1,i}, v_{k+1,i})$.
\end{enumerate}

The probabilistic decomposition in~\eqref{eq:thermal_macro} is closely related to the probabilistic structure arising in DSMC methods that employ \emph{virtual collisions} to efficiently and accurately sample collision events. In that setting, three mutually exclusive outcomes are considered, (1)~no virtual collision, (2)~a virtual collision that corresponds to a real collision, and~(3)~a virtual collision that does not result in a real collision, each occurring with a prescribed probability; see~\cite[Eqns.~(11)-(13)]{yang2023adjoint}. A key contribution of our previous work~\cite{yang2023adjoint} was the use of the score-function method~\cite{rubinstein1986score} to account for the sensitivity of such probability terms with respect to optimization parameters. This approach enables accurate gradient computation in adjoint DSMC methods, even when the underlying dynamics involve discrete probabilistic branching, as in the present thermal boundary setting.

\subsection{The Lagrangian}
Similar to steps in Section~\ref{sec:adjoint_DSMC_space}, we first write down the Lagrangian to include the constraints. We denote the expectations over the step size $\tau_k$ in the $k$th time step by $\mathbb{E}_{\tau_{k,i}}$,  and the expectations over the new velocities drawn from the thermal boundary condition by $\mathbb{E}_{g_{k,i},h_{k,i}} $. We define $\mathbb{E}^k = \mathbb{E}_{\{\tau_{k,i}\}} \mathbb{E}_{\{g_{k,i}\},\{h_{k,i}\}}$, the expectation over all random time steps and the wall emitted velocities, and $\bbE$ as $$\bbE = \bbE^0 \bbE^1 \cdots \bbE^{M-1}.$$

Due to the introduction of stochastic time steps, we need to modify $\cJ_1$ and $\cJ_3$ in~\eqref{J1} and~\eqref{J3}, respectively:
\begin{align}
\cJ_1^{t1} = & \frac{1}{N} \sum_{i=1}^N  \mathbb{E} \left[ r( x_{M,i}, v_{M,i}) \right] \label{eq:J1_T1} ,  \\
\cJ_{3}^{t1} = & \frac{1}{N}\sum_{k=0}^{M-1} \sum_{i=1}^N  \mathbb{E}^k \left[  \beta_{k+1,i}  \cdot  \left( v_{k+1,i}   -  \mathcal{R}^v \left( x_{k,i}\,,\, v_{k,i}'  \right) \right)   \bigg| \left( \cV_{k}, \cX_k \right) \right] +  \label{eq:J3_T1}  \\
&   \frac{1}{N}\sum_{k=0}^{M-1} \sum_{i=1}^N  \mathbb{E}^k \left[ \alpha_{k+1,i}  \cdot \left(  x_{k+1,i}   -  \mathcal{R}^x \left( x_{k,i}\,,\, v_{k,i}'  \right) \right)   \bigg|  \left( \cV_{k}, \cX_k \right)  \right] . \nonumber
\end{align}
The operator $\mathcal R(x,v)=(\mathcal R^x(x,v),\mathcal R^v(x,v))$ enforces the thermal BC as follows: Define the three branch maps $\mathcal R_\ell(x,v;\tau,g,h)=(\mathcal R^x_\ell,\mathcal R^v_\ell)$ by
\[
\mathcal R_1(x,v;\tau,g)=
\begin{bmatrix}
\left(1-\frac{\mathcal Pg}{\mathcal Pv}\right)L+\frac{\mathcal Pg}{\mathcal Pv}\,x+\tau\,\mathcal Pg\\[0.2em]
g
\end{bmatrix},\qquad
\mathcal R_2(x,v;\tau)=
\begin{bmatrix}
x+\tau\,\mathcal Pv\\ v
\end{bmatrix},
\]
\[
\mathcal R_3(x,v;\tau,h)=
\begin{bmatrix}
\left(1-\frac{\mathcal Ph}{\mathcal Pv}\right)R+\frac{\mathcal Ph}{\mathcal Pv}\,x+\tau\,\mathcal Ph\\[0.2em]
h
\end{bmatrix}.
\]
The event $\ell\in\{1,2,3\}$ is selected according to the probabilities $p_\ell(x,v';\Delta t,\varepsilon)$ induced by $\tau\sim\mathcal N(\Delta t,\varepsilon^2)$, namely $\ell=1$ if $x+\tau\mathcal Pv'<L$, $\ell=2$ if $L\le x+\tau\mathcal Pv'\le R$, and $\ell=3$ if $x+\tau\mathcal Pv'>R$.
Here, $\mathcal{P}$ is the projection onto the first component for the 1D spatial domain. Each term in $\cJ_{3}^{t1}$ should be interpreted following~\eqref{eq:thermal_macro} for a particular function $\varphi$.

\subsection{Adjoint equation derivation}

With these preparations, we will derive the optimality conditions for the Lagrangian
\[
\cJ = \cJ_1^{t1} +  \cJ_2 +   \cJ_{3}^{t1}\,,
\]
combining terms in~\eqref{eq:J1_T1}, \eqref{J2} and \eqref{eq:J3_T1}.

\subsubsection{Recovering forward DSMC} The collision and advection rules are derived from  the derivatives of $\mathcal J$  with   respect to $ \beta_{k,i}$ and $\alpha_{k,i}$ for $k = 0,\ldots, M-1$, and $i = 1,\ldots, N$. Similarly, we recover the initial conditions for $x_{0,i}$ and $v_{0,i}$, $i = 1,\ldots, N$.

\subsubsection{Gradient for $m$} The gradient computation still follows Equations~\eqref{eq:grad_m}.

\subsubsection{Final-time condition for the adjoint system} For each $i$,  we take the derivative of $\mathcal J$ with respect to the final velocity particle $v_{M,i}$ and the final location $x_{M,i}$. Recall that
$$
\cJ_1^{t1} =  \frac{1}{N} \mathbb{E} \left[ \sum_{i=1}^N r(x_{M,i}, v_{M,i}) \right]   \approx \frac{1}{N} \mathbb{E} \left[ \sum_{i=1}^N r(x_{M,i}, v_{M,i})
| (\cV_M, \cX_M) \right]  = \frac{1}{N}  \sum_{i=1}^N r(x_{M,i}, v_{M,i})\,.
$$
We then have
\[
\frac{\partial \cJ_{1}^{t1} }{\partial v_{M,i}}   =   \frac{1}{N} \partial_v r(x_{M,i}, v_{M,i}) \,,\quad  \frac{\partial \cJ_{1}^{t1} }{\partial x_{M,i}}   =   \frac{1}{N}
\partial_x r(x_{M,i}, v_{M,i}) \,.
\]
For $\cJ_{3}^{t1}$ enforcing the binary collision and the advection rules, since the expectation $\bbE^{M-1}$ does not depend on particles in $\mathcal{V}_M$ and spatial positions in $\cX_M$, we have
\[
\frac{\partial \cJ_{3}^{t1} }{\partial v_{M,i}}   =   \frac{1}{N} \bbE^{M-1}[ \beta_{M,i} |  \left( \mathcal{V}_{M-1} , \mathcal{X}_{M-1} \right)  ]  \,,\quad  \frac{\partial \cJ_{3}^{t1} }{\partial x_{M,i}}   =   \frac{1}{N}
\bbE^{M-1}[ \alpha_{M,i} | \left( \mathcal{V}_{M-1} , \mathcal{X}_{M-1} \right)] \,.
\]
The other terms  in  $\mathcal{J}$  do not depend on the particle positions $\cX_M$ and velocities $\mathcal{V}_M$ at the final time, and so do not contribute to the derivatives. Summarizing all the terms and setting $ \frac{\partial \cJ }{\partial v_{M,i}}  =  \frac{\partial \cJ }{\partial x_{M,i}}  = 0$, we obtain the same final conditions~\eqref{eq:J_v_final}-\eqref{eq:J_x_final}.

\subsubsection{Adjoint equations} \label{subsec:link_appendx_A}
Next, we derive the key equations for $\alpha_{k,i}$ and $\beta_{k,i}$ where $0\leq k \leq M-1$. For $\cJ_1^{t1}$, we observe that, for any $k = 0,\ldots, M-1$,
\begin{align}
\partial_{x_{k,i}} \cJ_1^{t1}
&= \frac{1}{N}\,\partial_{x_{k,i}}\, \bbE^{k}\cdots \bbE^{M-1}
\left[ \sum_{j=1}^N r(x_{M,j}, v_{M,j}) \,\middle|\, (\cV_{k}, \cX_{k}) \right] \nonumber \\
&= \frac{1}{N} \sum^{N}_{j=1}(\partial_{x_{k,i}}\bbE^{k}) \bbE^{k+1}\cdots \bbE^{M-1}
\left[r(x_{M,j}, v_{M,j}) \,\middle|\, (\cV_k,\cX_k)\right] \nonumber \\
&= \frac{1}{N} \sum^{N}_{j=1}\bbE^{k}\cdots \bbE^{M-1}
\left[  \delta_{j i}\, \left( \partial_{x_{k,i}} \log p_l(x_{k,j}, v'_{k,j})  \right) \, r(x_{M,j}, v_{M,j})
\,\middle|\, (\cV_k,\cX_k) \right]  \nonumber \\
&=  \frac{1}{N}\,\bbE^{k}\cdots \bbE^{M-1}
\left[
\left( \partial_{x_{k,i}} \log p_l(x_{k,i}, v'_{k,i}) \right)\,
r(x_{M,i}, v_{M,i})
\,\middle|\, (\cV_k,\cX_k)
\right]. \label{eq:J_final_t1_x_exact}
\end{align}
Here, we used the identity $\partial_\theta \int f(x) p(x,\theta)\,dx = \int f(x)\,(\partial_\theta \log p(x,\theta))\,p(x,\theta)\,dx$. In practice, we approximate the conditional expectation in~\eqref{eq:J_final_t1_x_exact} by a single-trajectory estimator, yielding
\begin{equation}
\label{eq:J_final_t1_x}
\partial_{x_{k,i}}\cJ_1^{t1} \approx \frac{1}{N}\,
(\partial_{x_{k,i}}\log p_l(x_{k,i},v'_{k,i}))\, r_i,
\end{equation}
where $r_i := r(x_{M,i}, v_{M,i})$.

Similarly,
\begin{align}
\partial_{v_{k,i}} \cJ_1^{t1}
&= \frac{1}{N}\,\partial_{v_{k,i}}\, \bbE^{k}\cdots \bbE^{M-1}
\left[ \sum_{j=1}^N r(x_{M,j}, v_{M,j}) \,\middle|\, (\cV_{k}, \cX_{k}) \right] \nonumber \\
&= \frac{1}{N} \sum^{N}_{j=1}(\partial_{v_{k,i}}\bbE^{k}) \bbE^{k+1}\cdots \bbE^{M-1}
\left[r(x_{M,j}, v_{M,j}) \,\middle|\, (\cV_k,\cX_k)\right] \nonumber \\
&= \frac{1}{N} \sum^{N}_{j=1}\bbE^{k}\cdots \bbE^{M-1}
\left[  \left(  \delta_{j i_1} + \delta_{j i} \right)\,
\left( \partial_{v_{k,i}} \log p_l(x_{k,j}, v'_{k,j})  \right)\,
r(x_{M,j}, v_{M,j})
\,\middle|\, (\cV_k,\cX_k) \right]  \nonumber \\
&= \frac{1}{N}\,\bbE^{k}\cdots \bbE^{M-1}
\Bigl[
\bigl(\partial_{v_{k,i}} \log p_l(x_{k,i}, v'_{k,i})\bigr)\, r_i
+
\bigl(\partial_{v_{k,i}} \log p_l(x_{k,i_1}, v'_{k,i_1})\bigr)\, r_{i_1}
\Bigm| (\cV_k,\cX_k)
\Bigr]. \nonumber
\end{align}
In practice, we approximate the conditional expectation in the last term above by a single-trajectory estimator and apply the chain rule through the collision map, obtaining
\begin{align}
\partial_{v_{k,i}} \cJ_1^{t1}
\approx \frac{1}{N}\Bigl[
\bigl(\partial_{v_{k,i}} v'_{k,i}\bigr)^{\top}
\bigl(\partial_{v'_{k,i}} \log p_l(x_{k,i}, v'_{k,i})\bigr)\, r_i
+
\bigl(\partial_{v_{k,i}} v'_{k,i_1}\bigr)^{\top}
\bigl(\partial_{v'_{k,i_1}} \log p_l(x_{k,i_1}, v'_{k,i_1})\bigr)\, r_{i_1}
\Bigr]. \label{eq:J_final_t1_v}
\end{align}
Here, particles with indices $i$ and $i_1$ form a collision pair at the $k$th time step. Explicit expressions for $\partial_v \log p_l(x,v)$ and $\partial_x \log p_l(x,v)$ are provided in~\ref{app:log}.

Next, we differentiate $\cJ_{3}^{t1}$ with respect to $x_{k,i}$ and $v_{k,i}$. Throughout, we view $\alpha_{k,i}$ and $\beta_{k,i}$ as column vectors, so Jacobians act on the adjoint variables through their transposes. We obtain
\begin{align}
\partial_{x_{k,i}} \cJ_{3}^{t1}
&= \frac{1}{N}\alpha_{k,i}
-\frac{1}{N}\bigl(D_x \mathcal{R}^x_{k,i}\bigr)^{\top}\alpha_{k+1,i}
-\frac{1}{N}\bigl(D_x \mathcal{R}^v_{k,i}\bigr)^{\top}\beta_{k+1,i},
\label{eq:dJ3_dx}
\\
\partial_{v_{k,i}} \cJ_{3}^{t1}
&= \frac{1}{N}\beta_{k,i}
-\frac{1}{N}\Bigl(\partial_{v_{k,i}} v'_{k,i}\Bigr)^{\top}\bigl(D_v \mathcal{R}^v_{k,i}\bigr)^{\top}\beta_{k+1,i}
-\frac{1}{N}\Bigl(\partial_{v_{k,i}} v'_{k,i_1}\Bigr)^{\top}\bigl(D_v \mathcal{R}^v_{k,i_1}\bigr)^{\top}\beta_{k+1,i_1}
\nonumber\\
&\quad
-\frac{1}{N}\Bigl(\partial_{v_{k,i}} v'_{k,i}\Bigr)^{\top}\bigl(D_v \mathcal{R}^x_{k,i}\bigr)^{\top}\alpha_{k+1,i}
-\frac{1}{N}\Bigl(\partial_{v_{k,i}} v'_{k,i_1}\Bigr)^{\top}\bigl(D_v \mathcal{R}^x_{k,i_1}\bigr)^{\top}\alpha_{k+1,i_1}.
\label{eq:dJ3_dv}
\end{align}

Here we use the shorthand Jacobian notation
\begin{equation}\label{eq:short_R}
\begin{aligned}
& D_x \cR^v_{k,i} := \frac{\partial \cR^v}{\partial x}(x_{k,i},v'_{k,i}),
\qquad  D_x \cR^x_{k,i} := \frac{\partial \cR^x}{\partial x}(x_{k,i},v'_{k,i}), \\
& D_v \cR^v_{k,i} := \frac{\partial \cR^v}{\partial v}(x_{k,i},v'_{k,i}),
\qquad  D_v \cR^v_{k,i_1} := \frac{\partial \cR^v}{\partial v}(x_{k,i_1},v'_{k,i_1}),\\
& D_v \cR^x_{k,i} := \frac{\partial \cR^x}{\partial v}(x_{k,i},v'_{k,i}),
\qquad  D_v \cR^x_{k,i_1} := \frac{\partial \cR^x}{\partial v}(x_{k,i_1},v'_{k,i_1}).
\end{aligned}
\end{equation}

These are components of the Jacobian of $\mathcal R(x,v)=(\mathcal R^x(x,v),\mathcal R^v(x,v))$. When a wall collision occurs, the sampled thermal velocity ($g$ or $h$) is treated as \emph{frozen} inside the derivative (i.e., we differentiate the map conditional on the sampled draw). Under this convention,
\begin{equation}
\label{eq:R_Jacobian}
D_{(x,v)} \mathcal{R}(x,v)
=
\begin{cases}
\begin{bmatrix}
D_x \mathcal R^x(x,v) & D_v \mathcal R^x(x,v) \\
0 & 0
\end{bmatrix},
& x<L,\; g\sim |M_{T_L}|, \\[0.25cm]
\begin{bmatrix}
I & \tau\,\mathcal P \\
0 & I
\end{bmatrix},
& L\le x \le R, \\[0.25cm]
\begin{bmatrix}
D_x \mathcal R^x(x,v) & D_v \mathcal R^x(x,v) \\
0 & 0
\end{bmatrix},
& x>R,\; h\sim -|M_{T_R}|.
\end{cases}
\end{equation}
Moreover, in the left-wall case ($x+\tau\mathcal P v<L$) we have
\[
D_x \mathcal R^x(x,v)=\frac{\mathcal P g}{\mathcal P v},
\qquad
D_v \mathcal R^x(x,v)=(L-x)\,\frac{\mathcal P g}{(\mathcal P v)^2}\,\mathcal P,
\]
and in the right-wall case ($x+\tau\mathcal P v>R$) the same formulas hold with $(L,g)$ replaced by $(R,h)$.
Setting $\partial_{x_{k,i}}\cJ = 0$ and $\partial_{v_{k,i}}\cJ = 0$, we obtain the adjoint equations
\begin{align}
\alpha_{k,i}
&=
\bigl(D_x \cR^x_{k,i}\bigr)^{\top}\alpha_{k+1,i}
+
\bigl(D_x \cR^v_{k,i}\bigr)^{\top}\beta_{k+1,i}
-
\partial_{x} \log p_l(x_{k,i},v'_{k,i}) \, r_i ,
\label{eq:alpha_adjoint}
\\[0.5ex]
\begin{bmatrix}
\beta_{k,i}\\
\beta_{k,i_1}
\end{bmatrix}
&=
B_{k,i}
\begin{bmatrix}
\bigl(D_v \cR^x_{k,i}\bigr)^{\top}\alpha_{k+1,i}
+
\bigl(D_v \cR^v_{k,i}\bigr)^{\top}\beta_{k+1,i}
-
\partial_{v} \log p_l(x_{k,i},v'_{k,i}) \, r_i
\\[0.5ex]
\bigl(D_v \cR^x_{k,i_1}\bigr)^{\top}\alpha_{k+1,i_1}
+
\bigl(D_v \cR^v_{k,i_1}\bigr)^{\top}\beta_{k+1,i_1}
-
\partial_{v} \log p_l(x_{k,i_1},v'_{k,i_1}) \, r_{i_1}
\end{bmatrix}.
\label{eq:beta_adjoint}
\end{align}

Finally, we have derived the adjoint system for the thermal BC, \eqref{eq:alpha_adjoint}-\eqref{eq:beta_adjoint} starting with the final conditions~\eqref{eq:J_v_final}-\eqref{eq:J_x_final}. We introduced a key ingredient of the framework by employing a random time step around $\Delta t$. Randomization is needed to regularize discontinuities arising from wall collisions, enabling gradient evaluation via the score-function method.

\begin{remark}
[Practical implementation] In practice, the randomization of the time step needs to be applied only to particles whose spatial bins lie within a small neighborhood of the boundary. Particles well inside the domain may be advanced using the deterministic step $\Delta t$ without affecting the adjoint formulation. This localized treatment significantly reduces unnecessary random sampling while preserving the differentiable structure required for adjoint sensitivity analysis. It is an easy adaptation of the adjoint DSMC algorithm presented in this section.
\end{remark}

\subsection{Derivatives with respect to thermal boundary parameters}

We now derive gradients of the objective function with respect to parameters appearing in the thermal boundary conditions. Let the thermal velocities at the left and right boundaries be generated through reparameterizations
\[
g = S_L(\theta_L,\epsilon), \qquad h = S_R(\theta_R,\epsilon'),
\]
where $\epsilon,\epsilon' \sim \mathcal Y$ are drawn from a parameter-independent reference distribution, and $S_L$ and $S_R$ are deterministic maps pushing forward $\mathcal Y$ to the half-Maxwellian flux distributions at the left and right walls, respectively. Typical choices of $\theta_L$ and $\theta_R$ include the wall temperatures and prescribed tangential wall velocities. Under this representation, the sensitivities of the sampled thermal velocities are
\begin{equation}
\partial_{\theta_L} g = \partial_{\theta_L} S_L(\theta_L,\epsilon),
\qquad
\partial_{\theta_R} h = \partial_{\theta_R} S_R(\theta_R,\epsilon').
\label{eq:dg_dh}
\end{equation}

{
For example, if the parameter $\theta_L$ or $\theta_R$ is the vector-valued temperature \(T=[T^1,T^2,T^3]^\top\), we use the following explicit reparameterization of the half-Maxwellian flux. Let \(U\sim\mathcal U(0,1)\) and \(Z_2,Z_3\sim\mathcal N(0,1)\) be independent of \(T\). At the left and right boundaries, respectively, we set
\[
g(T_L)=
\begin{bmatrix}
\sqrt{-2T_L^1\log U}\\
\sqrt{T_L^2}Z_2\\
\sqrt{T_L^3}Z_3
\end{bmatrix},
\qquad
h(T_R)=
\begin{bmatrix}
-\sqrt{-2T_R^1\log U'}\\
\sqrt{T_R^2}Z_2'\\
\sqrt{T_R^3}Z_3'
\end{bmatrix}.
\]
Consequently,
\[
\partial_{T_L}g
=
\frac12\operatorname{diag}
\left(\frac{g_1}{T_L^1},\frac{g_2}{T_L^2},\frac{g_3}{T_L^3}\right),
\qquad
\partial_{T_R}h
=
\frac12\operatorname{diag}
\left(\frac{h_1}{T_R^1},\frac{h_2}{T_R^2},\frac{h_3}{T_R^3}\right).
\]

}

Taking derivatives of the Lagrangian $\mathcal J$ with respect to $\theta_L$ and $\theta_R$ yields
\begin{align}
\partial_{\theta_L} \mathcal J
&=
-\frac{1}{N} \sum_{k=0}^{M-1} \sum_{i=1}^N
\mathbb E^k \!\left[
\beta_{k+1,i}^{\top}\,\partial_{\theta_L}\mathcal R^v(x_{k,i},v'_{k,i})
+
\alpha_{k+1,i}^{\top}\,\partial_{\theta_L}\mathcal R^x(x_{k,i},v'_{k,i})
\,\middle|\, (\mathcal V_k,\mathcal X_k)
\right], \label{eq:thermal_grad_L_expect} \\
\partial_{\theta_R} \mathcal J
&=
-\frac{1}{N} \sum_{k=0}^{M-1} \sum_{i=1}^N
\mathbb E^k \!\left[
\beta_{k+1,i}^{\top}\,\partial_{\theta_R}\mathcal R^v(x_{k,i},v'_{k,i})
+
\alpha_{k+1,i}^{\top}\,\partial_{\theta_R}\mathcal R^x(x_{k,i},v'_{k,i})
\,\middle|\, (\mathcal V_k,\mathcal X_k)
\right]. \label{eq:thermal_grad_R_expect}
\end{align}

After sampling the random variables at time step $k$, the conditional expectations collapse to single trajectory evaluations. Moreover, only particles undergoing thermal BC contribute to the gradients. Thus,
\begin{align}
\partial_{\theta_L} \mathcal J
&=
-\frac{1}{N} \sum_{k=0}^{M-1}
\sum_{\substack{i:\\ x_{k,i}+\tau\,\mathcal P v'_{k,i}<L}}
\Bigl[
\beta_{k+1,i}^{\top}\,\partial_{\theta_L}\mathcal R^v(x_{k,i},v'_{k,i})
+
\alpha_{k+1,i}^{\top}\,\partial_{\theta_L}\mathcal R^x(x_{k,i},v'_{k,i})
\Bigr],
\label{eq:thermal_grad_L}\\
\partial_{\theta_R} \mathcal J
&=
-\frac{1}{N} \sum_{k=0}^{M-1}
\sum_{\substack{i:\\ x_{k,i}+\tau\,\mathcal P v'_{k,i}>R}}
\Bigl[
\beta_{k+1,i}^{\top}\,\partial_{\theta_R}\mathcal R^v(x_{k,i},v'_{k,i})
+
\alpha_{k+1,i}^{\top}\,\partial_{\theta_R}\mathcal R^x(x_{k,i},v'_{k,i})
\Bigr].
\label{eq:thermal_grad_R}
\end{align}

When a left-wall collision occurs ($x+\tau\,\mathcal P v<L$), the reflection map depends on $\theta_L$ only through the sampled thermal velocity $g$. Treating the sampled draw as frozen inside the derivative, we obtain
{
\begin{equation}
\partial_{\theta_L}
\begin{bmatrix}
\mathcal R^x\\ \mathcal R^v
\end{bmatrix}
=
\partial_g \mathcal R \;\partial_{\theta_L} g
=
\begin{bmatrix}
\mathcal P\!\left(\dfrac{x-L}{\mathcal P v}+\tau\right)\\[0.3em]
I
\end{bmatrix}
\partial_{\theta_L} g .
\end{equation}
Similarly, when a right-wall collision occurs ($x+\tau\,\mathcal P v>R$),
\begin{equation}
\partial_{\theta_R}
\begin{bmatrix}
\mathcal R^x\\ \mathcal R^v
\end{bmatrix}
=
\partial_h \mathcal R \;\partial_{\theta_R} h
=
\begin{bmatrix}
\mathcal P\!\left(\dfrac{x-R}{\mathcal P v}+\tau\right)\\[0.3em]
I
\end{bmatrix}
\partial_{\theta_R} h .
\end{equation}
Here, $\partial_{\theta_L} g$ and $\partial_{\theta_R} h$ are given by~\eqref{eq:dg_dh}.
}

{
The reparameterization above also accommodates moving thermal walls. Let
\[
\Pi_\tau=I-\mathcal P^*\mathcal P
\]
denote the projection onto the tangential velocity space. The walls are fixed in the normal direction, so their prescribed velocities satisfy $\mathcal P U_L=\mathcal P U_R=0$. The diffusely re-emitted velocities in the laboratory frame can be written as
\[
g(U_L,T_L)=\Pi_\tau U_L+\widetilde g(T_L),
\qquad
h(U_R,T_R)=\Pi_\tau U_R+\widetilde h(T_R),
\]
where $\widetilde g$ and $\widetilde h$ are the corresponding half-Maxwellian samples in the wall frames. It follows that
\[
\partial_{U_L}g=\Pi_\tau,
\qquad
\partial_{U_R}h=\Pi_\tau.
\]
Because the wall velocity is purely tangential, $\mathcal P\Pi_\tau=0$. Therefore it does not directly affect the reflected particle position in the present one-dimensional spatial geometry:
\[
\partial_{U_L}\mathcal R^x
=
\partial_{U_R}\mathcal R^x
=0,
\qquad
\partial_{U_L}\mathcal R^v
=
\partial_{U_R}\mathcal R^v
=\Pi_\tau.
\]
Substitution into~\eqref{eq:thermal_grad_L}-\eqref{eq:thermal_grad_R} gives
\begin{align}
\nabla_{U_L}J
&=
-\frac{1}{N}
\sum_{k=0}^{M-1}
\sum_{\substack{
x_{k,i}+\tau\mathcal P v'_{k,i}<L}}
\Pi_\tau\beta_{k+1,i},
\label{eq:wall_velocity_grad_L}\\
\nabla_{U_R}J
&=
-\frac{1}{N}
\sum_{k=0}^{M-1}
\sum_{\substack{x_{k,i}+\tau\mathcal P v'_{k,i}>R}}
\Pi_\tau\beta_{k+1,i}.
\label{eq:wall_velocity_grad_R}
\end{align}
There is no direct $\alpha_{k+1,i}$ contribution because the tangential wall velocity does not enter $\mathcal R^x$. In addition, no new score-function term is required since the current boundary-crossing probabilities depend on the incoming particle state and the randomized time step, but not on the tangential velocity subsequently assigned at the wall. Indirect effects on later collisions and boundary interactions are included through the usual backward propagation of $\alpha$ and $\beta$.

For the symmetric Couette parameterization, for which we will show an example in Section~\ref{subsec:test4},
\[
U_L=-U_w\mathbf e_2,
\qquad
U_R=U_w\mathbf e_2.
\]
The chain rule therefore gives
\begin{align}
\frac{dJ}{dU_w}
=
-\mathbf e_2^\top\nabla_{U_L}J
+\mathbf e_2^\top\nabla_{U_R}J  =
\frac{1}{N}
\sum_{k=0}^{M-1}
\sum_{\substack{x_{k,i}+\tau\mathcal P v'_{k,i}<L}}
\beta_{k+1,i}^{2}
-
\frac{1}{N}
\sum_{k=0}^{M-1}
\sum_{\substack{x_{k,i}+\tau\mathcal P v'_{k,i}>R}}
\beta_{k+1,i}^{2},
\label{eq:couette_wall_speed_gradient}
\end{align}
where $\beta_{k+1,i}^{2}$ denotes the second component of the velocity adjoint $\beta_{k+1,i}$.}

\subsection{Summary}
In this section, we derived an adjoint DSMC formulation for the Boltzmann equation with thermal  boundary conditions. The main difficulty stems from the discontinuous velocity resampling induced by wall interactions, which obstructs classical adjoint differentiation. We resolved this by introducing a randomized time step around $\Delta t$, providing a regularization that smooths boundary interactions in expectation. The resulting adjoint computes the gradient of the regularized objective $J_\varepsilon$. Its relationship with the original deterministic DSMC objective is examined numerically through an $\varepsilon$-stability study in Section~\ref{sec:numerics}.

\section{Ensemble-Adjoint DSMC Method With Inflow Boundary Conditions}
\label{sec:adjoint_DSMC_inflow}

In this section, we derive an ensemble-adjoint DSMC method for the Boltzmann equation with inflow boundary conditions. We consider a one-dimensional spatial domain $\Omega=[L,R]$ and a three-dimensional velocity domain. Inflow boundary conditions are imposed at $x=L$ and $x=R$, and the incoming velocities are sampled from the corresponding boundary flux distributions, for example, from half-Maxwellian fluxes~\cite[Sec.~6]{pareschi2001introduction}.

The principal difficulty is that a perturbation of an inflow parameter changes not only the velocities of the newly injected particles, but also their later cell memberships, the local particle counts, and the collision schedules. A pathwise adjoint that differentiates one realized sequence of cell memberships and collision pairs treats these discrete choices as frozen and therefore does not, in general, approximate the derivative of the expected DSMC objective. Here, we instead differentiate the ensemble evolution. The resulting backward equation is the adjoint of the linearized Boltzmann collision operator, while the derivative of the inflow source is evaluated using a local score of the incoming flux distribution.

The ensemble collision formulation builds directly on the continuous adjoint framework introduced in our earlier work~\cite{caflisch2021adjoint}. In that work, the nonlinear collision operator was linearized around a general nonequilibrium distribution, its adjoint was derived in a weighted $L^2$ space, and the continuous adjoint was identified with a Fr\'echet derivative of the objective. The earlier derivation was presented for the spatially homogeneous Boltzmann equation. The present formulation applies the same linearized collision-adjoint structure locally in every spatial cell and couples it to free transport, absorbing outflow, and parameter-dependent inflow sources. In this sense, the new method extends the continuous, or optimize-then-discretize, side of~\cite{caflisch2021adjoint}, while using DSMC particles to perform the forward simulation and the Monte Carlo quadrature of the backward collision integrals.

For an inflow distribution $g(v,t;\theta)$ depending on boundary parameters $\theta$, the incoming particle flux is
\[
(-v\cdot\mathbf n)g(v,t;\theta)
\mathds{1}_{\{v\cdot\mathbf n<0\}},
\]
where $\mathbf n$ is the outward unit normal. Thus, over a boundary segment of area $\Delta S$ and a time interval $\Delta t$, the expected number of particles with velocities in $[v,v+dv]$ is
\[
dN =
\Delta S\,\Delta t\,(-v\cdot\mathbf n)g(v,t;\theta)
\mathds{1}_{\{v\cdot\mathbf n<0\}}\,dv.
\]
In the one-dimensional spatial setting, we take $\Delta S=1$. At each time step, let $N_L$ and $N_R$ denote the numbers of particles injected at the left and right boundaries, respectively. In the formulation considered below, $N_L$ and $N_R$ are evaluated at the boundary parameters and then held fixed during differentiation. Consequently, we differentiate the normalized distributions of the incoming particle states, but not the discontinuous dependence of the integer-valued injection counts on the boundary parameters.

Let $\mathfrak L_k$ and $\mathfrak R_k$ denote the particles injected during the interval $[t_k,t_{k+1}]$ through the left and right boundaries. At the left boundary, an incoming velocity $g_{k,i}$ and an entry time $\xi_{k,i}\sim\mathcal U([0,\Delta t])$ are sampled independently, and the particle state at time $t_{k+1}$ is
\begin{equation}
v_{k+1,i}=g_{k,i},
\qquad
x_{k+1,i}=L+\xi_{k,i}\,g_{k,i}\cdot e_1,
\qquad i\in\mathfrak L_k.
\label{eq:ensemble_inflow_left_state}
\end{equation}
Similarly, at the right boundary,
\begin{equation}
v_{k+1,i}=h_{k,i},
\qquad
x_{k+1,i}=R+\xi_{k,i}\,h_{k,i}\cdot e_1,
\qquad i\in\mathfrak R_k,
\label{eq:ensemble_inflow_right_state}
\end{equation}
where $h_{k,i}\cdot e_1<0$. Particles transported outside $[L,R]$ are removed from the active forward population and make no subsequent contribution to the collision dynamics or the terminal objective.

\subsection{Forward ensemble and linearized collision operator}

Let $f_k^j(v):=f(t_k,x\in\Omega_j,v)$ denote the particle distribution in spatial cell $\Omega_j$ at time $t_k$, and let
\[
N_k^j=|\mathcal V_k^j|,
\qquad
\rho_k^j=\frac{N_k^j}{N_0|\Omega_j|},
\qquad
\mu_k^j=\rho_k^j\int_{\mathbb S^2}q(\sigma)\,d\sigma.
\]
The forward DSMC solver uses the actual full-population value of $N_k^j$ to determine the local collision frequency. The expected number of collision pairs in cell $j$ during one time step is
\[
a_k^j=\frac{N_k^j\Delta t\,\mu_k^j}{2}.
\]
An integer realization of $a_k^j$ is obtained using the rounding rule of the forward DSMC algorithm, and the collision particles are selected uniformly without replacement from $\mathcal V_k^j$.

To derive the ensemble adjoint, we perturb the distribution according to $f_k\mapsto f_k+\delta f_k$. Equivalently, wherever $f_k>0$, we may write the multiplicative perturbation $\delta f_k=f_k\psi_k$. Following the notation of~\cite{caflisch2021adjoint}, the linearized operator associated with the multiplicative perturbation is
\begin{equation}
L[f_k]\psi_k
=
f_k^{-1}
\left[
Q(f_k,f_k\psi_k)+Q(f_k\psi_k,f_k)
\right].
\label{eq:multiplicative_linearized_collision_inflow}
\end{equation}
Since the Boltzmann collision operator is quadratic, the equivalent additive first variation is
\begin{equation}
\delta Q(f_k,f_k)
=Q(\delta f_k,f_k)+Q(f_k,\delta f_k)
=:\mathcal L_{f_k}\delta f_k.
\label{eq:linearized_collision_inflow}
\end{equation}
The two notations are related by
\[
\mathcal L_{f_k}\delta f_k=f_k\,L[f_k]\psi_k,
\qquad \delta f_k=f_k\psi_k.
\]
Both terms in~\eqref{eq:linearized_collision_inflow} are required. The first describes the collision of a perturbed target distribution with the unperturbed forward distribution, while the second describes the change in the collision-partner distribution and in the local collision rate. In particular, the dependence of the ensemble collision dynamics on the local particle count $N_k^j$ is contained in this linearization. This is the contribution that is absent when one differentiates a single DSMC trajectory while holding its cell memberships and collision schedule fixed.

\subsection{The ensemble adjoint equation}

Consider a terminal objective of the form
\begin{equation}
J
=
\frac{1}{\mathcal N}
\sum_{i\in\mathcal I_M}r(x_{M,i},v_{M,i}),
\qquad
\mathcal N=N_0+M(N_L+N_R),
\label{eq:ensemble_inflow_objective}
\end{equation}
where $\mathcal I_M$ is the set of particles active at the terminal time $T=M\Delta t$. We introduce the scalar ensemble adjoint
\[
\lambda_k(x,v)=\frac{\delta J}{\delta f_k}(x,v).
\]
The sign convention differs from the one adopted in our earlier work~\cite{caflisch2021adjoint}. This convention gives the gradient formula in
\eqref{eq:ensemble_inflow_gradient_left}-\eqref{eq:ensemble_inflow_gradient_right} below.
The terminal condition is
\begin{equation}
\lambda_M(x,v)=\frac{r(x,v)}{\mathcal N},
\qquad (x,v)\in\Omega\times\mathbb R^3,
\label{eq:ensemble_inflow_terminal_adjoint}
\end{equation}
and $\lambda_M=0$ outside $\Omega$.

We first pull the adjoint backward through free transport and absorbing outflow. Define
\begin{equation}
\widetilde\lambda_{k+1}(x,v)
=
\mathds{1}_{\{x+\Delta t\,v\cdot e_1\in\Omega\}}
\lambda_{k+1}(x+\Delta t\,v\cdot e_1,v).
\label{eq:ensemble_inflow_transport_adjoint}
\end{equation}
The indicator in~\eqref{eq:ensemble_inflow_transport_adjoint} sets the future value to zero when the corresponding particle leaves the domain.

For a collision between velocities $v$ and $w$, let $(v',w')$ denote the post-collision velocities. The adjoint of the linearized collision operator in cell $j$ is
\begin{align}
(\mathcal L_{f_k^j}^{*}\psi)(v)
&=
\int_{\mathbb R^3}f_k^j(w)
\int_{\mathbb S^2}q(\sigma)
\Big[
\psi(v')+\psi(w')-\psi(v)-\psi(w)
\Big]\,d\sigma\,dw.
\label{eq:ensemble_collision_adjoint_inflow}
\end{align}
This is the same four-term gain-loss structure as the operator $L^*[f]$ derived in~\cite{caflisch2021adjoint}. The only change is that it is applied to the local distribution $f_k^j$ in each spatial cell. In particular, when $\delta f_k=f_k\psi_k$, the duality relation from the earlier work becomes
\[
\bigl(L[f_k]\psi_k,\lambda_k\bigr)_{L^2(f_k\,dv)}
=
\bigl(\psi_k,L^*[f_k]\lambda_k\bigr)_{L^2(f_k\,dv)},
\]
which is the continuous identity approximated by the collision-partner Monte Carlo sum below. The discrete collision-transport splitting used by the DSMC implementation is then approximated by
\begin{equation}
\lambda_k(x,v)
=
\widetilde\lambda_{k+1}(x,v)
+\Delta t\,
(\mathcal L_{f_k^j}^{*}\widetilde\lambda_{k+1})(x,v),
\qquad x\in\Omega_j.
\label{eq:ensemble_adjoint_recursion_inflow}
\end{equation}
Equation~\eqref{eq:ensemble_adjoint_recursion_inflow} differentiates the expected collision evolution rather than the particular collision pairs selected in one forward realization. Thus, no derivative of an individual pair-selection index or hard cell-membership indicator appears in the backward equation.

\subsection{Local score for the inflow boundary source}

We next derive the derivative with respect to the inflow parameters. Let $b_L(v;\theta_L)$ and $b_R(v;\theta_R)$ denote the normalized left and right incoming flux densities. For a sufficiently regular function $F$, the score identity gives
\begin{equation}
\partial_{\theta_s}
\mathbb E_{b_s(\cdot;\theta_s)}[F(V)]
=
\mathbb E_{b_s(\cdot;\theta_s)}
\left[
F(V)\,\partial_{\theta_s}\log b_s(V;\theta_s)
\right],
\qquad s\in\{L,R\}.
\label{eq:local_score_identity_inflow}
\end{equation}

A particle injected during $[t_k,t_{k+1}]$ is already part of the particle distribution at time $t_{k+1}$. Its adjoint value therefore contains its complete future effect, including all later transport, collisions, and its possible terminal contribution. Applying~\eqref{eq:local_score_identity_inflow} locally to each boundary source gives
\begin{align}
\nabla_{\theta_L}J
&=
\sum_{k=0}^{M-1}N_L\,
\mathbb E\left[
\lambda_{k+1}
\bigl(L+\xi V\cdot e_1,V\bigr)
\nabla_{\theta_L}\log b_L(V;\theta_L)
\right],
\label{eq:ensemble_inflow_gradient_left}
\\
\nabla_{\theta_R}J
&=
\sum_{k=0}^{M-1}N_R\,
\mathbb E\left[
\lambda_{k+1}
\bigl(R+\xi V\cdot e_1,V\bigr)
\nabla_{\theta_R}\log b_R(V;\theta_R)
\right],
\label{eq:ensemble_inflow_gradient_right}
\end{align}
where $\xi\sim\mathcal U([0,\Delta t])$ is independent of $V$. There is no additional score associated with $\xi$, because its distribution is parameter-independent. The sign in
\eqref{eq:ensemble_inflow_gradient_left}-\eqref{eq:ensemble_inflow_gradient_right} is positive under the adjoint
definition~\eqref{eq:ensemble_inflow_terminal_adjoint}.

For the anisotropic zero-mean Maxwellian used in the numerical experiments, let $T_s=(T_{s,1},T_{s,2},T_{s,3})$. The normalized incoming flux density at the left boundary is
\begin{align}
b_L(v;T_L)
&=
\mathds{1}_{\{v_1>0\}}
\frac{v_1}{T_{L,1}}
\exp\left(-\frac{v_1^2}{2T_{L,1}}\right)
\prod_{\ell=2}^3
\frac{1}{\sqrt{2\pi T_{L,\ell}}}
\exp\left(-\frac{v_\ell^2}{2T_{L,\ell}}\right).
\label{eq:left_half_range_flux_density}
\end{align}
At the right boundary, $v_1>0$ is replaced by $v_1<0$ and $v_1$ by $|v_1|$. The componentwise temperature scores are
\begin{align}
\partial_{T_{s,1}}\log b_s(v;T_s)
&=
-\frac{1}{T_{s,1}}
+\frac{v_1^2}{2T_{s,1}^2},
\label{eq:normal_temperature_score_inflow}
\\
\partial_{T_{s,\ell}}\log b_s(v;T_s)
&=
-\frac{1}{2T_{s,\ell}}
+\frac{v_\ell^2}{2T_{s,\ell}^2},
\qquad \ell=2,3,
\label{eq:tangential_temperature_score_inflow}
\end{align}
for $s\in\{L,R\}$. Since the normal score depends on $v_1^2$, the same formula applies at both boundaries; the boundary orientation is represented by the support of the corresponding half-range distribution.

If the physical injection intensity is also differentiated, then its derivative must be added separately. For example, for a zero-mean Maxwellian inflow with intensity $N_s(T_{s,1})\propto\sqrt{T_{s,1}}$, the smooth intensity contribution adds
\(
\partial_{T_{s,1}}\log N_s(T_{s,1})
=\frac{1}{2T_{s,1}}
\)
to~\eqref{eq:normal_temperature_score_inflow}. A parameter-dependent normalization $\mathcal N$ would introduce the corresponding quotient-rule term. These are not included in the formulation considered in
\eqref{eq:ensemble_inflow_gradient_left}-\eqref{eq:ensemble_inflow_gradient_right}.

\subsection{Monte Carlo approximation of the backward equation}

For a moderate particle population, the adjoint field may be represented on the complete sequence of forward particle populations, but sometimes, storing every particle at every time level is impractical. We therefore introduce a memory-bounded adjoint formulation without modifying the physical forward simulation.

At every time level and in every cell $\Omega_j$, we retain a uniform sample $\widehat{\mathcal V}_k^j\subset\mathcal V_k^j$ containing at most $N_{\mathrm{res}}$ particles. The exact full-population cell count $N_k^j$ is stored separately and is used in the collision coefficient. We also retain uniform samples $\widehat{\mathfrak L}_k\subset\mathfrak L_k$ and $\widehat{\mathfrak R}_k\subset\mathfrak R_k$ of the injected particles. If $m_{L,k}=|\widehat{\mathfrak L}_k|$, the corresponding weight is $N_L/m_{L,k}$, with an analogous definition on the right.

The adjoint is represented by scalar values on the retained states. A local interpolant $\mathcal I_k[\lambda_k](x,v)$ supplies the adjoint at other non-retained states, including the post-collision states and the injected source states. In the numerical implementation, the interpolant is a weighted local affine fit in $(x,v_1,v_2,v_3)$ based on neighboring retained particles.

For a retained state $(x_{k,i},v_{k,i})\in\Omega_j$, we sample $N_{\mathrm{col}}$ collision-partner states $(x_{k,i_s},w_{k,i}^{(s)})$ from the same cell and scattering directions $\sigma_{k,i}^{(s)}$ from the normalized angular kernel. Let $(v_{k,i}^{\prime(s)},w_{k,i}^{\prime(s)})$ denote the corresponding post-collision velocities, and let
\(
Q_0=\int_{\mathbb S^2}q(\sigma)\,d\sigma
\). The collision-adjoint correction is approximated by
\begin{align}
\widehat C_{k,i}
&=
\Delta t\,\rho_k^j Q_0\,
\frac{1}{N_{\mathrm{col}}}
\sum_{s=1}^{N_{\mathrm{col}}}
\Big[
\widetilde\lambda_{k+1}
(x_{k,i},v_{k,i}^{\prime(s)})
+\widetilde\lambda_{k+1}
(x_{k,i_s},w_{k,i}^{\prime(s)})
\nonumber\\
&\hspace{5.2cm}
-\widetilde\lambda_{k+1}(x_{k,i},v_{k,i})
-\widetilde\lambda_{k+1}(x_{k,i_s},w_{k,i}^{(s)})
\Big].
\label{eq:MC_collision_adjoint_inflow}
\end{align}
The adjoint update for this retained particle is therefore
\begin{equation}
\widehat\lambda_{k,i}
=
\widetilde\lambda_{k+1}(x_{k,i},v_{k,i})
+\widehat C_{k,i}.
\label{eq:reservoir_adjoint_update_inflow}
\end{equation}
Although the velocity distribution inside a cell is represented by the retained states, the coefficient $\rho_k^j$ in
\eqref{eq:MC_collision_adjoint_inflow} is always computed using the exact
full-population count $N_k^j$. Thus, using a subset of the states does not artificially reduce the local collision frequency.

The sampled left-boundary gradient is
\begin{align}
\widehat{\nabla_{\theta_L}J}
&=
\sum_{k=0}^{M-1}
\frac{N_L}{m_{L,k}}
\sum_{i\in\widehat{\mathfrak L}_k}
\mathcal I_{k+1}[\widehat\lambda_{k+1}]
(x_{k+1,i},v_{k+1,i})
\nabla_{\theta_L}\log b_L(v_{k+1,i};\theta_L),
\label{eq:sampled_left_score_gradient_inflow}
\end{align}
and the right-boundary gradient is obtained analogously. The expansion weights make these source sums unbiased conditional on the reconstructed adjoint field.

The complete procedure is summarized below.
\begin{enumerate}
\item Run the full-population forward DSMC simulation. At every time level, store a subset of the particles, the exact cell
counts, and weighted samples of the two inflow sources.
\item Initialize the terminal adjoint using
\eqref{eq:ensemble_inflow_terminal_adjoint}.
\item For $k=M-1,\ldots,0$, construct the interpolant of the adjoint at
$t_{k+1}$, accumulate the two local boundary-score contributions, pull the
adjoint backward through transport and absorbing outflow~\eqref{eq:ensemble_inflow_transport_adjoint}, and apply the
Monte Carlo collision correction~\eqref{eq:MC_collision_adjoint_inflow}.
\item Sum the boundary-source contributions over all time steps to obtain the
left and right parameter gradients; see~\eqref{eq:sampled_left_score_gradient_inflow}.
\end{enumerate}

For a scalar objective, only one scalar backward field is required. After one backward adjoint solve, the adjoint is evaluated once at each sampled incoming state. The gradient with respect to each boundary parameter is then obtained by multiplying these adjoint values by the corresponding source-score components and summing over the incoming samples. Consequently, the dominant forward- and backward-DSMC costs are independent of the number of inflow parameters.

\section{Numerical Examples} \label{sec:numerics}
All numerical examples in this section are conducted on a 1D spatial domain $\Omega = [0, 1] \subset \mathbb{R}$ and a 3D velocity domain $\mathbb{R}^3$. {The Knudsen number in all examples is chosen to be $1$.} We consider four numerical examples and verify the adjoint gradients by comparison with centered finite-difference approximations. The examples include two-sided thermal boundaries, mixed thermal-specular boundaries, two-sided inflow, and a high-Mach-number Couette flow driven by moving thermal walls.

\subsection{Heat conduction: two-sided thermal BC}\label{subsec:test1}
Our first example concerns heat conduction. We enforce thermal BC at the left boundary $x=0$ and the right boundary point $x=1$, both with half-Maxwellian flux distributions. The temperatures for the two half-Maxwellian fluxes are respectively
\begin{equation}\label{eq:flux_temp}
T_L = [T_L^1,T_L^2,T_L^3 ] =  [0.6,0.6,0.6]\,,\qquad T_R= [T_R^1,T_R^2,T_R^3 ]=  [0.9,0.9,0.9]\,.
\end{equation}
The temperature difference induces a non-equilibrium heat flux across the domain. This setup provides a canonical test of thermal BC in kinetic simulations and allows us to assess whether the adjoint method correctly captures the sensitivity of energy transport with respect to boundary temperature parameters.

The initial distribution for the Boltzmann equation is set as
{
$$
f_0(x,v) = \left(\frac{1}{\sqrt{2\pi}\sigma}\right)^3\exp\left(-\frac{|v|^2}{2\sigma^2}\right)\,,\quad \sigma = 0.7\,.
$$
}
Let $N$ denote the total number of particles used in the forward DSMC. We draw i.i.d.~samples $x_{0,i}\sim \mathcal{U}([0,1])$ and $v_{0,i}\sim \mathcal{N}(0,\sigma^2)$, $1\leq i \leq N$, to approximate the initial distribution $f_0$. The total forward simulation time is $T = 1$ with $\Delta t = 0.1$. We discretize the spatial domain $[0,1]$ into equal-size cells with $\Delta x \in \{0.05, 0.1, 0.2\}$. Here, we consider $T_L$ and $T_R$ as variables for an optimization problem with the objective function
\[
J_1(T_L,T_R) = \int_0^1 \int_{\bbR^3}  |v|^2 e^{-(x-0.2)^2} f(x,v,T) dv dx\,,
\]
where $f(x,v,T)$ is the Boltzmann equation solution at the final simulation time $T = 1$.

The objective $J_1$ measures the localized kinetic energy near $x=0.2$, with the Gaussian weight emphasizing sensitivity to energy transport from the left boundary. Thus, the gradient $\nabla_{T_L,T_R} J_1$ quantifies how changes in boundary temperatures influence the spatial redistribution of energy inside the domain.

{
To illustrate the macroscopic state associated with the gradient test, Figure~\ref{fig:test1-macroscopic-profiles} shows the normalized number density $\widehat n=N_j/(N_0\Delta x)$, the scalar temperature $T=\langle |v-u|^2\rangle/3$, and the streamwise bulk velocity $u_1$ at the initial and final times.  The profiles are obtained from one forward realization with $N_0=2\times10^6$ particles and $\Delta x=0.05$.  Initially, the density is uniform, the temperature is approximately $0.49=0.7^2$, and the bulk velocity is zero up to sampling error.  At $T=1$, the density decreases from approximately $1.08$ near the left wall to $0.90$ near the right wall, while the temperature increases from approximately $0.61$ to $0.70$ across the domain in response to the wall temperatures $T_L=0.6$ and $T_R=0.9$.  The streamwise velocity remains small throughout the domain, with $|u_1|<9\times10^{-3}$, indicating that the dominant nonequilibrium effect is thermal transport.

\begin{figure}
    \centering
    \includegraphics[width=\textwidth]{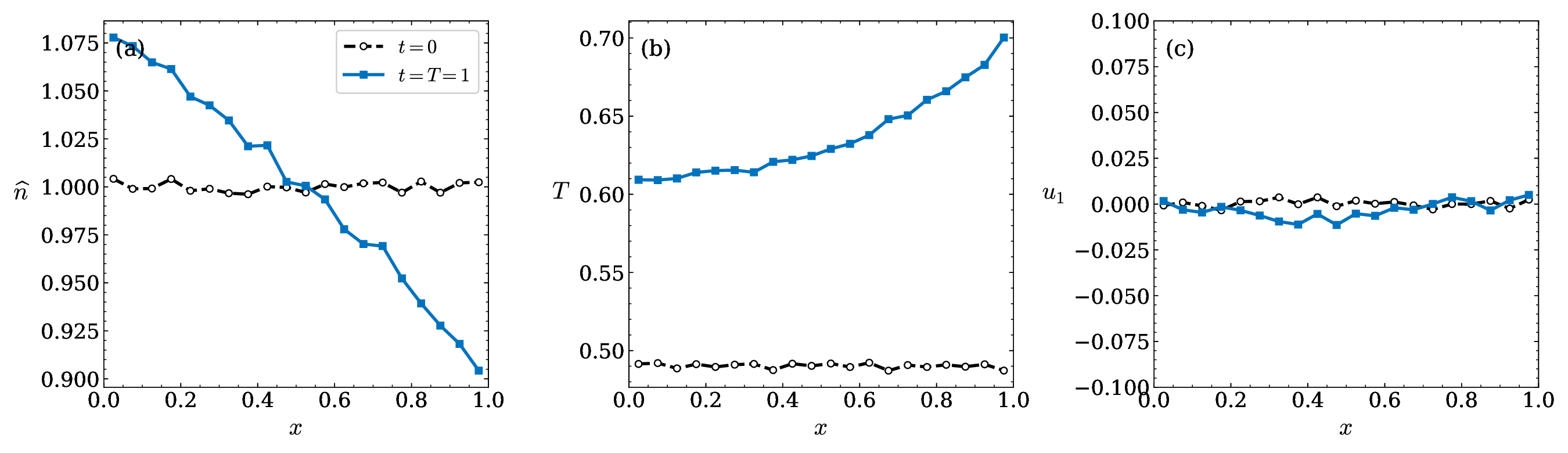}
    \caption{Two-sided thermal boundary conditions: profiles of
    (a) normalized number density $\widehat n$, (b) scalar temperature $T$,
    and (c) streamwise bulk velocity $u_1$ at $t=0$ and $t=T=1$.
    The simulation uses $N_0=2\times10^6$ and $\Delta x=0.05$.}
    \label{fig:test1-macroscopic-profiles}
\end{figure}
}

Following the algorithm derived in~Section~\ref{sec:adjoint_DSMC_thermal}, we first modify the forward DSMC algorithm by replacing the deterministic time step size $\Delta t$ with a random variable $\tau \sim \mathcal{N}(\Delta t, \varepsilon^2)$, where we set $\varepsilon = \Delta t / 10 = 0.01$. { This selective stochastic regularization smooths otherwise discontinuous boundary-interaction events while introducing only a small perturbation of the deterministic forward scheme for the chosen value of \(\varepsilon\).} To improve sampling efficiency, the random time step is applied only to particles that may interact with the boundary during the current update. Specifically, at each time level $t_k$, we examine the particle's current position $x_{k,i}$ and post-collision velocity $v_{k,i}'$. If
{
\[
\operatorname{dist}(x_{k,i}, \partial\Omega) < (\Delta t + 3\varepsilon) |\mathcal{P}v_{k,i}'|\,,
\]
}
then the particle is considered sufficiently close to the boundary that a boundary-crossing event could occur under a randomized update. In this case, the deterministic time step $\Delta t$ is replaced by a random variable $\tau \sim \mathcal{N}(\Delta t,\varepsilon^2)$. For all other particles, whose trajectories remain well separated from the boundary over the current time increment, the deterministic time step $\Delta t$ is retained. This selective stochastic regularization smooths otherwise discontinuous boundary interaction events, thereby enabling differentiation through boundary reflections while converging to the deterministic forward scheme as \(\varepsilon\to0\).

{
To assess the sensitivity of the regularized adjoint gradient to the regularization width, we repeat Test~1 with $N=10^6$, $\Delta x=0.1$, and $\varepsilon\in[0.005,0.02]$, using 48 independent runs for each value. Figure~\ref{fig:test1-epsilon-stability} shows the resulting mean objective and adjoint gradients with respect to $T_L$ and $T_R$, with error bars representing twice the standard error of the mean. The objective varies by less than $0.07\%$ over the tested interval, while the gradient components exhibit no systematic dependence on $\varepsilon$ within their Monte Carlo uncertainty. These results indicate that the regularized adjoint gradient is stable over the range of $\varepsilon$ relevant to the numerical experiments, including the chosen value $\varepsilon=\Delta t/10=0.01$.

\begin{figure}
    \centering
    \includegraphics[width=\textwidth]{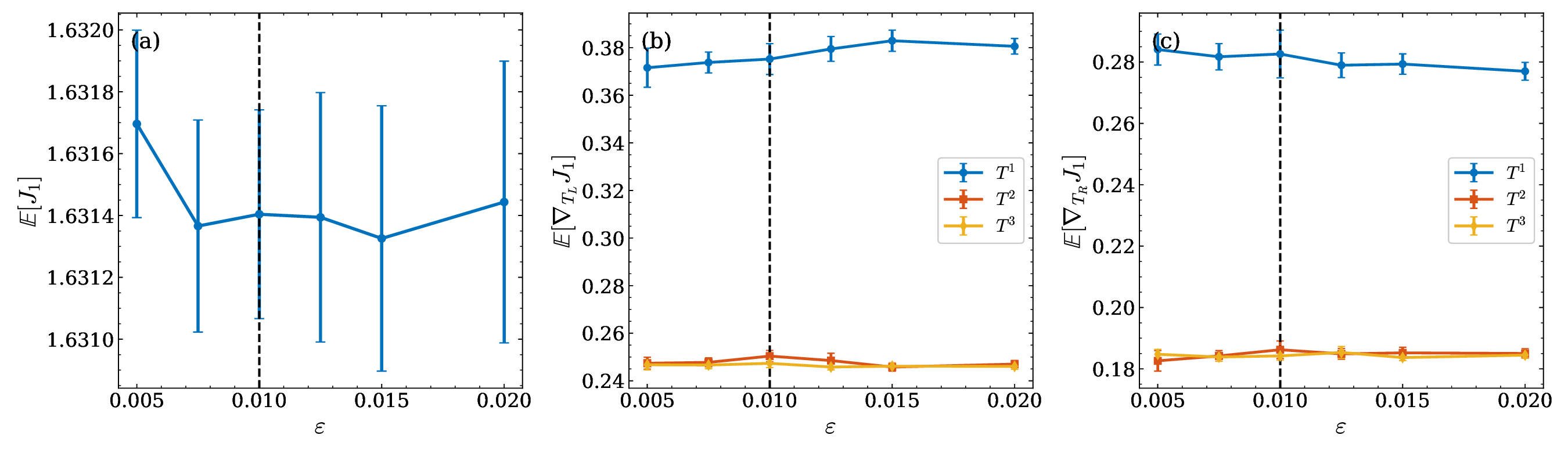}
    \caption{Dependence of (a) the mean objective $J_1$, (b) the mean adjoint
    gradient with respect to the left-wall temperature $T_L$, and (c) the mean
    adjoint gradient with respect to the right-wall temperature $T_R$ on the
    regularization width $\varepsilon$. Error bars represent twice the standard
    error of the mean over 48 independent runs. The black dashed line marks
    $\varepsilon=\Delta t/10=0.01$, used in the main numerical experiments.
    The study uses $N=10^6$ and $\Delta x=0.1$.}
    \label{fig:test1-epsilon-stability}
\end{figure}

}

Next, using this modified forward DSMC algorithm, we implement the adjoint DSMC approach to compute the gradient with respect to the boundary temperature parameters $T_L$ and $T_R$, denoted by $\nabla^\text{AJ}_{T_L} J_1$ and $\nabla^\text{AJ}_{T_R} J_1$, following the derived gradient formulae~\eqref{eq:thermal_grad_L}-\eqref{eq:thermal_grad_R}. To apply these equations, we compute the adjoint variable $\{\alpha_{k,i}\}$ and $\{\beta_{k,i}\}$ following the adjoint DSMC rule given in~\eqref{eq:alpha_adjoint} and~\eqref{eq:beta_adjoint}. The goal of the numerical experiment is therefore twofold: (i) to verify that the adjoint gradient agrees with a finite-difference benchmark, and (ii) to confirm that its variance exhibits the expected Monte Carlo scaling with respect to the number of particles $N$.

\begin{figure}
    \centering
    \subfloat[$|\nabla^\text{AJ}_{T_L} J_1 -\nabla^\text{FD}_{T_L} J_1|/|\nabla^\text{FD}_{T_L} J_1|$ for $\Delta x\in \{0.05, 0.1, 0.2\}$.]{
    \includegraphics[width = 0.33\textwidth]{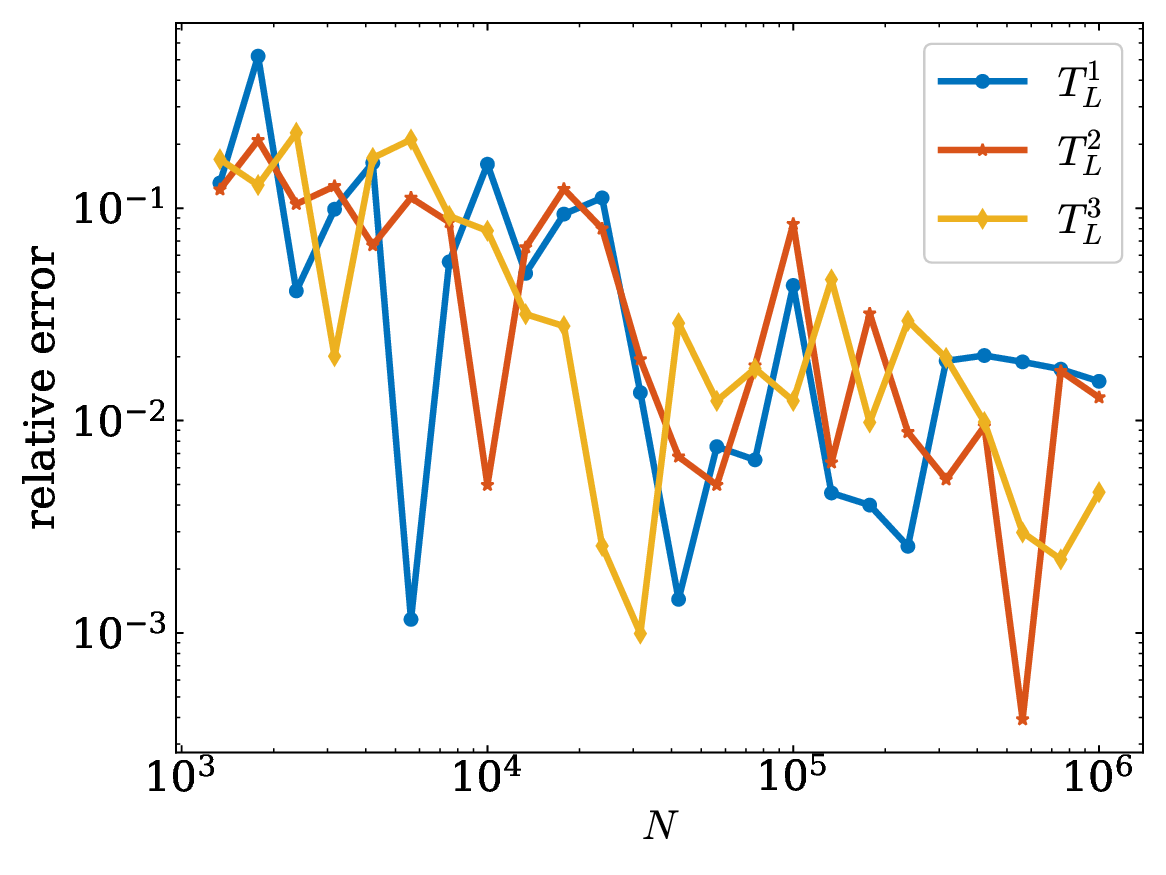}
    \includegraphics[width = 0.33\textwidth]{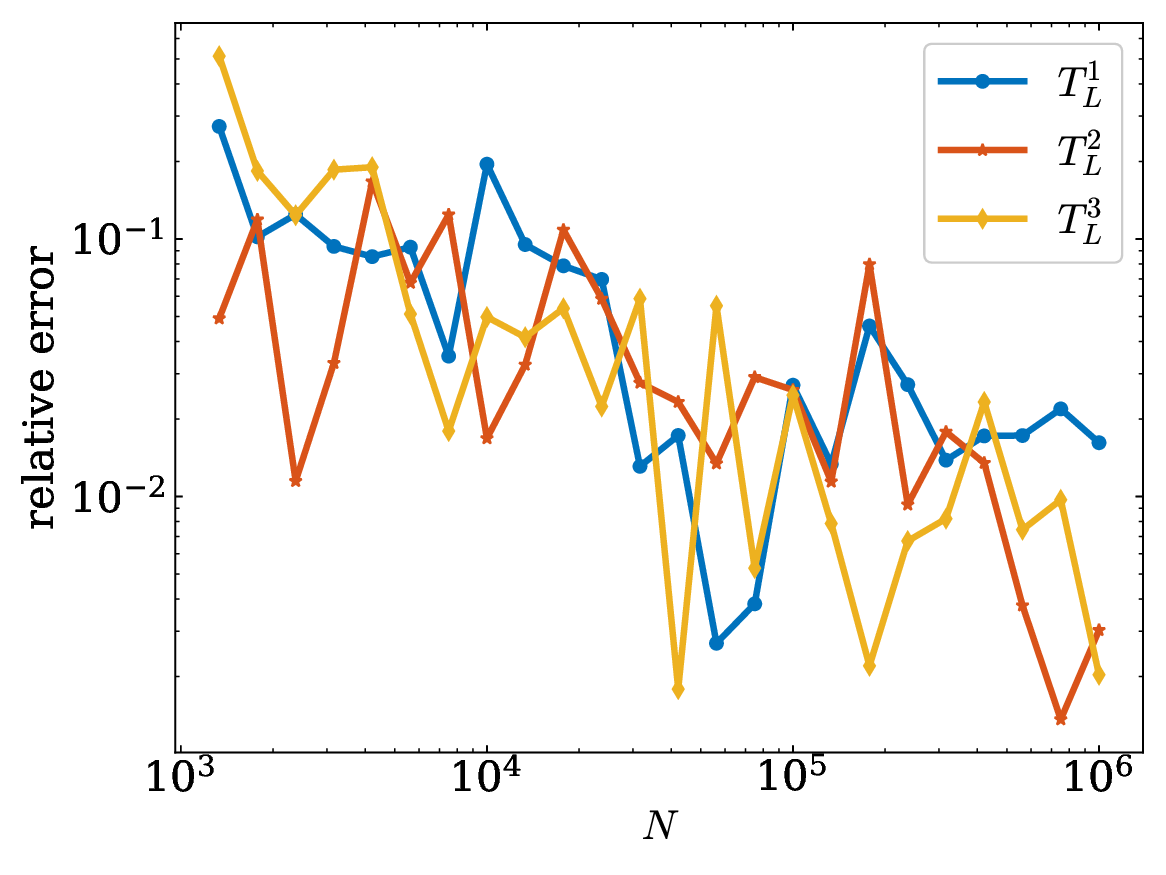}
    \includegraphics[width = 0.33\textwidth]{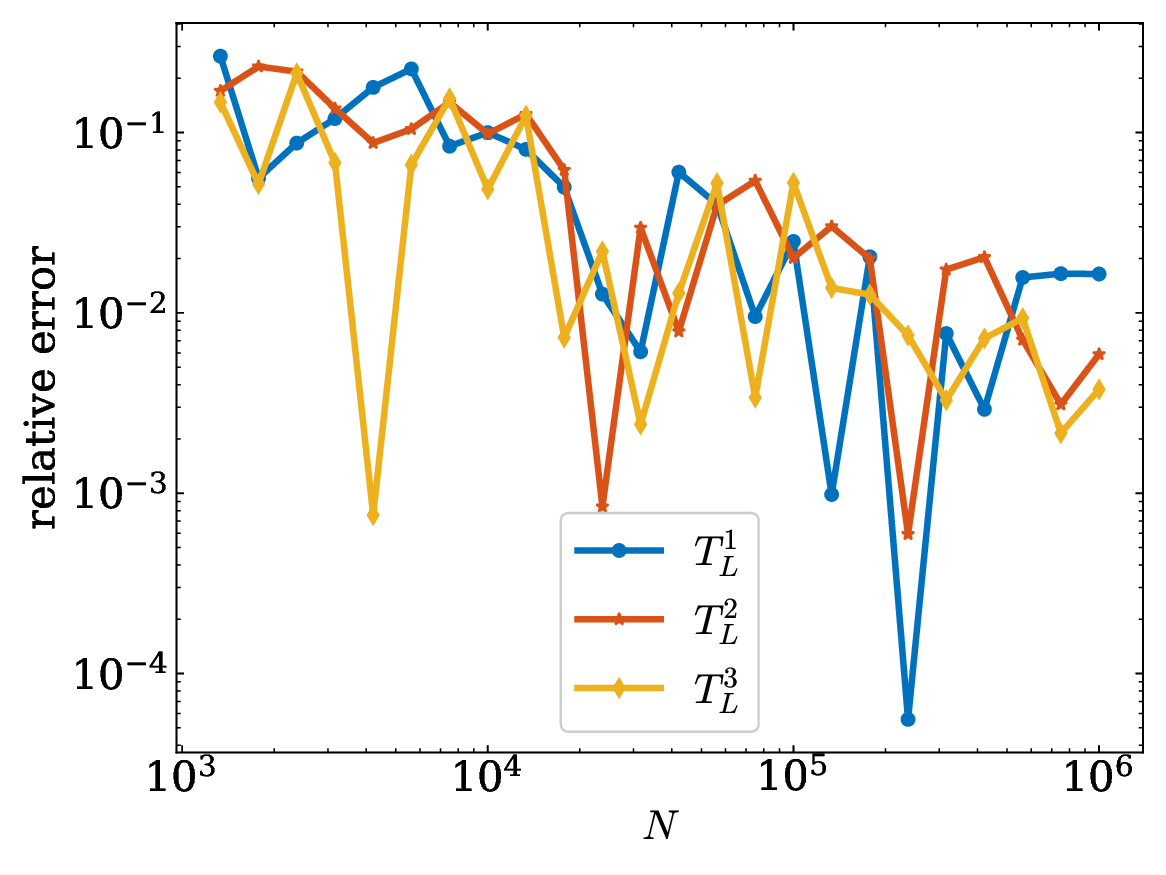}
    \label{fig:test1-TL-error}}\\
    \subfloat[$|\nabla^\text{AJ}_{T_R} J_1 -\nabla^\text{FD}_{T_R} J_1|/|\nabla^\text{FD}_{T_R} J_1|$ for $\Delta x\in \{0.05, 0.1, 0.2\}$.]{
    \includegraphics[width = 0.33\textwidth]{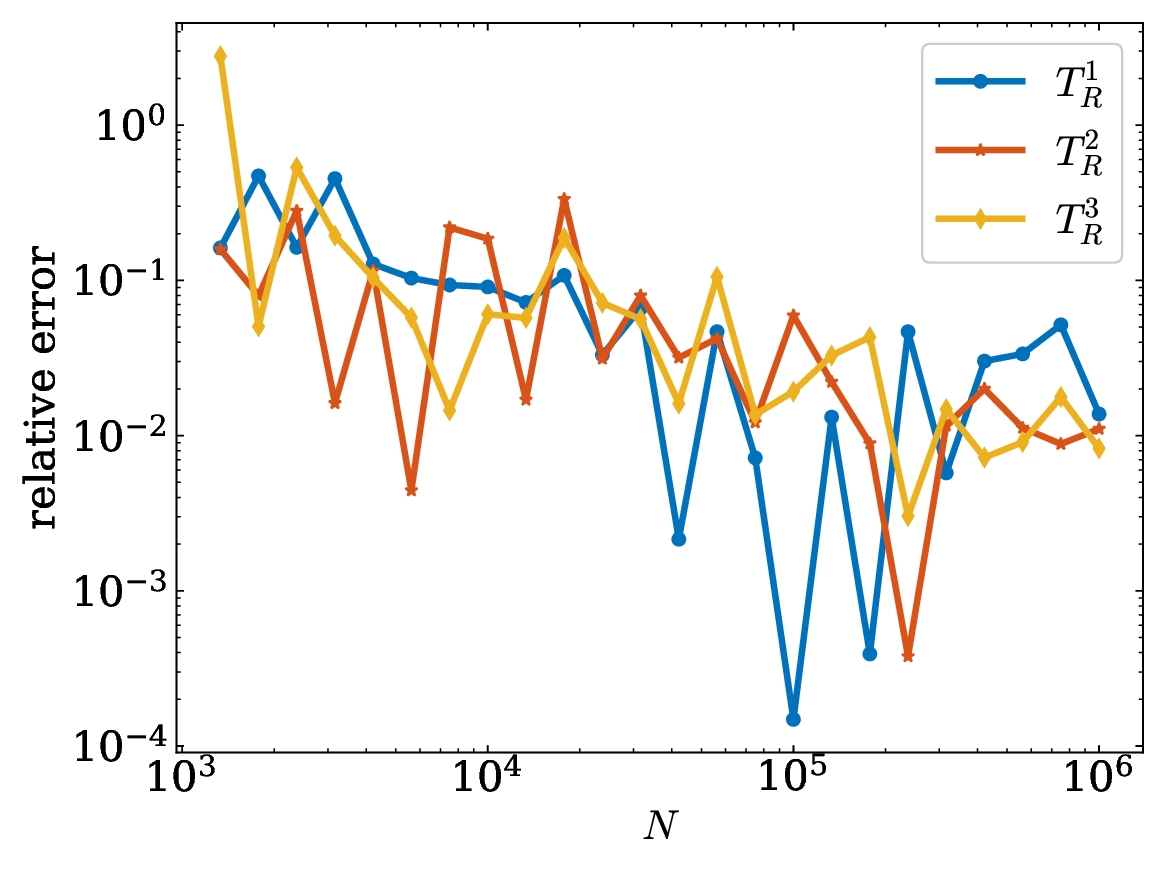}
    \includegraphics[width = 0.33\textwidth]{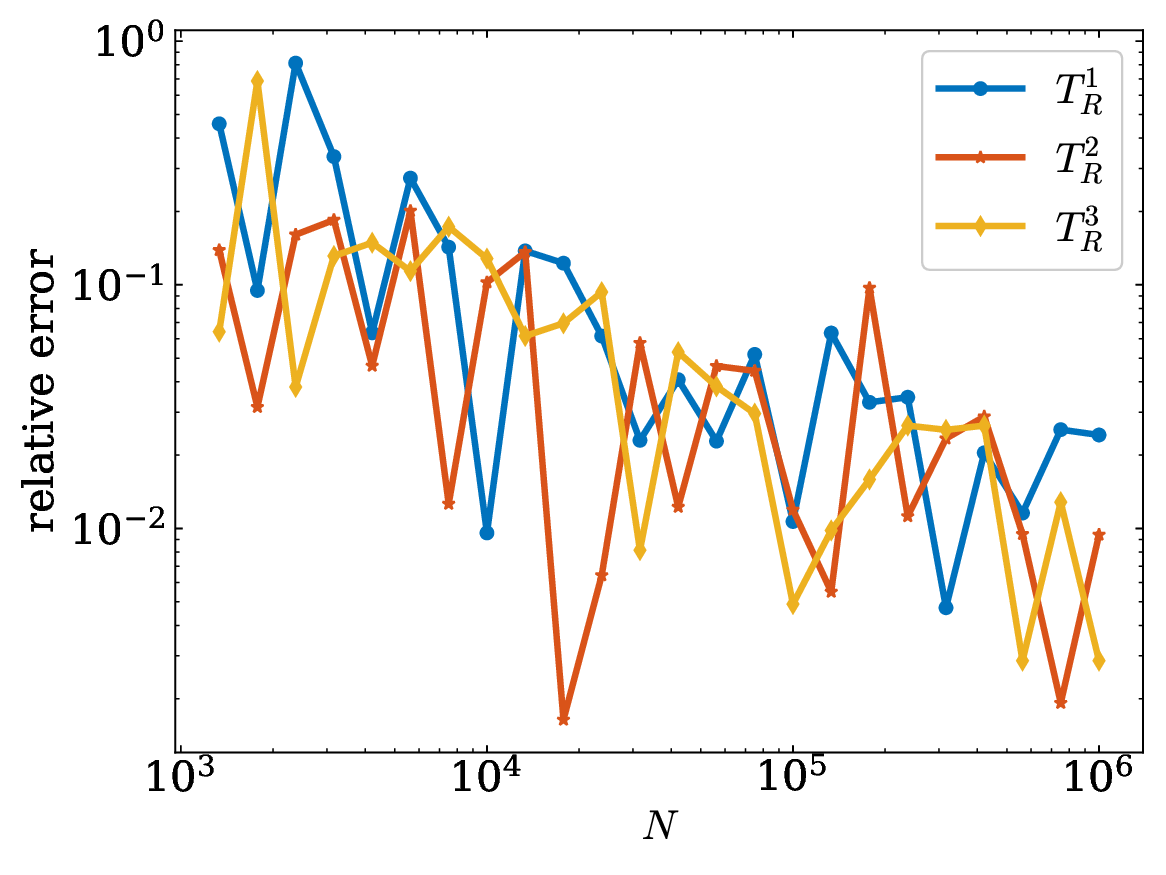}
    \includegraphics[width = 0.33\textwidth]{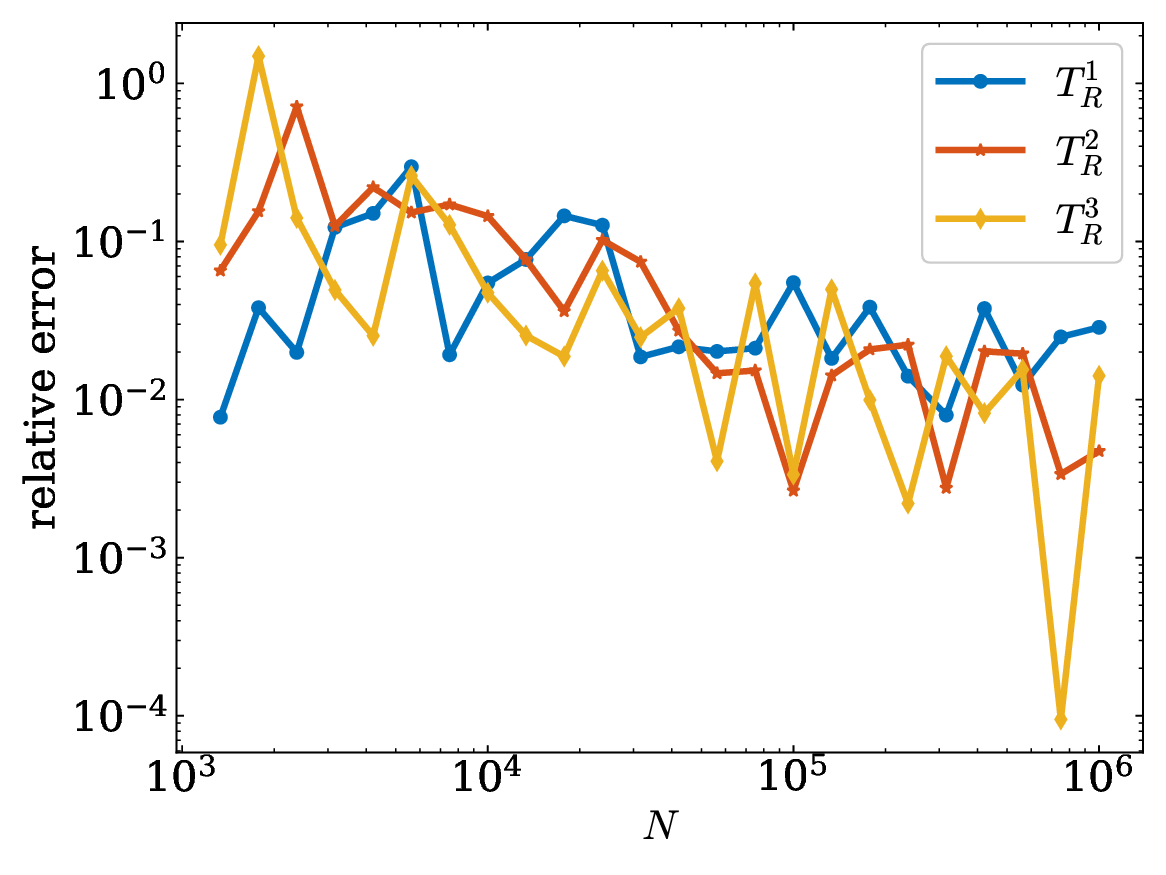}
    \label{fig:test1-TR-error}
    }\\
    \subfloat[standard deviation in $\nabla^\text{AJ}_{T_L} J_1$ for $\Delta x\in \{0.05, 0.1, 0.2\}$]{\includegraphics[width = 0.33\textwidth]{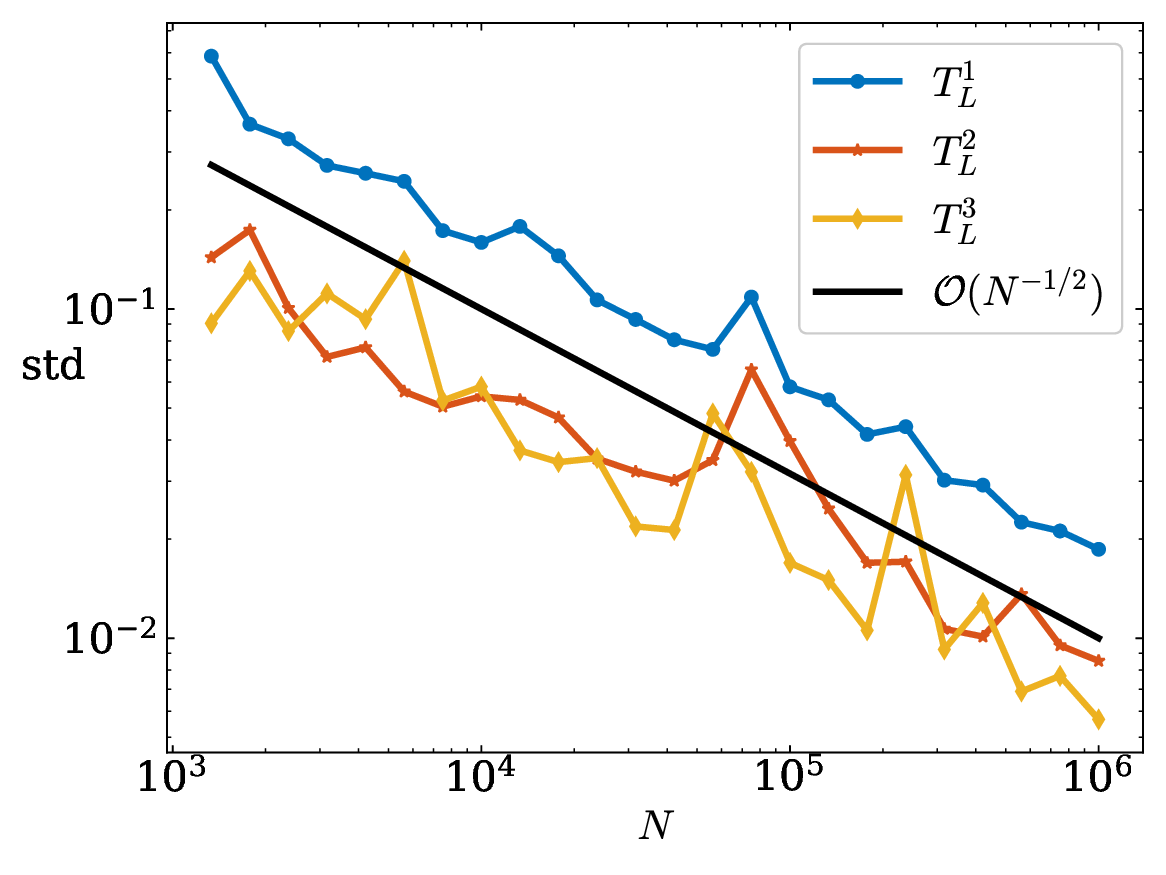}
    \includegraphics[width = 0.33\textwidth]{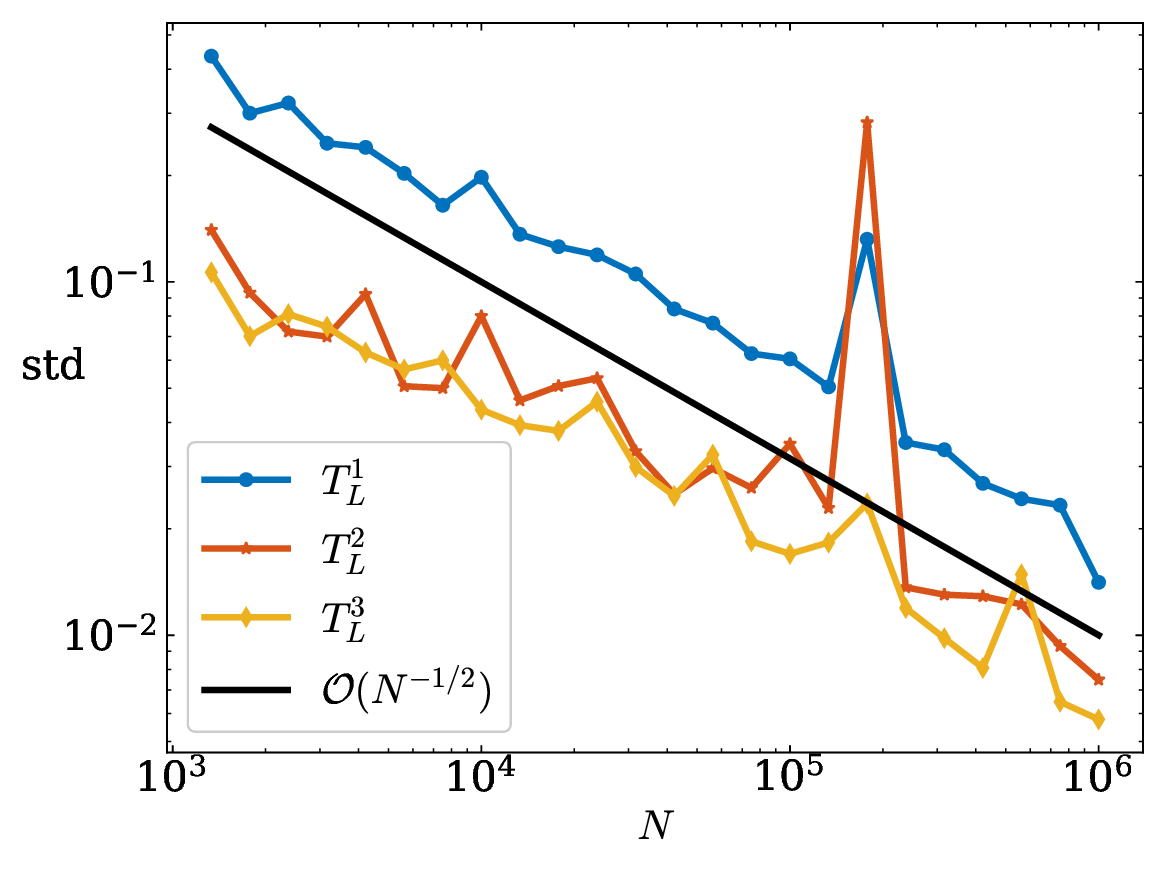}
    \includegraphics[width = 0.33\textwidth]{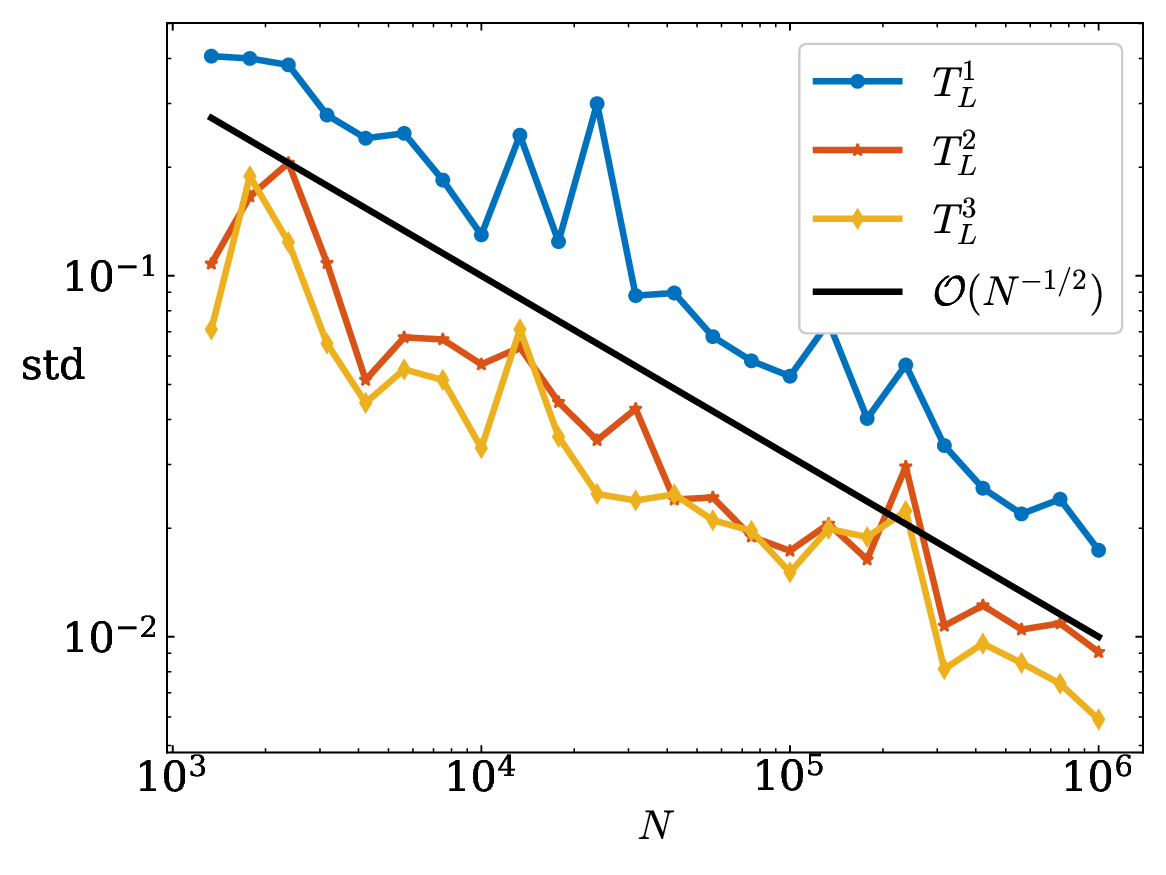}
    \label{fig:test1-TL-std}
    }\\
    \subfloat[standard deviation in $\nabla^\text{AJ}_{T_R} J_1$ for $\Delta x\in \{0.05, 0.1, 0.2\}$]{\includegraphics[width = 0.33\textwidth]{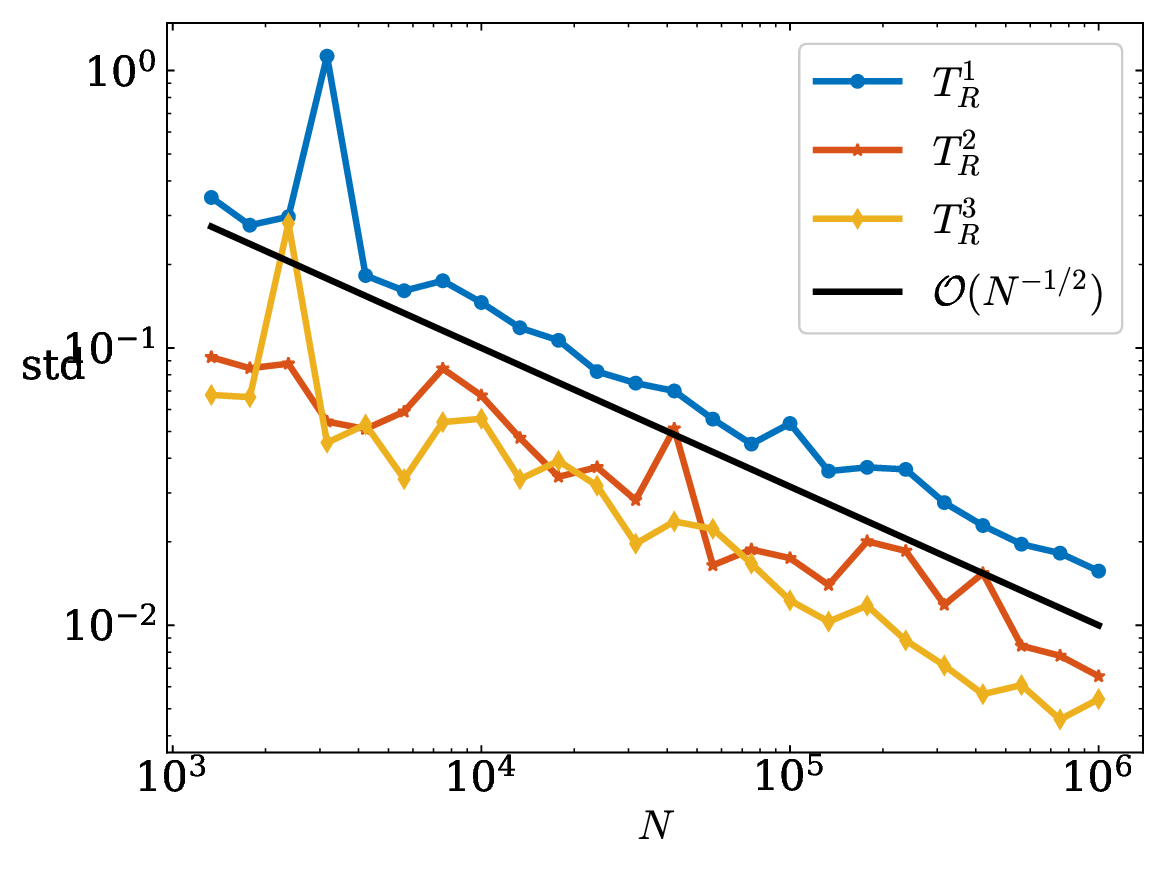}
    \includegraphics[width = 0.33\textwidth]{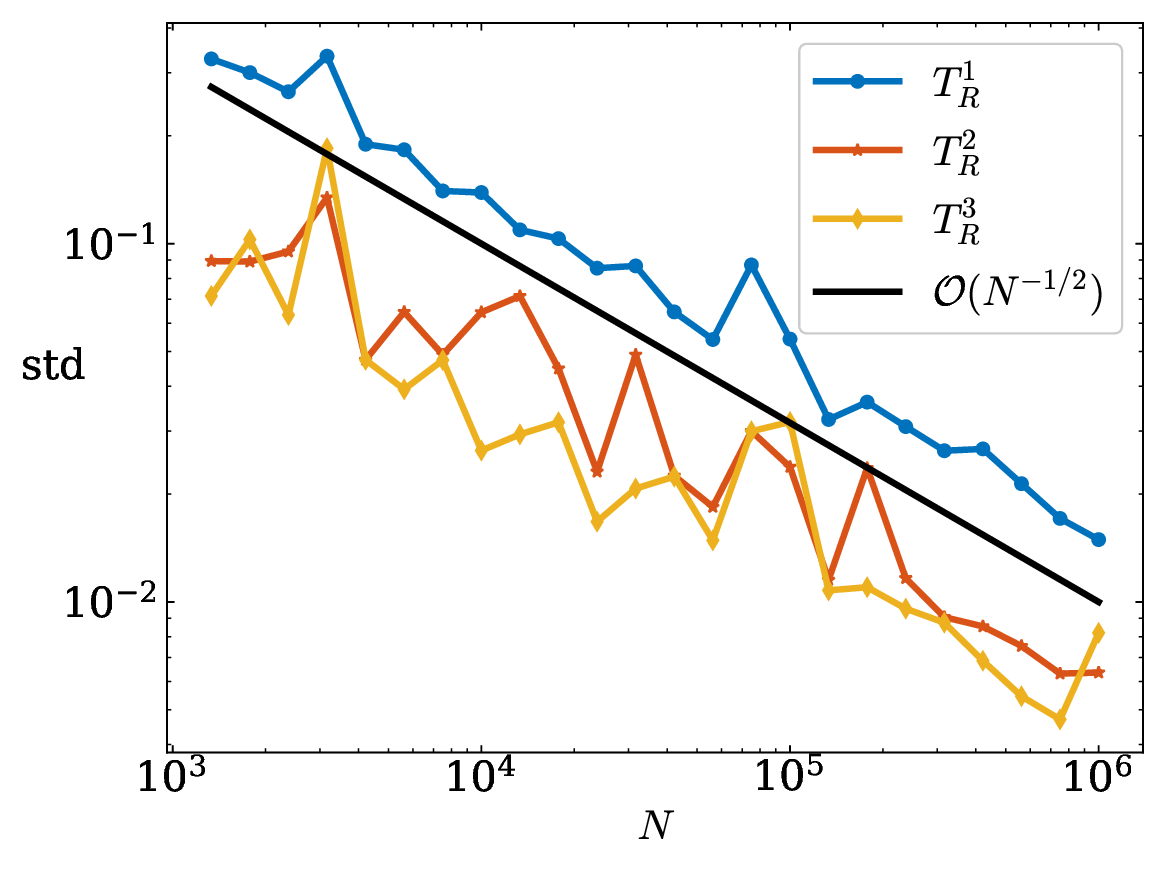}
    \includegraphics[width = 0.33\textwidth]{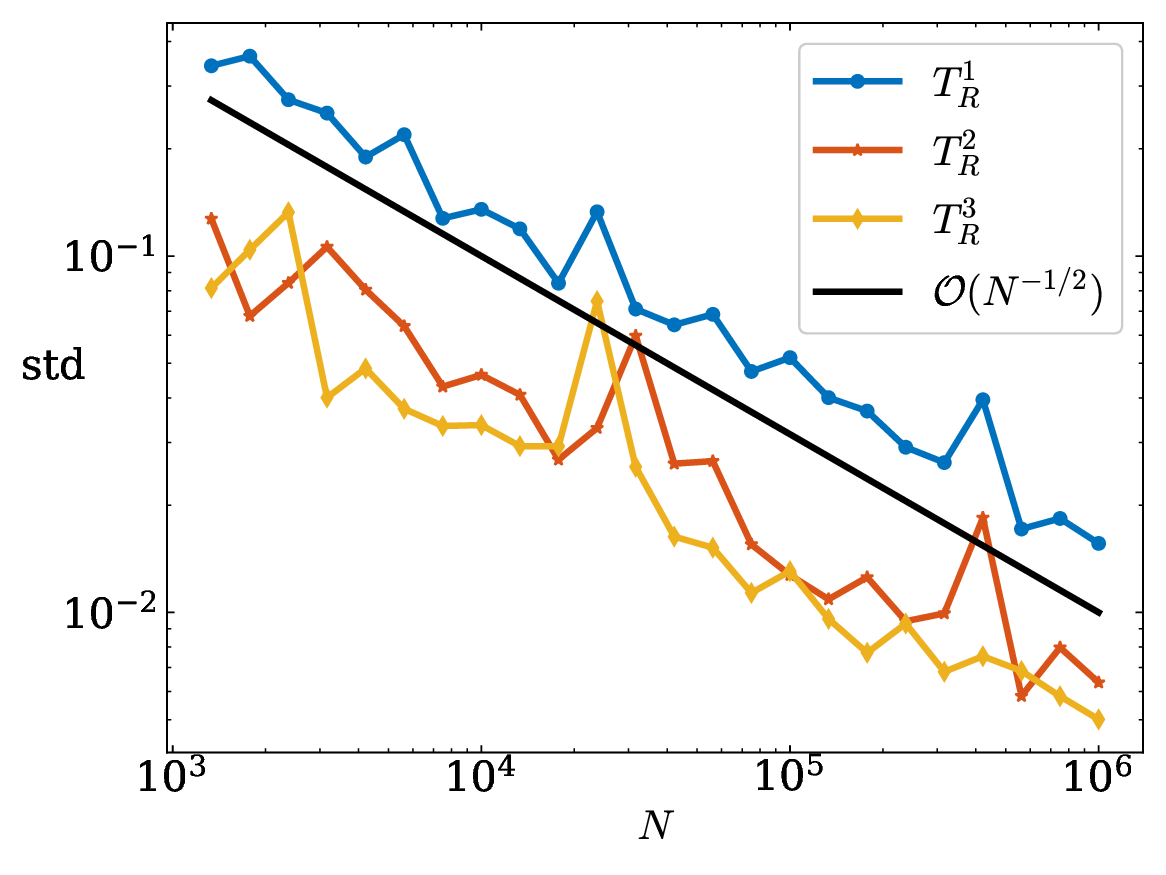}
    \label{fig:test1-TR-std}}
    \caption{Thermal-thermal
 BCs (\Cref{subsec:test1}): (a) and (b) show the relative error between the adjoint gradient and the finite-difference gradient for $T_L$ and $T_R$, respectively, as the number of particles $N$ increases; (c) and (d) show the standard deviation of the adjoint gradient for $T_L$ and $T_R$, respectively. The standard deviation decreases with \(N\) and is broadly consistent with the Monte Carlo \(N^{-1/2}\) reference trend. The three columns correspond to different spatial bin sizes: from left to right, $\Delta x = 0.05$, $0.1$, and $0.2$, respectively.
 }
    \label{fig:test1}
\end{figure}

To verify the accuracy of the computed adjoint gradients, we compare them to gradients approximated using finite difference perturbations. For instance, to approximate $\partial_{T_L^1} J_1$, we employ a central difference scheme given by:
\[
\nabla^\text{FD}_{T_L^1} J_1 := \frac{J_1(T_L^1 + \delta) - J_1(T_L^1 - \delta)}{2\delta} \approx \partial_{T_L^1} J_1\,,
\]
where $\delta = 0.025$ is fixed for all components of $T_L$ and $T_R$. The finite-difference gradients are computed using the original DSMC scheme with deterministic time step $\Delta t$, whereas the adjoint gradient corresponds to the $\varepsilon$-regularized scheme. Thus, the comparison is between an adjoint approximation of $\nabla J_\varepsilon$ and a finite-difference approximation of $\nabla J_0$. Agreement between the two estimators indicates that the regularization bias is small for the chosen value of $\varepsilon$, while the separate $\varepsilon$-stability experiment examines this dependence directly.

We range the number of particles $N$ from $10^{3}$ to $10^6$ to examine the accuracy of the adjoint gradient. To further reduce random error, we conduct $96$ independent runs for both the adjoint gradient and the finite-difference gradient and compare the mean value of these experimental outcomes. We also use these $96$ runs to approximate the standard deviation. In the heat conduction example, we analyze the performance of the adjoint gradient computation in comparison to the finite-difference gradient approximation as the number of particles $N$ increases. Figures~\ref{fig:test1-TL-error}-\ref{fig:test1-TR-error} illustrate the relative error between the adjoint gradient and the finite-difference gradient for the boundary temperature parameters $T_L$ and $T_R$, respectively. As $N$ increases, the error decreases, demonstrating consistency between the adjoint and finite-difference methods. Additionally, Figures~\ref{fig:test1-TL-std}-\ref{fig:test1-TR-std} show the standard deviation (std) of the adjoint gradient for $T_L$ and $T_R$, respectively. This follows the expected Monte Carlo convergence rate of $\mathcal{O}(1/\sqrt{N})$ (plotted in black solid lines for reference). %

The observed decay of both the bias (relative error) and variance confirms that the adjoint DSMC method produces statistically consistent gradient estimates, even in the presence of thermal boundary reflections. In particular, the $\mathcal{O}(1/\sqrt{N})$ variance scaling demonstrates that the stochastic boundary regularization does not degrade Monte Carlo efficiency. These results validate the robustness and practical feasibility of the proposed adjoint framework for sensitivity analysis in spatially inhomogeneous Boltzmann simulations with thermal BCs.

\subsection{Mixed reflecting BC: thermal and specular reflections}\label{subsec:test2}
In our second test, we consider a mixed BC: at the left boundary, $x=0$,  thermal BC is enforced with a non-isotropic left half-Maxwellian flux distribution with the temperature $T_L  =  [0.6,0.5,0.8]$ and at the right boundary $x=1$, the specular reflection is enforced (see Section~\ref{subsec:R}). This mixed configuration creates an asymmetric setting in which one boundary injects thermal fluctuations while the other preserves kinetic energy through deterministic reflection. As a result, the dynamics combine stochastic boundary resampling and deterministic mirror interactions.

We use the same objective function as in Section~\ref{subsec:test1} but with the final time $T = 0.5$, which is denoted by $J_2$. The discretization parameters are $\Delta t = 0.05$, and $\Delta x = 0.1$. The initial condition is a uniform distribution over the spatial domain $\Omega = [0,1]$ and a Maxwellian distribution in velocity with temperature $T_0^2\cdot [1,1,1]$ where $T_0 = 1$. Similar to the example in~\Cref{subsec:test1}, a random time step is applied only to particles that may interact with the left boundary (associated with the thermal BC) during the current update.

{
Figure~\ref{fig:test2-macroscopic-profiles} presents the corresponding macroscopic profiles for the mixed thermal-specular configuration, again using $N_0=2\times10^6$ and $\Delta x=0.05$.  The initial density and temperature are approximately uniform, with $\widehat n\simeq1$ and $T\simeq1$, and the initial bulk velocity is statistically indistinguishable from zero.  By $T=0.5$, particles accumulate near the thermal left wall, where the normalized density reaches approximately $1.14$, while the density over most of the remaining domain is slightly below one.  The temperature is lowest near the thermal wall, where it is approximately $0.80$, and increases toward approximately $0.94$ near the specular boundary.  The resulting asymmetric transient produces a negative streamwise bulk velocity in the interior, with a minimum of approximately $-0.072$, whereas the initial velocity profile fluctuates only weakly around zero.

\begin{figure}
    \centering
    \includegraphics[width=\textwidth]{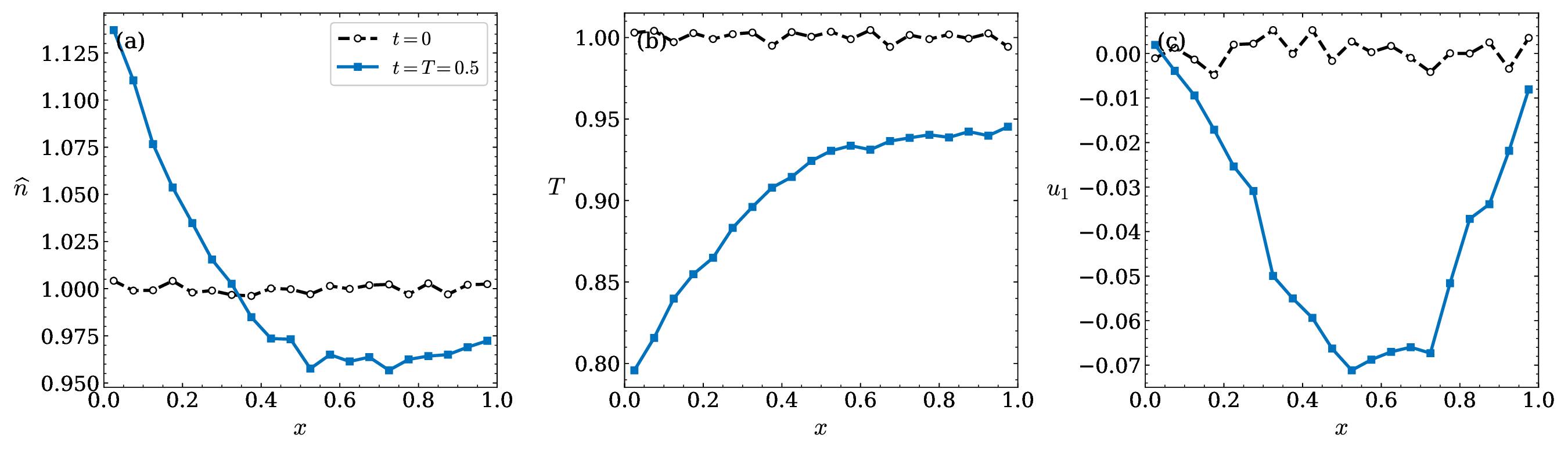}
    \caption{Mixed thermal and specular boundary conditions: profiles of
    (a) normalized number density $\widehat n$, (b) scalar temperature $T$,
    and (c) streamwise bulk velocity $u_1$ at $t=0$ and $t=T=0.5$.
    The simulation uses $N_0=2\times10^6$ and $\Delta x=0.05$.}
    \label{fig:test2-macroscopic-profiles}
\end{figure}

}
In this example, the objective depends both on boundary-driven energy injection (through $T_L$) and on the initial thermal state (through $T_0$). Consequently, the gradients $\nabla_{T_L} J$ and $\nabla_{T_0} J$ measure distinct physical sensitivities: the former captures the influence of boundary heating, while the latter quantifies how perturbations of the initial kinetic energy propagate over time under mixed reflection mechanisms.

\begin{figure}
    \centering
    \subfloat[$|\nabla^\text{AJ}_{T_0} J_1 -\nabla^\text{FD}_{T_0} J_1|/|\nabla^\text{FD}_{T_0} J_1|$]{\includegraphics[width = 0.48\textwidth]{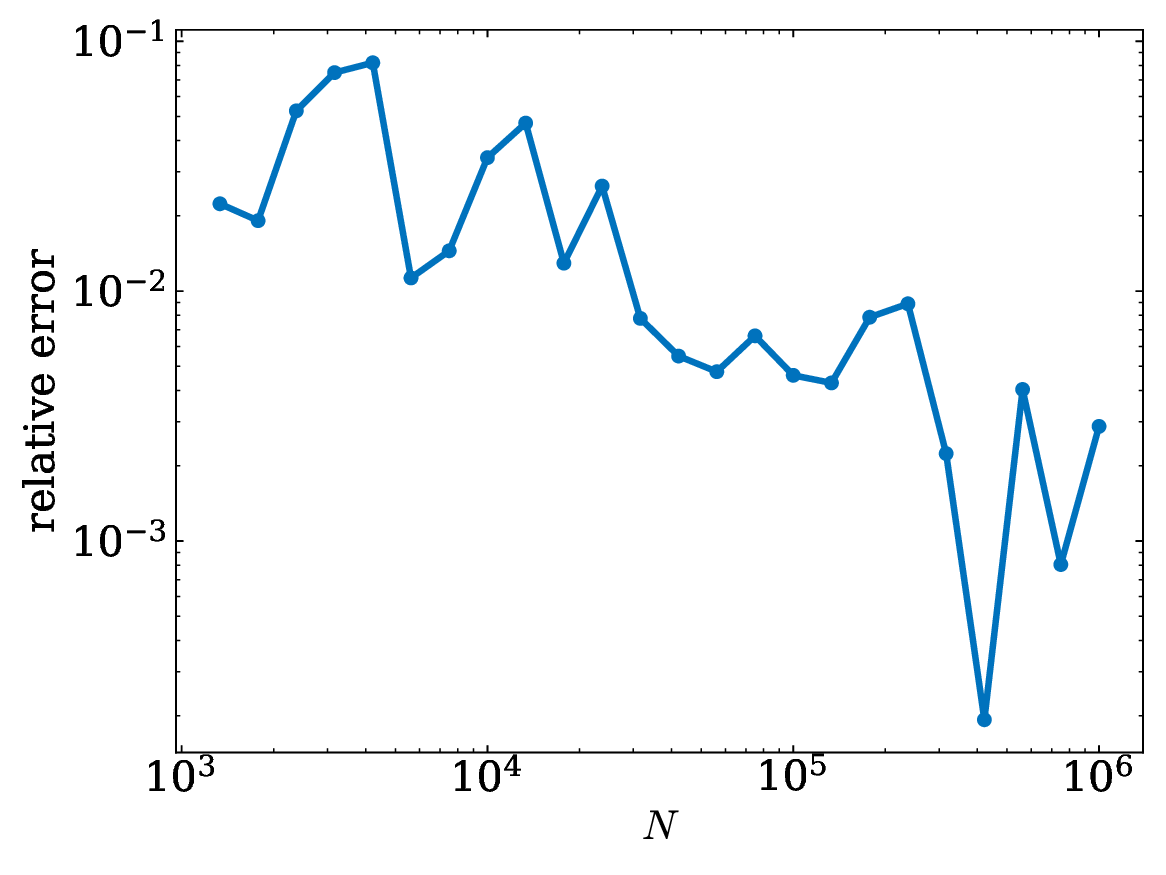}\label{fig:test2-T0-error}}
    \subfloat[std in $\nabla^\text{AJ}_{T_0} J_1,\nabla^\text{FD}_{T_0} J_1$]{\includegraphics[width = 0.48\textwidth]{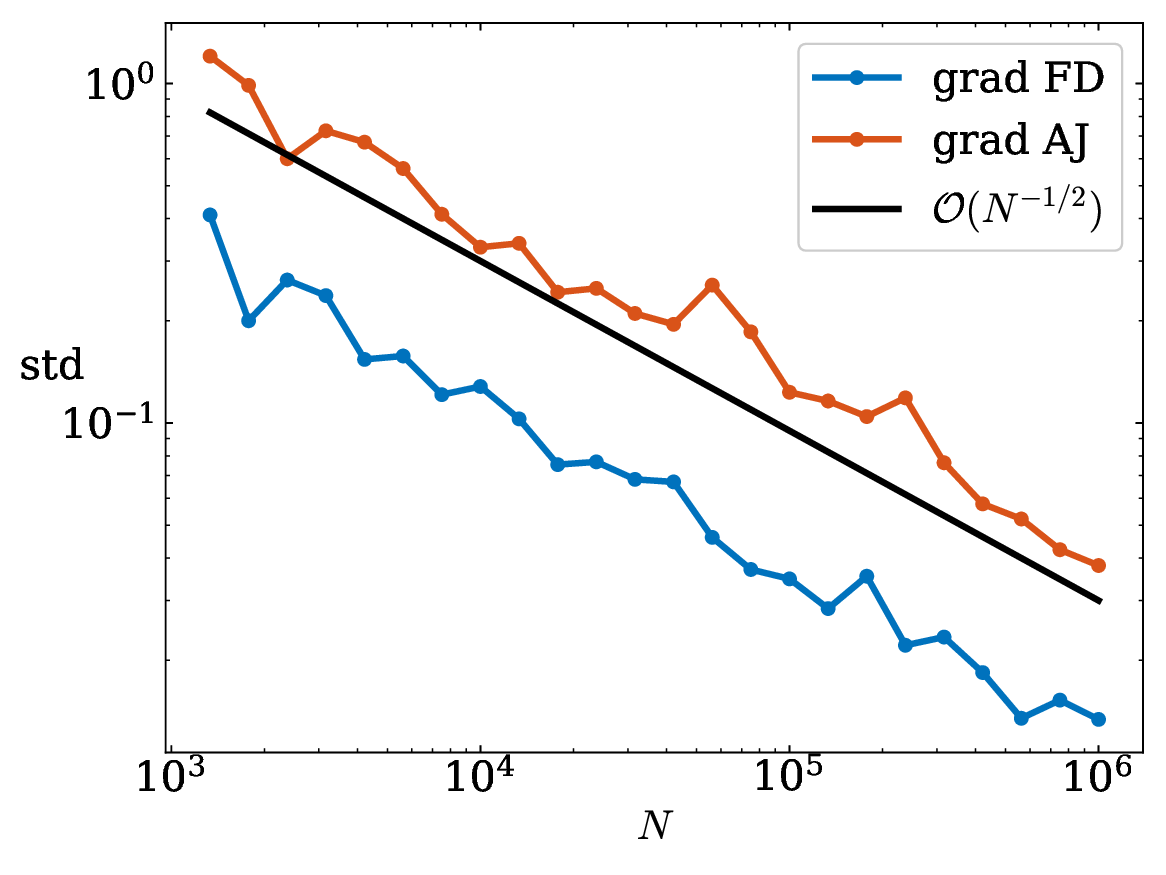}\label{fig:test2-T0-std}}\\
    \subfloat[$|\nabla^\text{AJ}_{T_L} J_1 -\nabla^\text{FD}_{T_L} J_1|/|\nabla^\text{FD}_{T_L} J_1|$]{\includegraphics[width = 0.48\textwidth]{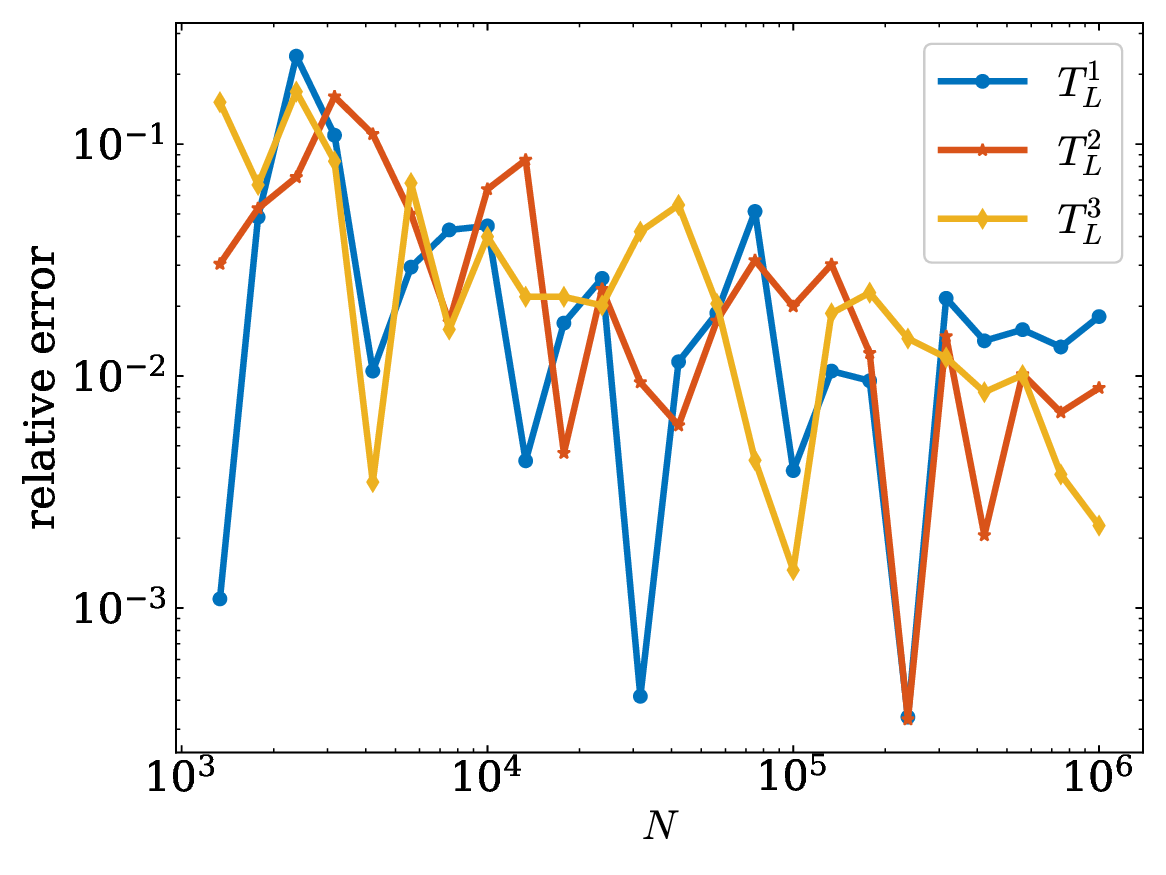}\label{fig:test2-TL-error}}
    \subfloat[standard deviation in $\nabla^\text{AJ}_{T_L} J_1$]{\includegraphics[width = 0.48\textwidth]{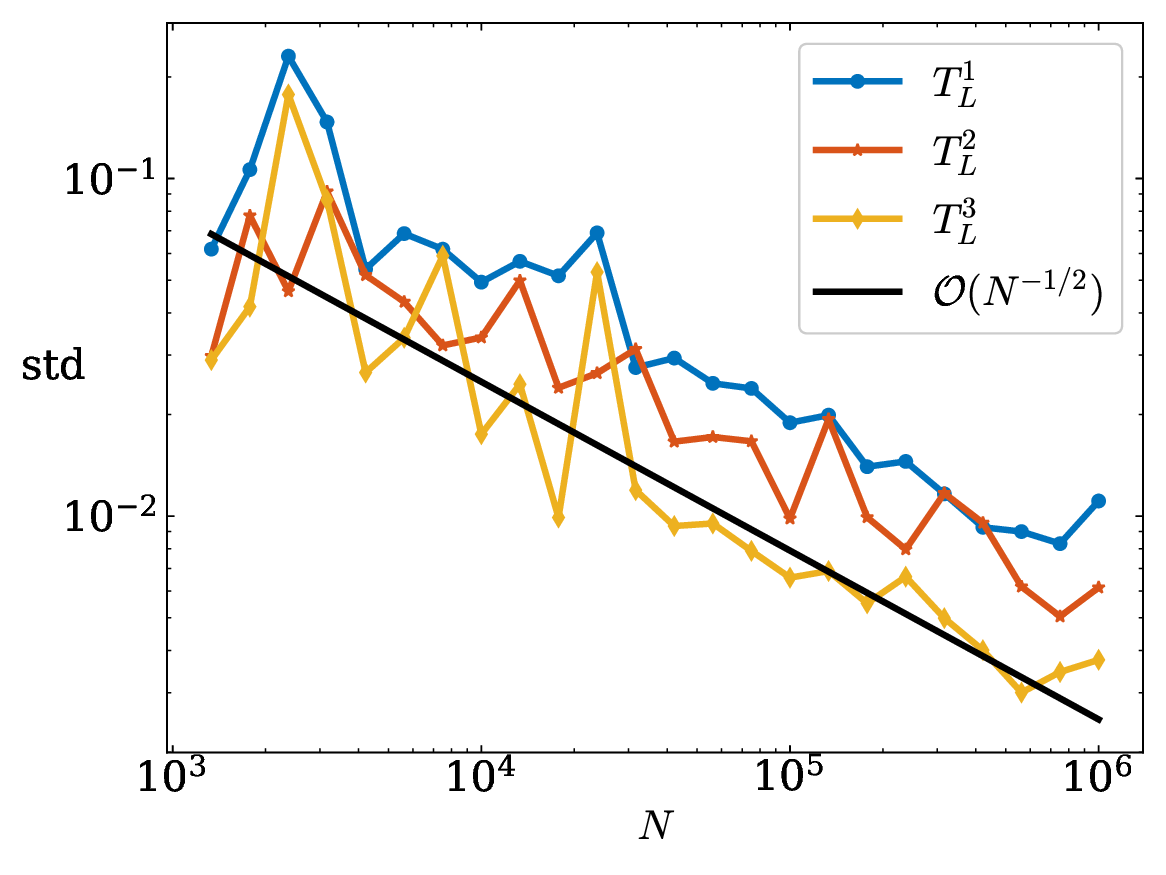}\label{fig:test2-TL-std}}\\
    \subfloat[$|\nabla^\text{AJ}_{a_x} J_1 -\nabla^\text{FD}_{a_x} J_1|$]{\includegraphics[width=0.48\textwidth]{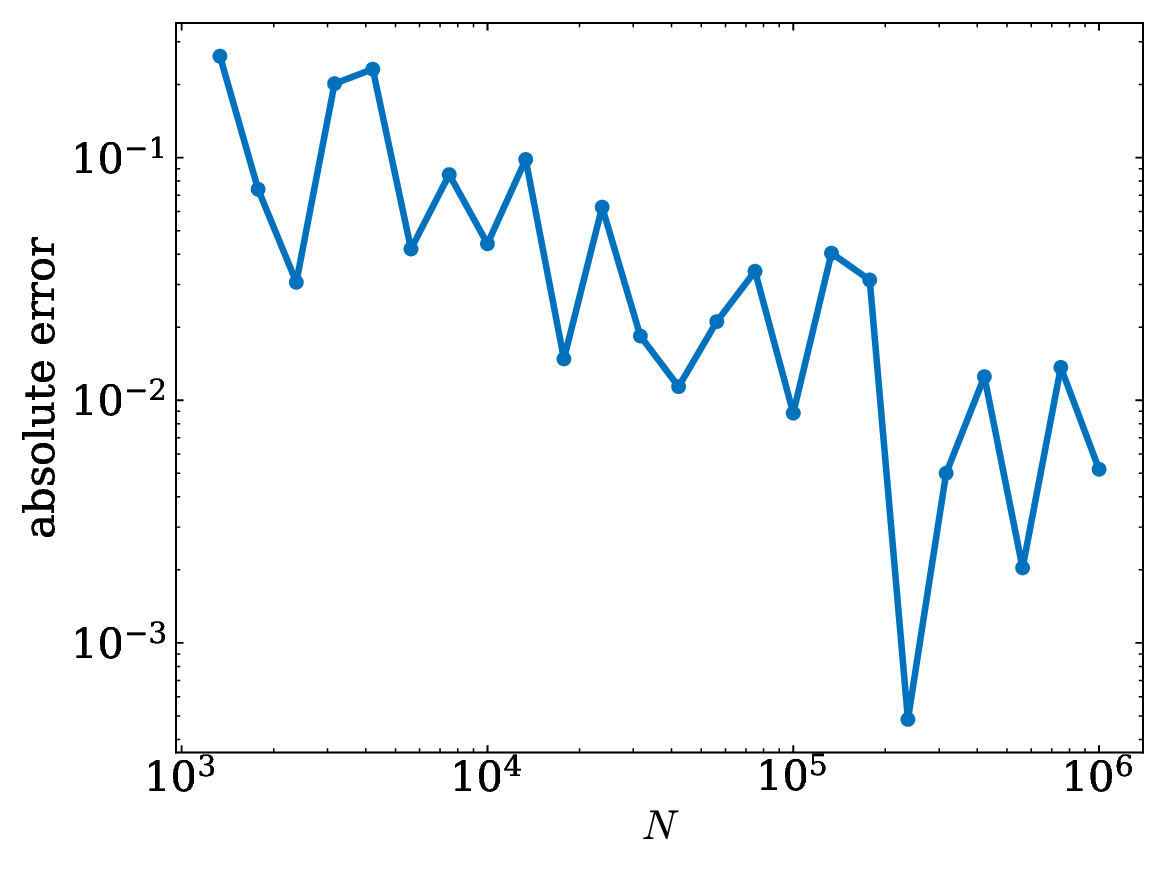}
    \label{fig:test2-Xinit-error}}
    \subfloat[std in $\nabla_{a_x}^{\mathrm{AJ}}J_1$ and $\nabla_{a_x}^{\mathrm{FD}}J_1$]{\includegraphics[width=0.48\textwidth]{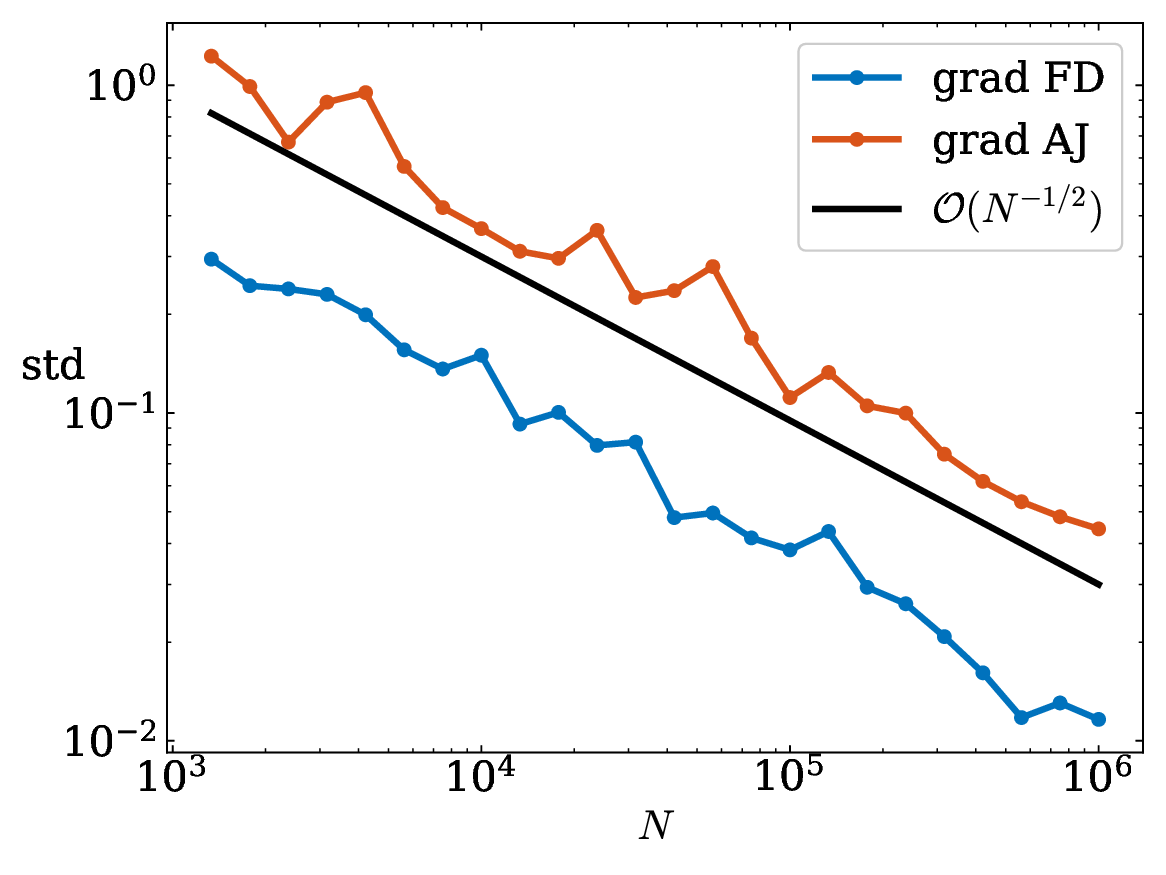}\label{fig:test2-Xinit-std}
}
\caption{Thermal and specular reflections BCs (\Cref{subsec:test2}): (a), (c) and (e) show the error between the adjoint gradient and the finite-difference gradient for $T_0$, $T_L$ and $a_x$, respectively, as the number of particles $N$ increases; (b), (d) and (f) show the standard deviation of the gradients with respect to these parameters, respectively. The standard deviation decays with a rate of $\mathcal{O}(\frac{1}{\sqrt{N}})$.}
    \label{fig:test2}
\end{figure}

We use the adjoint DSMC method to compute the gradients with respect to the thermal BC parameter $T_L$ and the initial distribution's parameter $T_0$. The results are then compared with the finite-difference approximation with the perturbation $\delta = 0.05$. The number of particles in the forward and adjoint DSMC ranges from $N=10^3$ to $N=10^6$. We use $42$ i.i.d.~runs to approximate the standard deviation of all gradients, and compare the difference between the mean values. Because specular reflection is deterministic and smooth with respect to particle states, no stochastic regularization is required at the right boundary. Therefore, this experiment isolates the effect of the thermal boundary regularization at $x=0$ and allows us to verify that the adjoint treatment remains accurate even when different reflection mechanisms coexist.

We plot both the error decay and standard deviation change in the finite-difference gradients $\nabla_{T_0}^\text{FD}$, $\nabla_{T_L}^\text{FD}$ and the adjoint DSMC gradients $\nabla_{T_0}^\text{AJ}$, $\nabla_{T_L}^\text{AJ}$ in Figure~\ref{fig:test2}. The observed convergence confirms that the adjoint DSMC method remains statistically consistent in the presence of mixed boundary conditions. In particular, the preservation of the $\mathcal{O}(1/\sqrt{N})$ variance scaling demonstrates that combining stochastic (thermal) and deterministic (specular) reflections does not introduce additional variance beyond standard Monte Carlo sampling error. These results further support the robustness of the adjoint framework for sensitivity analysis in heterogeneous boundary configurations.

{

To further examine sensitivity with respect to the initial spatial distribution, we introduce a scalar shape parameter $a_x>0$ and define
\[
\rho_0(x;a_x)=a_xx^{a_x-1},
\qquad x\in[0,1].
\]
The value $a_x=1$ recovers the uniform spatial distribution used in the original setup. Samples from this family are generated through the inverse-transform representation
\[
X_{0,i}=U_i^{1/a_x},
\qquad U_i\sim\mathcal{U}(0,1),
\]
for which
\[
\frac{\partial X_{0,i}}{\partial a_x}
= -\frac{X_{0,i}\log U_i}{a_x^2}.
\]
After propagating the position adjoint variables backward to the initial time, we combine them with these derivatives to obtain $\nabla_{a_x}^{\mathrm{AJ}}J_2$. We compare this result with the centered finite-difference approximation
\[
\nabla_{a_x}^{\mathrm{FD}}J_2
=
\frac{J_2(a_x+\delta)-J_2(a_x-\delta)}{2\delta},
\qquad \delta=0.05.
\]
The same uniform samples $\{U_i\}$ and random seeds are used for both perturbations to provide common random numbers.

Figure~\ref{fig:test2-Xinit-error} shows the absolute discrepancy
\[
\left|
\nabla_{a_x}^{\mathrm{AJ}}J_2
-
\nabla_{a_x}^{\mathrm{FD}}J_2
\right|
\]
as the number of particles increases. We report the absolute error rather than the relative error because the estimated reference initial-spatial-distribution gradient is itself very small, with a magnitude on the order of $10^{-2}$. Dividing by such a small reference gradient greatly amplifies Monte Carlo fluctuations and can produce large, statistically unstable relative errors even when the absolute difference between the two estimators is small. Although the absolute-error curve exhibits sampling fluctuations, its overall magnitude decreases as $N$ increases, falling from approximately $1.32\times10^{-1}$ at the smallest particle count to $6.03\times10^{-3}$ at $N=10^6$. This behavior indicates improving agreement between the adjoint and finite-difference sensitivities as the particle resolution is increased.

}

{
\subsection{Two-sided inflow boundary conditions}
\label{subsec:test3}
Our third numerical example considers inflow boundary conditions at both ends of the spatial domain $\Omega=[0,1]$. In contrast to the reflecting boundary conditions considered in the preceding examples, particles that leave the computational domain are removed, while new particles are introduced through both boundaries during every time step. Consequently, the number of active particles in the domain is not conserved.

The initial particle positions and velocities are sampled independently according to
\[
X_{0,i}\sim\mathcal U([0,1]),
\qquad
V_{0,i}\sim\mathcal N(0,I_3),
\qquad
1\leq i\leq N_0,
\]
where $N_0$ denotes the number of particles at the initial time. The inflow temperature vectors are
\[
T_L=[4,4,4]^\top,
\qquad
T_R=[4,4,4]^\top.
\]
Thus, both incoming distributions are isotropic and have the same  temperature, but their components are treated as independent boundary parameters in the gradient computation.

At the left boundary, incoming particles satisfy $v^1>0$, whereas particles entering through the right boundary satisfy $v^1<0$. For a temperature vector $T=[T^1,T^2,T^3]^\top$, the incoming half-Maxwellian flux samples are generated as
\[
v^1 =
\sqrt{T^1}\sqrt{-2\log(1-U)}
\quad\text{at }x=0,
\qquad
v^1 = -\sqrt{T^1}\sqrt{-2\log(1-U)}
\quad\text{at }x=1,
\]
and
\[
v^2=\sqrt{T^2}Z_2,
\qquad
v^3=\sqrt{T^3}Z_3,
\]
where $U\sim\mathcal U(0,1)$ and $Z_2,Z_3\sim\mathcal N(0,1)$ are independent. Independently of the incoming velocity, we sample a residence time $\xi\sim\mathcal U([0,\Delta t])$. The positions of the injected particles at the end of the time step are therefore $x_L=v_L^1\xi_L$ and $x_R=1+v_R^1\xi_R$, where $\xi_L$ and $\xi_R$ are independent of their corresponding incoming velocities.

The computational inflow number density is set to $n_{\mathrm{in}}^{(p)}=1000N_0$. At the  boundary temperatures, the numbers of particles introduced during one time step are
\[
N_L =
\operatorname{round}\!\left(
\frac{\sqrt{T_L^1}}{\sqrt{2\pi}}\,
\Delta x\,\Delta t\,n_{\mathrm{in}}^{(p)}
\right),
\qquad
N_R =
\operatorname{round}\!\left(
\frac{\sqrt{T_R^1}}{\sqrt{2\pi}}\,
\Delta x\,\Delta t\,n_{\mathrm{in}}^{(p)}
\right).
\]
The values of $N_L$ and $N_R$ are evaluated at the given temperatures and then held fixed in both the ensemble-adjoint and finite-difference calculations.

The spatial and temporal discretizations are $\Delta x=0.025$, $\Delta t=0.025$, $M=20$ and $ T=M\Delta t=0.5$. In each spatial cell $\Omega_j$, the local particle count, density, and collision frequency are
\[
N_k^j=|\mathcal V_k^j|,
\qquad
\rho_k^j=\frac{N_k^j}{N_0|\Omega_j|},
\qquad
\mu_k^j
=
\rho_k^j\int_{\mathbb S^2}q(\sigma)\,d\sigma.
\]
The number of collision pairs is obtained by unbiased stochastic rounding of $N_k^j\Delta t\,\mu_k^j/2$, and the collision particles are selected uniformly without replacement within the cell. In the numerical calculation, $
\int_{\mathbb S^2}q(\sigma)\,d\sigma=1$.

We consider $N_0=2\times10^6$. The injection counts are $N_L=N_R=997{,}356$. The normalization used in the particle approximation of the objective is therefore $ N_{\mathrm{norm}} = N_0+M(N_L+N_R) = 41{,}894{,}240$.

Figure~\ref{fig:test3-macroscopic-profiles} shows the macroscopic state generated by the two-sided inflow. The simulation starts with $N_0=2\times10^6$ particles and contains $32{,}033{,}761$ active particles at $T=0.5$. The normalized number density increases from its initial value $\widehat n\simeq1$ to values between approximately $15.35$ and $18.43$. The final scalar temperature ranges from approximately $3.73$ near the boundaries to $4.92$ in the interior. The streamwise bulk velocity is positive near the left boundary and negative near the right boundary, reaching a magnitude of approximately $0.466$, while its spatial mean remains close to zero. The approximate symmetry of the density and temperature profiles and the antisymmetry of the velocity profile are consistent with the equal left and right inflow temperatures and oppositely directed incoming fluxes.

\begin{figure}[t]
\centering
\includegraphics[width=\textwidth]
{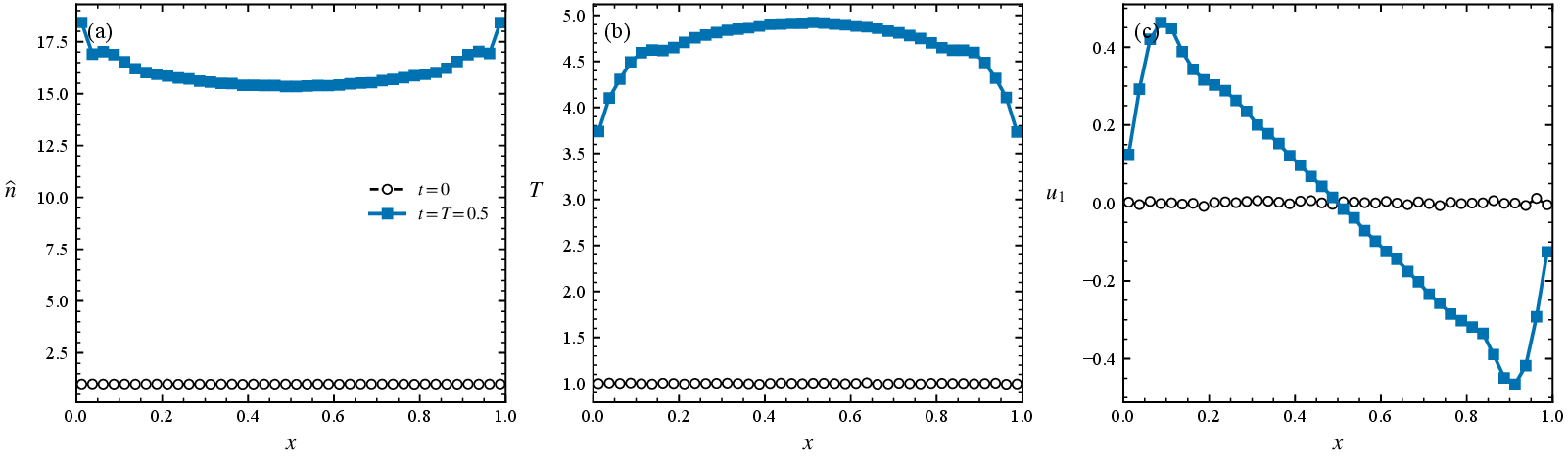}
\caption{Two-sided inflow boundary conditions: profiles of
(a) normalized number density $\widehat n$, (b) scalar temperature $T$, and
(c) streamwise bulk velocity $u_1$ at $t=0$ and $t=T=0.5$. The simulation
uses $N_0=2\times10^6$, $\Delta x=0.025$, and the local-density-dependent
collision frequency.}
\label{fig:test3-macroscopic-profiles}
\end{figure}

For this experiment, we consider the localized kinetic-energy objective
\begin{equation}
J_3(T_L,T_R)
=
\int_0^1\int_{\mathbb R^3}
|v|^2\exp\!\left[-36(x-0.5)^2\right]
f(x,v,T)\,dv\,dx.
\label{eq:test3_objective}
\end{equation}
Its particle approximation is
\[
\widehat J_3
=
\frac{1}{N_{\mathrm{norm}}}
\sum_{i\in\mathcal I_M}
|V_{M,i}|^2
\exp\!\left[-36(X_{M,i}-0.5)^2\right],
\]
where $\mathcal I_M$ denotes the particles that remain active at the terminal time. The Gaussian weight is centered at the midpoint of the domain and has the boundary value $e^{-9}$. Thus, the objective function measures kinetic energy transported through the interior and reduces the direct contribution from particles near the outflow boundaries.

We compute the gradient using the ensemble-adjoint formulation described in Section~\ref{sec:adjoint_DSMC_inflow}. The backward calculation propagates a scalar distributional adjoint through free transport and the adjoint of the locally linearized collision operator. The forward simulation uses the complete particle population to evaluate the exact cell counts and local collision frequencies, while a subset of the particle states is retained for the backward calculation. At every injection time, the boundary contribution is evaluated by multiplying the adjoint values at the sampled incoming states by the analytic temperature score of the half-Maxwellian flux distribution and summing over the source samples.

In this test, we focus on the two boundary-normal temperature derivatives
\[
G_L=\frac{\partial J_3}{\partial T_L^1},
\qquad
G_R=\frac{\partial J_3}{\partial T_R^1}.
\]
These components control both the incoming normal-speed distribution and the distribution of the injected particle positions. The centered finite-difference approximations are
\[
G_L^{\mathrm{FD}} =
\frac{
J_3(T_L+h e_1,T_R)-J_3(T_L-h e_1,T_R)
}{2h},
\]
and
\[
G_R^{\mathrm{FD}} =
\frac{
J_3(T_L,T_R+h e_1)-J_3(T_L,T_R-h e_1)
}{2h},
\]
where $h=0.025$. The values of $N_L$ and $N_R$ are retained in all perturbed forward calculations, and common random numbers are used within each centered finite-difference pair.

The validation uses eight independently seeded collision-and-inflow realizations, all sharing the same initial particle sample. Each realization uses $N_0=2\times10^6$, $32{,}000$ retained particles per spatial cell, and $800{,}000$ stored boundary-source samples per side. The ensemble-adjoint means and standard errors are
\[
G_L^{\mathrm{AJ}} = 0.15171883\pm0.00012830,
\qquad
G_R^{\mathrm{AJ}} = 0.15157001\pm0.00023059,
\]
whereas the corresponding finite-difference estimates are
\[
G_L^{\mathrm{FD}} = 0.14986940\pm0.00753680,
\qquad
G_R^{\mathrm{FD}} = 0.15028809\pm0.00631838.
\]
Here, each displayed uncertainty is one standard error across the eight realizations. The relative differences between the adjoint and finite-difference means are $1.23\%$ and $0.85\%$ at the left and right boundaries, respectively.

Figure~\ref{fig:test3-gradient-validation}(a) compares the two gradient estimators, with error bars showing twice the standard error. Panel (b) shows the realization-wise paired differences $G_{\mathrm{AJ}}-G_{\mathrm{FD}}$. Their $95\%$ Student-$t$ confidence intervals are $ [-0.01598,\,0.01968]$ and $ [-0.01349,\,0.01605]$ at the left and right boundaries, respectively. Both intervals contain zero, so the observed mean discrepancies are not statistically significant at the $5\%$ level. The substantially smaller error bars of the adjoint estimates also demonstrate that, for this calculation, the ensemble-adjoint estimator has much lower sampling variability than the centered finite difference.

\begin{figure}[t]
\centering
\includegraphics[width=0.92\textwidth]
{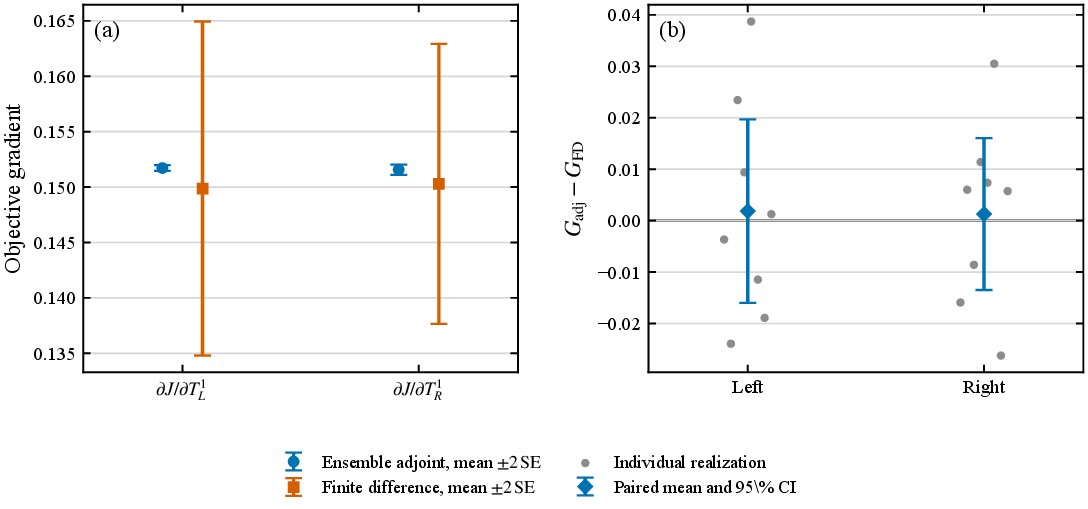}
\caption{Ensemble-adjoint validation for the  inflow boundary parameters. Each of the eight independently seeded realizations uses
$N_0=2\times10^6$, finite-difference step $h=0.025$, and $32{,}000$  retained
particles per cell. Panel (a) compares the ensemble-adjoint and
finite-difference gradients; the error bars show two standard errors. Panel (b)
shows the paired differences $G_{\rm AJ}-G_{\rm FD}$ for the individual
realizations together with the paired mean and its $95\%$ Student-$t$
confidence interval. Both confidence intervals contain zero.}
\label{fig:test3-gradient-validation}
\end{figure}

We additionally examine the dependence of the adjoint estimates on the number of retained particles in each spatial cell. The results with $2{,}000$, $8{,}000$, and $32{,}000$ retained particles per cell  are shown in Figure~\ref{fig:test3-reservoir-convergence}. Increasing the number of retained particles moves both boundary-gradient estimates toward the eight-realization finite-difference means.

A further calculation with $64{,}000$ retained particles per cell and approximately $10^6$ boundary-source samples per side continues this trend. This largest calculation contains only one realization, so it is shown without a sampling error bar.

\begin{figure}[t]
\centering
\includegraphics[width=0.92\textwidth]
{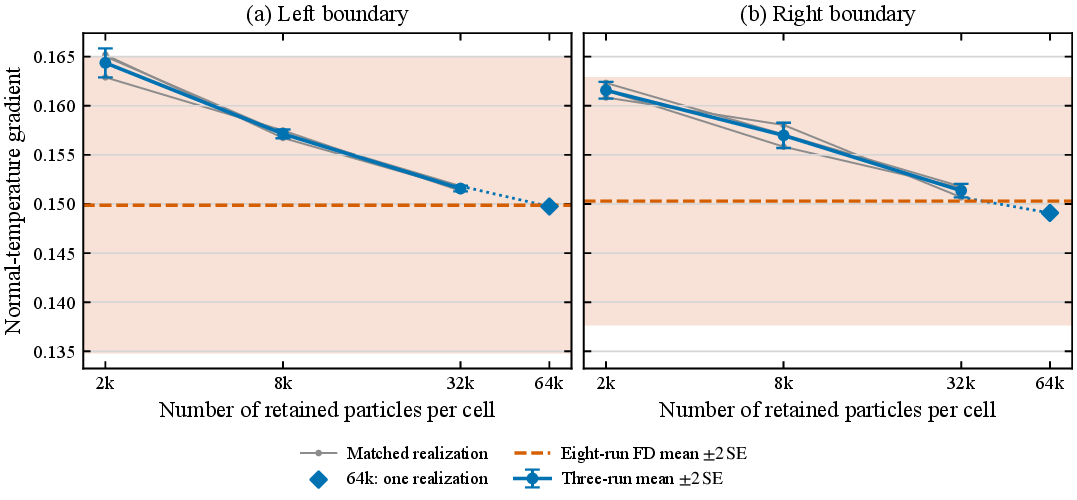}
\caption{Convergence of the ensemble-adjoint gradients regarding (a) left and (b) right boundary-normal inflow temperatures with respect to the number of retained particles.
The points through $32{,}000$ retained particles per cell show the mean and
two standard errors from three matched realizations; thin gray curves show
the individual matched realizations. The $64{,}000$-particle result is from
one realization and is shown without a sampling error bar. The dashed line
and shaded band denote the centered finite-difference mean and two standard
errors from eight independently seeded realizations. The number of stored
boundary-source samples is refined together with the number of retained particles.}
\label{fig:test3-reservoir-convergence}
\end{figure}
}

{

\subsection{High-Mach planar Couette flow}\label{subsec:test4}

To examine the adjoint DSMC method in a regime where the sensitivity signal is large relative to the statistical fluctuations, we consider a planar Couette-flow problem with two oppositely moving thermal walls. As in the preceding examples, the spatial domain is $\Omega=[0,1]$ and particles carry three-dimensional velocities. Both boundaries are modeled as diffuse thermal walls with the same isotropic temperature,
\[
T_w=1,\qquad T_L=T_R=T_w[1,1,1]^\top.
\]
The walls move in opposite directions tangential to the spatial domain with the speed of:
\[
U_L=[0,-U_w,0]^\top,
\qquad
U_R=[0,U_w,0]^\top,
\qquad U_w=4.
\]
Taking the ratio of the wall speed to the equilibrium sound speed, with $\gamma=5/3$, gives
\[
\mathrm{Ma}_w
=
\frac{U_w}{\sqrt{\gamma T_w}}
=
\frac{4}{\sqrt{5/3}}
\approx 3.10.
\]
Thus, the prescribed boundary motion generates a high-Mach-number shear flow whose mean tangential velocity is substantially larger than the equilibrium velocity fluctuations.

For a particle diffusely re-emitted from the left wall, the velocity is sampled according to
\[
V_1=\sqrt{-2T_w\log(1-\xi_1)},\qquad V_2=-U_w+\sqrt{T_w}\,\xi_2,\qquad V_3=\sqrt{T_w}\,\xi_3,
\]
where $\xi_1\sim\mathcal U(0,1)$ and $\xi_2,\xi_3\sim\mathcal N(0,1)$. At the right wall, the sign of the normal velocity is reversed and the tangential mean velocity is $+U_w$:
\[
V_1=-\sqrt{-2T_w\log(1-\xi_1)},\qquad V_2=U_w+\sqrt{T_w}\,\xi_2,\qquad V_3=\sqrt{T_w}\,\xi_3.
\]
Consequently, the derivative of a re-emitted tangential velocity with respect to the wall-speed magnitude is $-1$ at the left wall and $+1$ at the right wall. The corresponding directional adjoint gradient is therefore
\[
\frac{\partial J_4}{\partial U_w}
= -\frac{\partial J_4}{\partial U_{L,2}} +\frac{\partial J_4}{\partial U_{R,2}}.
\]

The gas is initially uniform in space and Maxwellian in velocity:
\[
X_i(0)\sim\mathcal U([0,1]),
\qquad
V_i(0)\sim\mathcal N(0,I_3),
\qquad 1\leq i\leq N.
\]
We use the spatial-cell width $\Delta x=0.05$, the time-step size $\Delta t=0.05$, and $10$ time steps, giving a final time $T=0.5$. For the adjoint calculation, the regularization width is $\varepsilon=\Delta t/10=0.005$. As in the preceding thermal-wall tests, the randomized time increment is used only for particles that may interact with a thermal wall during the current update.

Figure~\ref{fig:test4-profiles} shows the normalized number density, scalar temperature, and tangential bulk velocity at the initial and final times, computed using $N=2\times10^5$ particles. The density remains approximately uniform, whereas viscous heating increases the gas temperature from its initial value of approximately $1$ to approximately $2.62$. More importantly, the initially vanishing tangential bulk velocity develops a clear, approximately linear shear profile. At the cell centers nearest the two walls, the final bulk velocities are approximately $-1.98$ and $1.99$, respectively. The large and smoothly resolved velocity profile confirms that the macroscopic Couette signal is not obscured by thermal sampling noise.

\begin{figure}
    \centering
    \includegraphics[width=\textwidth]
    {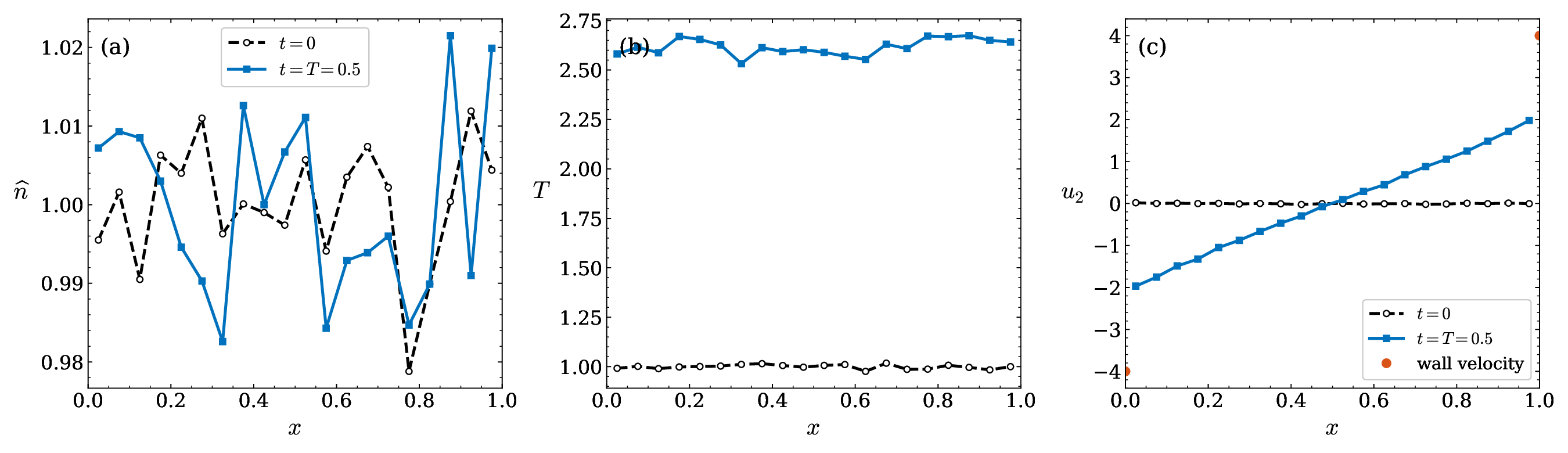}
    \caption{High-Mach Couette-flow macroscopic profiles computed with
    $N=2\times10^5$, $\Delta x=0.05$, and $T=0.5$:
    (a) normalized number density, (b) scalar temperature, and
    (c) tangential bulk velocity. The black dashed and blue solid curves
    represent the profiles at $t=0$ and $t=T$, respectively. The orange
    markers in panel (c) indicate the prescribed wall velocities
    $U_{L,2}=-4$ and $U_{R,2}=4$.}
    \label{fig:test4-profiles}
\end{figure}

To obtain a sensitivity with a comparably strong signal, we define the objective as the final mean tangential kinetic energy:
\[
J_4(U_w) =
\frac{1}{N}
\sum_{i=1}^{N} V_{i,2}(T)^2.
\]
This objective directly measures the kinetic energy generated by the moving walls. Its adjoint terminal condition is
\[
\frac{\partial J_4}{\partial V_{i,2}(T)}
=
\frac{2V_{i,2}(T)}{N},
\qquad
\frac{\partial J_4}{\partial X_i(T)}=0.
\]
Because $U_w$ changes the two wall velocities in opposite directions, the contributions from the two boundaries add in $\partial J_4/\partial U_w$ rather than cancel.

We compare the adjoint gradient with the centered finite-difference approximation
\[
\nabla_{U_w}^{\mathrm{FD}}J_4
=
\frac{J_4(U_w+\delta)-J_4(U_w-\delta)}{2\delta},
\qquad
\delta=0.05.
\]
The positive and negative perturbations use common random numbers. The finite-difference simulations employ the original deterministic time step, while the adjoint calculation uses the locally regularized boundary treatment. Thus, agreement between the two estimators is not produced by applying the same regularization to both calculations.

We vary the number of particles from $N=10^3$ to $N=2\times10^5$ and perform $96$ independent realizations for every particle count. In addition to the relative difference between the mean gradients, we measure the statistical noise using the coefficient of variation
\[
\operatorname{CV}(G)
=
\frac{\operatorname{std}(G)}
{|\mathbb E[G]|}.
\]
Unlike an absolute standard deviation, this quantity measures the fluctuations relative to the magnitude of the sensitivity signal and is therefore particularly suitable for assessing the noise level in this example.

The results are presented in Figure~\ref{fig:test4-gradient}. The adjoint and finite-difference mean gradients are statistically compatible over the resolved range; their relative difference is below \(2\%\) for every \(N\geq2000\). At $N=2\times10^5$, we obtain
\[
\nabla_{U_w}^{\mathrm{AJ}}J_4=3.0283,
\qquad
\nabla_{U_w}^{\mathrm{FD}}J_4=3.0011,
\]
corresponding to a relative difference of approximately $0.91\%$. At the same particle count, the coefficient of variation is approximately $1.87\%$ for the adjoint estimator and $5.33\%$ for the finite-difference estimator. Therefore, the sensitivity is resolved with relatively small statistical uncertainty, and the adjoint estimator is noticeably less noisy than the finite-difference benchmark.

\begin{figure}
    \centering
    \includegraphics[width=\textwidth]
    {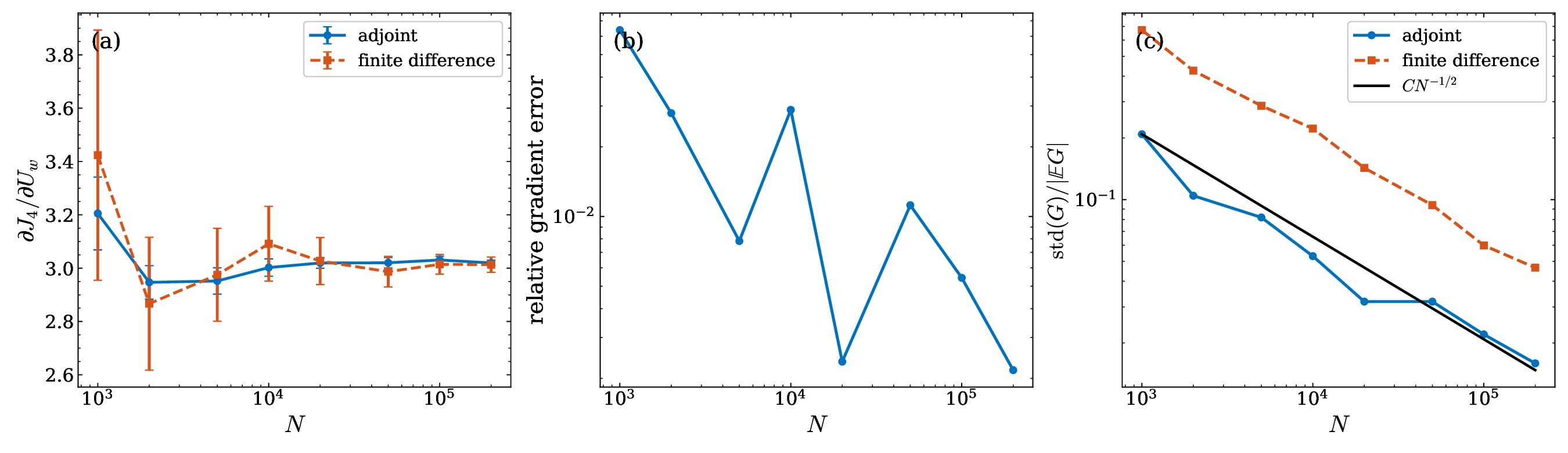}
    \caption{Adjoint-gradient validation for the high-Mach Couette-flow
    problem: (a) mean adjoint and finite-difference estimates of
    $\partial J_4/\partial U_w$, with error bars equal to twice the standard
    error of the mean; (b) relative difference between the mean adjoint and
    finite-difference gradients; and (c) coefficients of variation of the two
    gradient estimators. The black line in panel (c) shows an
    $N^{-1/2}$ reference rate. Results are based on $96$ independent
    realizations for every particle count.}
    \label{fig:test4-gradient}
\end{figure}

This experiment complements the preceding near-equilibrium tests by examining a strongly driven flow with a clearly resolved macroscopic response. The agreement with the finite-difference benchmark demonstrates that the boundary adjoint remains accurate when the thermal walls possess large tangential velocities. Moreover, the small coefficient of variation shows that the observed agreement is not masked by statistical noise. The comparison also illustrates a practical advantage of the adjoint estimator: in the large-particle regime it has substantially smaller relative fluctuations than the centered finite-difference approximation.
}

{

\section{Discussions}\label{sec:discussion}

This section discusses extensions to general collision kernels, the memory requirements of the adjoint formulations, the roles of objective selection and problem identifiability, and the treatment of objectives involving intermediate times.

\subsection{Adjoint DSMC for a different collision kernel} The boundary adjoint relations derived in this work are associated with the transport and boundary maps, and therefore do not rely on the Maxwell-molecule collision model. If a more realistic molecular potential, such as hard spheres, is used, then the boundary part of the adjoint remains unchanged provided the same boundary reflection law is imposed. The collision part of the adjoint, however, must be replaced by the corresponding general-collision update, which includes the velocity-dependent collision rate and the associated score-function and adjoint Jacobian terms.

\subsection{Memory consideration for Adjoint DSMC}

The pathwise adjoint formulations for periodic, specular, and thermal boundary conditions require access to the realized forward trajectory. The particle states, collision pairings, boundary events, and sampled boundary variables may be either stored or reconstructed via checkpointing. Once this information is fixed, the corresponding backward sweep is deterministic.

The ensemble-adjoint inflow method has a different storage pattern. The forward calculation uses the complete active particle population, but the backward calculation stores only the exact cell counts, a prescribed subset of particle states in each cell, and samples of the boundary sources. The scalar adjoint field is reconstructed from these retained states, and its collision contribution is evaluated by Monte Carlo quadrature. 

\subsection{Gradient calculation versus optimization landscape}
We emphasize that the adjoint DSMC method developed here provides gradients for a prescribed objective function; it does not by itself guarantee uniqueness of the associated inverse problem. If the objective depends only on a moment of the final-time distribution, then different particle distributions, or different parameter values, may in principle produce the same objective value. As in many PDE-constrained optimization problems, uniqueness and convexity depend on the chosen parameterization, the available observations, and any regularization terms. Additional observables, time-integrated objectives, or regularization may be used to improve identifiability, but a general uniqueness theory is beyond the scope of this work.

\subsection{Choice of the objective function} In the present work, the adjoint variables are not used to reweight the current particle cloud into an arbitrary target distribution. Instead, they propagate the sensitivity of a prescribed scalar objective with respect to the model or boundary parameters. For the moment-based objectives considered here, the target is a finite-dimensional quantity, not a full probability distribution, and therefore a support condition between two distributions is not required. If one wishes to formulate a full distribution-matching inverse problem, then the objective should be designed to avoid support mismatch, for example by using smooth observables, regularized discrepancies, time-integrated information, or continuation strategies.

\subsection{Objectives involving intermediate times}
\label{sec:time_integrated_objective}

For clarity, the preceding derivations use an objective depending only on the particle distribution at the final time. The same adjoint framework also applies to objectives involving the solution at intermediate times. For example, consider the time-integrated moment
\begin{equation}
\label{eq:time_integrated_objective}
\cJ_{\rm int}(m)
=
\int_0^T
\int_{\Omega}
\int_{\mathbb R^3}
r(x,v,t)f(x,v,t;m)
\,\mathrm{d}v\,\mathrm{d}x\,\mathrm{d}t .
\end{equation}
Let $t_k=k\Delta t$, and let $\omega_k$ denote the weights of a temporal quadrature rule. For example, $
\omega_k=\Delta t$, $0\leq k\leq M-1$, and $\omega_M=0$.
A DSMC approximation of~\eqref{eq:time_integrated_objective} is
\begin{equation}
\label{eq:time_integrated_objective_DSMC}
\cJ_{\rm int}^N(m)
=
\frac{1}{N}
\sum_{k=0}^{M}
\omega_k
\sum_{i\in\mathcal I_k}
r(x_{k,i},v_{k,i},t_k),
\end{equation}
where $\mathcal I_k$ denotes the set of particles present in the spatial domain at time $t_k$. For periodic, specular reflecting, and thermal BCs, $\mathcal I_k=\{1,\ldots,N\}$. For inflow BC, the sum is restricted to particles that have entered the domain and have not yet exited.

For the pathwise formulations associated with periodic, specular, and thermal boundary conditions, the running objective produces source terms in the particle adjoints $\alpha_{k,i}$ and $\beta_{k,i}$, as described below. Here, $\widehat\alpha_{k,i}$ and $\widehat\beta_{k,i}$ denote, respectively, the position and velocity adjoint values obtained at time $t_k$ by propagating the adjoints from later time levels backward through the collision, transport, and boundary updates, before adding the contribution of the running objective at $t_k$.
\begin{align}
\alpha_{k,i}
&=
\widehat\alpha_{k,i}
-
\omega_k
\partial_x r(x_{k,i},v_{k,i},t_k),
\label{eq:running_alpha}
\\
\beta_{k,i}
&=
\widehat\beta_{k,i}
-
\omega_k
\partial_v r(x_{k,i},v_{k,i},t_k).
\label{eq:running_beta}
\end{align}
The terminal conditions become
\begin{align}
\alpha_{M,i}
&=
-\omega_M
\partial_x r(x_{M,i},v_{M,i},t_M),
\\
\beta_{M,i}
&=
-\omega_M
\partial_v r(x_{M,i},v_{M,i},t_M).
\end{align}
Thus, if $\omega_M=0$, the terminal adjoint variables are zero, and the objective enters the adjoint system through the source terms in \eqref{eq:running_alpha}-\eqref{eq:running_beta}. The terminal-time objective considered in Sections~\ref{sec:adjoint_DSMC_space},~\ref{sec:adjoint_DSMC_thermal} and~\ref{sec:adjoint_DSMC_inflow} is recovered by taking $\omega_M=1$ and $\omega_k=0$ for $k<M$.

For the ensemble-adjoint inflow formulation, the same principle is expressed at the distributional level. The derivative of the running observable enters as a source in the backward equation for $\lambda_k(x,v)$. Consequently, the adjoint value used in each local inflow-score contribution contains all subsequent running and terminal objective contributions.

For thermal boundary conditions, the same running source terms are added to the adjoint equations. In addition, the quantity multiplying a score-function term at time $t_k$ must contain the portion of the objective affected by the randomized boundary transition at that time. Thus, the terminal objective $r_i$ appearing in the thermal-boundary score-function terms in~\eqref{eq:alpha_adjoint}-\eqref{eq:beta_adjoint} need also be replaced.

Once the modified adjoint variables at the initial time have been computed, the gradient formulas retain the same form.  Hence, incorporating observations at intermediate times does not require any change to the forward DSMC simulation and requires only additional source terms and, for thermal boundaries, a modification of the score term.
}

\section{Conclusions}\label{sec:conclusion}

We have developed complementary adjoint DSMC formulations for the spatially inhomogeneous Boltzmann equation with periodic, specular reflecting, diffuse thermal, and prescribed inflow boundary conditions. Periodic and specular boundaries are treated using a pathwise particle adjoint conditional on the realized event history.

For diffuse thermal boundaries, a randomized-time regularization smooths the wall hit-or-miss event and permits its sensitivity to be represented by score-function terms. Reparameterized half-Maxwellian samples provide gradients with respect to wall temperatures and tangential wall velocities. For fixed $\varepsilon$, the method differentiates the regularized objective $J_\varepsilon$, which converges to the deterministic objective as $\varepsilon\to0$. The numerical stability study indicates that the gradient is insensitive to $\varepsilon$ over the range used in the experiments.

For prescribed inflow, we introduced an ensemble adjoint based on the locally linearized collision operator and a local likelihood-ratio score for the incoming boundary source. This construction accounts for the effect of perturbed incoming states on subsequent cell populations, local collision frequencies, and collision evolution. The present implementation conditions on fixed nominal injection counts; derivatives of a parameter-dependent physical injection intensity must be added separately.

Our results provide a practical, mathematically consistent adjoint framework for Boltzmann-constrained optimization and, more broadly, for particle-based adjoint methods for transport-dominated PDEs under realistic boundary conditions.

Our adjoint framework should  naturally extend to more general stochastic reflection models, such as the Cercignani--Lampis boundary condition, provided the reflection law is implemented via an explicit sampling procedure or a known conditional density for the post-reflection velocity of the form $p(v'|v,x)$.

{A near-term direction of our work is the development and public release of a documented software package implementing the proposed adjoint DSMC framework.}

\section*{Acknowledgment}
Y.~Y.~was supported in part by the National Science Foundation under award DMS-2409855 and DMS-2540324 and by Office of Naval Research Grant under award N00014-24-1-2088. The authors thank Denis Silantyev for constructive discussions.

\section*{AI Disclosure Statement}
During the preparation of this work, the authors used generative AI to improve readability and language tone. After using this tool, the authors reviewed and edited the content as needed and took full responsibility for the publication's content.

\bibliographystyle{plain}
\bibliography{references}

@article{Sapienza2024,
  title={Differentiable programming for differential equations: A review},
  author={Sapienza, Facundo and Bolibar, Jordi and Sch{\"a}fer, Frank and Groenke, Brian and Pal, Avik and Boussange, Victor and Heimbach, Patrick and Hooker, Giles and P{\'e}rez, Fernando and Persson, Per-Olof and others},
  journal={arXiv preprint arXiv:2406.09699},
  year={2024}
}

@article{Baydin2018,
  title={Automatic differentiation in machine learning: a survey},
  author={Baydin, Atilim Gunes and Pearlmutter, Barak A and Radul, Alexey Andreyevich and Siskind, Jeffrey Mark},
  journal={Journal of machine learning research},
  volume={18},
  number={153},
  pages={1--43},
  year={2018}
}

@book{GriewankWalther2008,
  title={Evaluating derivatives: principles and techniques of algorithmic differentiation},
  author={Griewank, Andreas and Walther, Andrea},
  year={2008},
  publisher={SIAM}
}

@article{yuan2025adjoint,
  title={Adjoint shape optimization from the continuum to free-molecular gas flows},
  author={Yuan, Ruifeng and Wu, Lei},
  journal={Journal of Computational Physics},
  pages={114102},
  year={2025},
  publisher={Elsevier}
}

@article{yuan2024design,
  title={A design optimization method for rarefied and continuum gas flows},
  author={Yuan, Ruifeng and Wu, Lei},
  journal={Journal of Computational Physics},
  volume={517},
  pages={113366},
  year={2024},
  publisher={Elsevier}
}

@article{guan2024topology,
  title={Topology optimization of rarefied gas flows using an adjoint discrete velocity method},
  author={Guan, Kaiwen and Yamada, Takayuki},
  journal={Journal of Computational Physics},
  volume={511},
  pages={113111},
  year={2024},
  publisher={Elsevier}
}

@article{caflisch2024adjoint,
  title={{Adjoint Monte Carlo Method}},
  author={Caflisch, Russel and Yang, Yunan},
  journal={Active Particles, Volume 4: Theory, Models, Applications},
  pages={461--505},
  year={2024},
  publisher={Springer}
}

@article{caflisch2021adjoint,
  title={Adjoint {DSMC} for nonlinear {B}oltzmann equation constrained optimization},
  author={Caflisch, Russel and Silantyev, Denis and Yang, Yunan},
  journal={Journal of Computational Physics},
  volume={439},
  pages={110404},
  year={2021},
  publisher={Elsevier}
}

@article{mohamed2020monte,
  title={{Monte Carlo gradient estimation in machine learning}},
  author={Mohamed, Shakir and Rosca, Mihaela and Figurnov, Michael and Mnih, Andriy},
  journal={Journal of Machine Learning Research},
  volume={21},
  number={132},
  pages={1--62},
  year={2020}
}

@article{rubinstein1986score,
  title={The score function approach for sensitivity analysis of computer simulation models},
  author={Rubinstein, Reuven Y},
  journal={Mathematics and Computers in Simulation},
  volume={28},
  number={5},
  pages={351--379},
  year={1986},
  publisher={Elsevier}
}

@incollection{biegler2003large,
  title={Large-scale {PDE}-constrained optimization: an introduction},
  author={Biegler, Lorenz T and Ghattas, Omar and Heinkenschloss, Matthias and van Bloemen Waanders, Bart},
  booktitle={Large-Scale {PDE}-Constrained Optimization},
  pages={3--13},
  year={2003},
  publisher={Springer}
}

@inproceedings{lovbak2022reversible,
  title={{Reversible random number generation for adjoint Monte Carlo simulation of the heat equation}},
  author={L{\o}vbak, Emil and Blondeel, Fr{\'e}d{\'e}ric and Lee, Adam and Vanroye, Lander and Van Barel, Andreas and Samaey, Giovanni},
  booktitle={International Conference on Monte Carlo and Quasi-Monte Carlo Methods in Scientific Computing},
  pages={451--468},
  year={2022},
  organization={Springer}
}

@article{ni2025ergodic,
  title={Ergodic and foliated kernel-differentiation method for linear responses of random systems},
  author={Ni, Angxiu},
  journal={Journal of Nonlinear Science},
  volume={35},
  number={5},
  pages={90},
  year={2025},
  publisher={Springer}
}

@article{wang2016efficiency,
  title={{Efficiency of the Girsanov transformation approach for parametric sensitivity analysis of stochastic chemical kinetics}},
  author={Wang, Ting and Rathinam, Muruhan},
  journal={SIAM/ASA Journal on Uncertainty Quantification},
  volume={4},
  number={1},
  pages={1288--1322},
  year={2016},
  publisher={SIAM}
}

@article{bal2011importance,
  title={{Importance sampling and adjoint hybrid methods in Monte Carlo transport with reflecting boundaries}},
  author={Bal, Guillaume and Langmore, Ian},
  journal={arXiv preprint arXiv:1104.2550},
  year={2011}
}

@phdthesis{hoogenboom1977adjoint,
  title={{Adjoint Monte Carlo methods in neutron transport calculations}},
  author={Hoogenboom, Jan Eduard},
  year={1977},
  school={Delft University Press}
}

@book{hinze2008optimization,
  title={{Optimization with PDE constraints}},
  author={Hinze, Michael and Pinnau, Ren{\'e} and Ulbrich, Michael and Ulbrich, Stefan},
  volume={23},
  year={2008},
  publisher={Springer Science \& Business Media}
}

@article{babovsky1989convergence,
  title={A convergence proof for {N}anbu’s simulation method for the full {Boltzmann} equation},
  author={Babovsky, Hans and Illner, Reinhard},
  journal={{SIAM Journal on Numerical Analysis}},
  volume={26},
  number={1},
  pages={45--65},
  year={1989},
  publisher={SIAM}
}

@article{babovsky1986simulation,
  title={On a simulation scheme for the {Boltzmann} equation},
  author={Babovsky, Hans and Neunzert, H},
  journal={{Mathematical Methods in the Applied Sciences}},
  volume={8},
  number={1},
  pages={223--233},
  year={1986},
  publisher={Wiley Online Library}
}

@article{nanbu1980direct,
  title={Direct simulation scheme derived from the {Boltzmann} equation. {I.} Monocomponent gases},
  author={Nanbu, Kenichi},
  journal={Journal of the Physical Society of Japan},
  volume={49},
  number={5},
  pages={2042--2049},
  year={1980},
  publisher={The Physical Society of Japan}
}

@book{cercignani2000rarefied,
  title={Rarefied gas dynamics: from basic concepts to actual calculations},
  author={Cercignani, Carlo},
  volume={21},
  year={2000},
  publisher={{Cambridge University Press}}
}

@article{bird1994molecular,
  title={Molecular gas dynamics and the direct simulation of gas flows},
  author={Bird, Graeme A},
  journal={Molecular gas dynamics and the direct simulation of gas flows},
  year={1994}
}

@article{bird1970direct,
  title={Direct simulation and the {B}oltzmann equation},
  author={Bird, GA},
  journal={The Physics of Fluids},
  volume={13},
  number={11},
  pages={2676--2681},
  year={1970},
  publisher={American Institute of Physics}
}

@inproceedings{pareschi2001introduction,
  title={{An introduction to Monte Carlo method for the Boltzmann equation}},
  author={Pareschi, Lorenzo and Russo, Giovanni},
  booktitle={ESAIM: Proceedings},
  volume={10},
  pages={35--75},
  year={2001},
  organization={EDP Sciences}
}

@article{cercignani1988boltzmann,
  title={The {Boltzmann} equation and its applications. 1988},
  author={Cercignani, Carlo},
  journal={Applied Mathematical Sciences},
  year={1988}
}

@article{li2022monte,
  title={{Monte Carlo Gradient in Optimization Constrained by Radiative Transport Equation}},
  author={Li, Qin and Wang, Li and Yang, Yunan},
  journal={SIAM Journal on Numerical Analysis},
  volume={61},
  number={6},
  pages={2744--2774},
  year={2023},
  publisher={SIAM}
}

@article{yang2023adjoint,
  title={{Adjoint DSMC for nonlinear spatially-homogeneous Boltzmann equation with a general collision model}},
  author={Yang, Yunan and Silantyev, Denis and Caflisch, Russel},
  journal={Journal of Computational Physics},
  pages={112247},
  year={2023},
  publisher={Elsevier}
}

\appendix
\section{Formula of $\partial_x \log p_l(x,v)$ and  $\partial_v \log p_l(x,v)$ in Section~\ref{subsec:link_appendx_A}}\label{app:log}
This appendix provides detailed calculations of the score-function terms introduced in Section~\ref{subsec:link_appendx_A} for the adjoint system associated with the thermal boundary condition.

For the 1D domain $\Omega=[L,R]$, define
\[
t_i^L = \frac{L-x_{k,i}}{\mathcal P v'_{k,i}},
\qquad
t_i^R = \frac{R-x_{k,i}}{\mathcal P v'_{k,i}} .
\]
If $\mathcal P v'_{k,i}<0$, then
\[
(p_1,p_2,p_3) =
\bigl(1-F_{\Delta t,\varepsilon}(t_i^L),\;
F_{\Delta t,\varepsilon}(t_i^L),\; 0\bigr),
\]
while if $\mathcal P v'_{k,i}>0$, then
\[
(p_1,p_2,p_3) =
\bigl(0,\;
F_{\Delta t,\varepsilon}(t_i^R),\; 1-F_{\Delta t,\varepsilon}(t_i^R)\bigr),
\]
where $F_{\Delta t,\varepsilon}$ denotes the cumulative distribution function of $\mathcal N(\Delta t,\varepsilon^2)$. We also denote by $f_{\Delta t,\varepsilon}$ its probability density function. Here, $p_l$ is a short-hand notation for $p_l(x,v)$, $l= 1,2,3$.

Define the hazard and reverse-hazard ratios
\[
\lambda(t):=\frac{f_{\Delta t,\varepsilon}(t)}{F_{\Delta t,\varepsilon}(t)},\qquad
\bar\lambda(t):=\frac{f_{\Delta t,\varepsilon}(t)}{1-F_{\Delta t,\varepsilon}(t)}.
\]
For $\mathcal P v<0$ set $t^L=\frac{L-x}{\mathcal P v}$, and for $\mathcal P v>0$ set $t^R=\frac{R-x}{\mathcal P v}$. Then
\begin{align*}
\nabla_x \log p_1(x,v) &= \mathds{1}_{\mathcal P v<0}\;\frac{1}{\mathcal P v}\,\bar\lambda(t^L),\\
\nabla_x \log p_2(x,v) &= -\mathds{1}_{\mathcal P v<0}\;\frac{1}{\mathcal P v}\,\lambda(t^L) -\mathds{1}_{\mathcal P v>0}\;\frac{1}{\mathcal P v}\,\lambda(t^R),\\
\nabla_x \log p_3(x,v) &= \mathds{1}_{\mathcal P v>0}\;\frac{1}{\mathcal P v}\,\bar\lambda(t^R),
\end{align*}
and
\begin{align*}
\nabla_v \log p_1(x,v) &= \mathds{1}_{\mathcal P v<0}\;\frac{t^L}{\mathcal P v}\,\bar\lambda(t^L)\,\mathcal P^*,\\
\nabla_v \log p_2(x,v) &= -\mathds{1}_{\mathcal P v<0}\;\frac{t^L}{\mathcal P v}\,\lambda(t^L)\,\mathcal P^* -\mathds{1}_{\mathcal P v>0}\;\frac{t^R}{\mathcal P v}\,\lambda(t^R)\,\mathcal P^*,\\
\nabla_v \log p_3(x,v) &= \mathds{1}_{\mathcal P v>0}\;\frac{t^R}{\mathcal P v}\,\bar\lambda(t^R)\,\mathcal P^* .
\end{align*}

\end{document}